\documentclass[11pt,reqno]{article}

\usepackage{amsmath,amsthm,amsfonts,amssymb,bm,enumerate,color,tikz,hyperref}
\usetikzlibrary{matrix,arrows,positioning,calc}
\usepackage{tikzpagenodes}
\usepackage[margin=1.1in]{geometry}
\usepackage{etex}
\usepackage{cite}

\usepackage{pifont}
\usepackage{makecell} 
\usepackage{multirow}
\usepackage{latexsym}
\usepackage{mathtools}

\usepackage{ragged2e}
\usepackage{graphics,graphicx,psfrag}
\usepackage{amscd}
\usepackage{epsfig}
\usepackage[all]{xy}
\usepackage{float}
\usepackage{enumitem}
\usetikzlibrary{decorations.pathmorphing}
\usepackage{setspace}
\newtheoremstyle{break}
  {0.5cm}%
  {0.5cm}%
  {\itshape}%
  {}%
  {\bfseries}%
  {\vspace{3pt}}%
  {\newline}%
  {}%
\theoremstyle{break}

\newtheorem{definition}{Definition}[subsection]
\newtheorem{theorem}[definition]{Theorem}
\newtheorem{lemma}[definition]{Lemma}
\newtheorem{corollary}[definition]{Corollary}
\newtheorem{proposition}[definition]{Proposition}
\newtheorem{remark}[definition]{Remark}
\newtheorem{example}[definition]{Example}

\newenvironment{myproof}[1][\proofname]{%
  \begin{proof}[\textbf{\textit{Proof}}]$ $\par\vspace{-2pt}
}{%
  \end{proof}
}

\numberwithin{equation}{section}

\newcommand{\Hom}{\mathsf{Hom}}
\newcommand{\Ext}{\mathsf{Ext}}
\newcommand{\Tor}{\mathsf{Tor}}
\newcommand{\Flat}{\mathsf{Flat}}
\newcommand{\DProj}{\mathsf{DProj}}
\newcommand{\DInj}{\mathsf{DInj}}
\newcommand{\pd}{\mathsf{pd}}
\newcommand{\id}{\mathsf{id}}
\newcommand{\too}{\longrightarrow}

\newcommand{\A}{\mathcal{A}}
\newcommand{\B}{\mathcal{B}}

\newcommand{\C}{\mathcal{C}}
\newcommand{\F}{\mathcal{F}}
\newcommand{\G}{\mathcal{G}}

\newcommand{\X}{\mathcal{X}}
\newcommand{\Y}{\mathcal{Y}}
\newcommand{\Sb}{\mathcal{S}}

\newcommand{\Gpd}{\mathsf{Gpd}}

\newcommand{\Proj}{\mathsf{Proj}}
\newcommand{\Inj}{\mathsf{Inj}}

\newcommand{\Gid}{\mathsf{Gid}}
\newcommand{\resdim}{\mathsf{resdim}}
\newcommand{\Thick}{\mathsf{Thick}}

\newcommand{\coresdim}{\mathsf{coresdim}}
\newcommand{\add}{\mathsf{add}}

\newcommand{\Ker}{\mathrm{Ker}}

\newcommand{\Coker}{{\rm CoKer}}
\newcommand{\op}{{}^{\rm op}}
\newcommand{\GProj}{\mathsf{GProj}}
\newcommand{\GInj}{\mathsf{GInj}}

\newcommand{\uno}{\mathsf{(1)}}
\newcommand{\dos}{\mathsf{(2)}}
\newcommand{\tres}{\mathsf{(3)}}
\newcommand{\cuatro}{\mathsf{(4)}}
\newcommand{\cinco}{\mathsf{(5)}}
\newcommand{\ai}{\mathsf{(a)}}
\newcommand{\bi}{\mathsf{(b)}}
\newcommand{\ci}{\mathsf{(c)}}
\newcommand{\di}{\mathsf{(d)}}

\newcommand{\iroman}{\mathsf{(i)}}
\newcommand{\iiroman}{\mathsf{(ii)}}
\newcommand{\iiiroman}{\mathsf{(iii)}}
\newcommand{\scpuno}{\mathsf{(scp1)}}
\newcommand{\scpdos}{\mathsf{(scp2)}}
\newcommand{\scptres}{\mathsf{(scp3)}}
\newcommand{\scpcuatro}{\mathsf{(scp4)}}
\newcommand{\scpcinco}{\mathsf{(scp5)}}

\hypersetup{
    bookmarks=true,         % show bookmarks bar?
    unicode=false,          % non-Latin characters in Acrobat?s bookmarks
    pdftoolbar=true,        % show Acrobat?s toolbar?
    pdfmenubar=true,        % show Acrobat?s menu?
    pdffitwindow=false,     % window fit to page when opened
    pdfstartview={FitH},    % fits the width of the page to the window
    pdftitle={My title},    % title
    pdfauthor={Author},     % author
    pdfsubject={Subject},   % subject of the document
    pdfcreator={Creator},   % creator of the document
    pdfproducer={Producer}, % producer of the document
    pdfkeywords={keyword1} {key2} {key3}, % list of keywords
    pdfnewwindow=true,      % links in new window
    colorlinks=true,       % false: boxed links; true: colored links
    linkcolor=black,          % color of internal links (change box color with linkbordercolor)
    citecolor=magenta,        % color of links to bibliography
    filecolor=cyan,      % color of file links
    urlcolor=violet           % color of external links
}

\usepackage{fancyhdr}

\fancypagestyle{plain}{%
\fancyhf{}% Clear all headers/footers
\fancyhead[L]{Becerril-Mendoza-P\'erez-Santiago}\fancyhead[C]{}\fancyhead[R]{Frobenius pairs in abelian categories}
\fancyfoot[L]{}\fancyfoot[C]{}\fancyfoot[R]{Page \thepage}
\pagestyle{fancy}
\thispagestyle{plain}
}

\usepackage{tocloft}

\usepackage{titlesec}

\titleformat{\subsubsection}[runin]
  {\normalfont\bfseries}{\thesubsection}{1em}{}

\begin{document}

{\huge Frobenius pairs in abelian categories:} 

\vspace{0cm}

{\sc correspondences with cotorsion pairs, exact model categories, and Auslander-Buchweitz contexts}

\vspace{0cm}

{\footnotesize
\vskip3mm \noindent Victor Becerril \\ 
Instituto de Matem\'aticas \\
Universidad Nacional Aut\'onoma de M\'exico \\ 
Circuito Exterior, Ciudad Universitaria \\
C.P. 04510, M\'exico, D.F. MEXICO \\ 
{\tt mathvick06@gmail.com}

\vskip3mm \noindent Octavio Mendoza Hern\'andez \\ 
Instituto de Matem\'aticas \\ 
Universidad Nacional Aut\'onoma de M\'exico \\ 
Circuito Exterior, Ciudad Universitaria \\
C.P. 04510, M\'exico, D.F. MEXICO \\ 
{\tt omendoza@matem.unam.mx}

\vskip3mm \noindent Marco A. P\'erez \\ 
Instituto de Matem\'aticas \\ 
Universidad Nacional Aut\'onoma de M\'exico \\ 
Circuito Exterior, Ciudad Universitaria \\
C.P. 04510, M\'exico, D.F. MEXICO \\ 
{\tt maperez@im.unam.mx}

\vskip3mm \noindent Valente Santiago\\
Departamento de Matem\'aticas. Facultad de Ciencias\\
Universidad Nacional Aut\'onoma de M\'exico \\
Circuito Exterior, Ciudad Universitaria\\
M\'exico D.F. 04510, M\'EXICO.\\
{\tt valente.santiago.v@gmail.com}

}

\vfill

{\small The authors thank Project PAPIIT-Universidad Nacional Aut\'onoma de M\'exico IN102914}

\

\hfill \today

\thispagestyle{empty}

\newpage

%%%%%%%%%%%%%%%%%%%%%%%%%%%%%%%%%%%%%
%%%%%%%%%%%%%%%%%%%%%%%%%%%%%%%%%%%%% 

~\

\vspace{6cm}

\thispagestyle{empty}

\begin{abstract} 
\noindent In this work, we revisit Auslander-Buchweitz Approximation Theory and find some relations with cotorsion pairs and model category structures. From the notions of relatives generators and cogenerators in Approximation Theory, we introduce the concept of left Frobenius pairs $(\X,\omega)$ in an abelian category $\C$. We show how to construct from $(\X,\omega)$ a projective exact model structure on $\X^\wedge$, as a result of Hovey-Gillespie Correspondence applied to two compatible and complete cotorsion pairs in $\X^\wedge$. These pairs can be regarded as examples of what we call cotorsion pairs relative to a thick sub-category of $\C$. We establish some correspondences between Frobenius pairs, relative cotorsion pairs, exact model structures and Auslander-Buchweitz contexts. Finally, some applications of these results are given in the context of Gorenstein homological algebra by generalizing some existing model structures on the categories of modules over Gorenstein and Ding-Chen rings, and by encoding the stable module category of a ring as a certain homotopy category. We also present some connections with perfect cotorsion pairs, covering classes, and cotilting modules.
\end{abstract}

\newpage

%%%%%%%%%%%%%%%%%%%%%%%%%%%%%%%%%%%%%
%%%%%%%%%%%%%%%%%%%%%%%%%%%%%%%%%%%%% 

\tableofcontents

\listoffigures

\newpage

%%%%%%%%%%%%%%%%%%%%%%%%%%%%%%%%%%%%%
%%%%%%%%%%%%%%%%%%%%%%%%%%%%%%%%%%%%% 

\section{Introduction}

The study of homological dimensions, which are obtained by replacing the projective or injective modules by certain sub-categories, was initiated by Maurice Auslander and Ragnar O. Buchweitz in their seminal paper \emph{\cite{AuB}}, which was the starting point for what is now called \textit{relative homological algebra}. Of course, the existence of approximations is the prerequisite for computing relative dimensions. In recent years, a powerful machinery for producing approximations was developed, see \emph{\cite{ET, EJ, GT}} for instance. So it is not surprising that \textit{Auslander-Buchweitz Approximation Theory} (to which we refer as ``\textit{AB Theory}'' for short) provides a good setting for investigating relative Gorenstein projective objects in abelian categories. 
 
The main purpose of this paper is to use AB Theory in order to develop, in the general setting provided by an abelian category $\C$, the \textit{theory of left and right Frobenius pairs}, and give some applications in \textit{Gorenstein Homological Algebra}, \textit{Model Category Theory}, \textit{Homotopy Theory}, and \textit{Cotilting Theory}.   

This paper is organized as follows. We begin recalling in \emph{Section~\ref{sec:ABtheory}} some concepts and results from AB Theory, along with basic notions from relative homological algebra such as \textit{resolutions} and \textit{homological dimensions}. We also present the notion of (left and right) Frobenius pairs (see \emph{Definition~\ref{def:Frobenius_pair}}), which will constitute the main subject studied in this work. In \emph{Section~\ref{sec:relative_pairs}}, we recall the concept of cotorsion pairs in exact categories. In the particular case where $\mathcal{S} \subseteq \C$ is a thick sub-category, $\Sb$ can be regarded as an exact category, and a \textit{complete} cotorsion pair in $\mathcal{S}$ is what we call an \textit{$\Sb$-cotorsion pair} (see \emph{Definition~\ref{def:left_cotorsion_pair}}). We later provide in \emph{Proposition~\ref{CP1}} an alternative description of $\Sb$-cotorsion pairs, and use it to induce relative cotorsion pairs from the notions of \textit{generator} and \textit{cogenerator} in AB Theory. Motivated by the interplay between cotorsion pairs and model categories, we show how to obtain from a \textit{strong} left Frobenius pair $(\X,\omega)$ two \textit{compatible} and \textit{complete} cotorsion pairs in the sub-category $\Sb := \X^\wedge$ of objects in $\C$ with finite \textit{resolution dimension} with respect to $\X$ (which turns out to be thick). We then apply in \emph{Section~\ref{sec:model_structures}} a result known as \textit{Hovey-Gillespie Correspondence} (see \emph{Theorem~\ref{theo:Hovey-Gillespie}}) to obtain an \textit{exact} model structure on $\X^\wedge$, which we call the \textit{projective Auslander-Buchweitz model structure}, where $\X$, $\X^\wedge$ and $\omega^\wedge$ are the classes of cofibrant, fibrant, and trivial objects, respectively. This model structure can be applied to particular choices of $\X$ in order to encode some known results and to present new ones. The most remarkable applications of this point will be generalizations of \textit{abelian} model structures in Gorenstein and Ding-Chen Homological Algebra. Moreover, we see in \emph{Proposition~\ref{prop:homotopy_AB_structures}} that the homotopy category of this model structure represents, in some sense, a generalization of the stable module category $\mathsf{Stmod}(R)$ of a ring $R$. Finally, in \emph{Section~\ref{sec:AB_contexts}} we recall in \emph{Definition~\ref{def:AB-context}} the notion of \textit{Auslander-Buchweitz contexts}, and present some (one-to-one) correspondences with Frobenius pairs, relative cotorsion pairs, and exact model structures (see \emph{Theorems~\ref{theo:correspondence_with_AB-context}}, \emph{\ref{theo:correspondence_with_G}}, \emph{\ref{theo:correspondence_G_P_one_to_one}}, \emph{\ref{theo:correspondence_Hovey_pairs_vs_Frobenius_pairs}}, and \emph{\ref{theo:Frobenius_pairs_vs_model_structures}}). As an application of these correspondences, we will present in a categorical context, an important theorem by M. Auslander and I. Reiten, which establishes a bijective correspondence between basic cotilting modules in $\mathsf{mod}(\Lambda)$ (the category of finitely generated left modules over an Artin algebra $\Lambda$), resolving pre-covering sub-categories $\F \subseteq \mathsf{mod}(\Lambda)$ such that $\F^\wedge = \mathsf{mod}(\Lambda)$, and coresolving pre-enveloping sub-categories $\G \subseteq \mathsf{mod}(\Lambda)$ with finite injective dimension. 

The results presented in this paper have their corresponding dual statements (some of which are self-dual). For the sake of simplicity, we do not state dual statements. For pedagogical reasons, we only make an exception for those results whose dual versions are also used throughout the paper.    

%%%%%%%%%%%%%%%%%%%%%%%%%%%%%%%%%%%%%
%%%%%%%%%%%%%%%%%%%%%%%%%%%%%%%%%%%%% 
 
\section{Auslander-Buchweitz Approximation Theory}\label{sec:ABtheory}

We start this section by collecting all the background material that will be necessary in the sequel. First, we introduce some general notation. Next, we recall the notions of relative projective dimension and resolution dimension of a given class of objects in an abelian category $\C$. Finally, we also recall  definitions and basic properties  we need from AB Theory. In all that follows, we are taking as a main reference the paper \emph{\cite{AuB}}.

We remark that M. Auslander and R.-O. Buchweitz in some of the results in \emph{\cite{AuB}} worked with a resolving and additively closed sub-category $\X \subseteq \C$, which is also closed under direct summands (in $\C$). In a very carefully revision of the proofs of those results, we can see that some of the properties assumed for $\X$ are not used. In order to give nice applications of AB Theory to Gorenstein Homological Algebra, we give a review by putting in each statement the minimum needed hypothesis.

%%%%%%%%%%%%%%%%%%%%%%%%%%%%%%%%%%
%%%%%%%%%%%%%%%%%%%%%%%%%%%%%%%%%%

\subsection{Notation} 

Throughout this paper, $\C$ will be an abelian category and $\X \subseteq \C$ a class of objects of $\C$ (which can be also though as a full sub-category of $\C$). The notation $M \in \C$ means that $M$ is an object of $\C$. We denote by $\pd(M)$ the \textbf{projective dimension} of $M$. Similarly, $\id(M)$ denotes the \textbf{injective dimension} of $M$. For any non-negative integer $n$, we set 
\[
\Proj_n(\C) := \{M \in \C\;:\;\pd(M) \leq n\}.
\] 
In particular, $\Proj(\C):=\Proj_0(\C)$ is the class of projective objects in $\C$. The classes  $\Inj_n(\C)$  and $\Inj(\C)$ are defined dually. 

Let $\X$ and $\mathcal{Y}$ be two classes of objects of $\C$, and $M$ and $N$ be objects in $\C$. We set the following notation for each non-negative integer $i \geq 0$: 
\begin{itemize}[itemsep=2pt,topsep=0pt]
\item $\Ext^i_{\C}(\X, N) = 0$ if $\Ext^i_{\C}(X,N) = 0$ for every $X \in \X$.

\item $\Ext^i_{\C}(M, \mathcal{Y}) = 0$ if $\Ext^i_{\C}(M, Y) = 0$ for every $Y \in \mathcal{Y}$.

\item $\Ext^i_{\C}(\X, \mathcal{Y}) = 0$ if $\Ext^i_{\C}(X,Y) = 0$ for every $X \in \X$ and $Y \in \mathcal{Y}$.
\end{itemize}

We denote by $\add\,(\X)$ the class of all objects isomorphic to direct summands of finite direct sums of objects in $\mathcal{X}$. Moreover, for each positive integer $i > 0$, we denote 
\[
\X^{\perp_i} := \{ N \in \C\;:\;\Ext^i_\C(\X, N) = 0 \}\quad {\rm and}\quad  \X^\perp := \bigcap_{i>0}\,\X^{\perp_i}.
\] 
Dually, we have the classes ${}^{\perp_i}\X$ and ${}^{\perp}\X.$

%%%%%%%%%%%%%%%%%%%%%%%%%%%%%%%%%%
%%%%%%%%%%%%%%%%%%%%%%%%%%%%%%%%%%

\subsection{Resolving and coresolving classes} 

It is said that $\X$ is a \textbf{pre-resolving} class if it is closed under extensions and kernels of epimorphisms in $\X$. A pre-resolving class is said to be  \textbf{resolving} if it contains $\Proj(\C)$. If the dual properties hold true, then we get \textbf{pre-coresolving} and \textbf{coresolving} classes. A \textbf{left thick} (respectively, \textbf{right thick}) class is a pre-resolving (respectively, pre-coresolving) class which is closed under direct summands in $\C$. A class is \textbf{thick} if it is both a right thick and left thick class. We denote by $\Thick\,(\X)$, $\Thick^{-}(\X)$ and $\Thick^{+}(\X)$ the smallest thick, left thick, and right thick full sub-categories of $\C$, respectively, containing the class $\X$. A \textbf{left saturated} (respectively, \textbf{right saturated}) class is a resolving (respectively, coresolving) class which is closed under direct summands in $\C$. A \textbf{saturated} class is  both a right saturated and left saturated class. For example, $\Proj(\C)$  and ${}^{\perp}\X$ are left saturated classes of $\C$, while $\Inj(\C)$ and $\X^\perp$ are right saturated classes of $\C$. In the case $\C$ is the category $\mathsf{Mod}(R)$ (or $\mathsf{Mod}(R^{\rm op})$) of left (respectively, right) $R$-modules, for the sake of simplicity, we will denote by $\Proj(R)$ and $\Proj(R^{\rm op})$ the classes of projective left and right $R$-modules, respectively.

\begin{remark} 
The concepts presented in the previous paragraph have their analogous in exact categories. 
\end{remark}

%%%%%%%%%%%%%%%%%%%%%%%%%%%%%%%%%%
%%%%%%%%%%%%%%%%%%%%%%%%%%%%%%%%%%

\subsection{Relative homological dimensions} 

Given a class $\X \subseteq \C$ and $M \in \C$, the \textbf{relative projective dimension} of $M$ with respect to $\X$ is defined as 
\[
\pd_{\X}\,(M) := \mathrm{min}\,\{n \geq 0\,:\,\Ext_\C^j(M,\X) = 0  \text{ for every } j > n\}.
\] 
Dually, we denote by $\id_{\X}\,(M)$ the  {\bf{relative injective dimension}} of $M$ with respect to $\X.$ Furthermore, for any class $\Y \subseteq \C$, we set 
\[
\pd_\X\,(\Y) := \mathrm{sup}\,\{\pd_\X\,(Y)\;:\; Y\in\Y\}\text{ and }\id_\X\,(\Y):=\mathrm{sup}\,\{\id_\X\,(Y)\;:\; Y\in\Y\}.
\]
It can be seen that $\pd_\X\,(\Y)=\id_\Y\,(\X).$ If $\X=\C$, we just write $\pd\,(\Y)$ and $\id\,(\Y)$.

%%%%%%%%%%%%%%%%%%%%%%%%%%%%%%%%%%
%%%%%%%%%%%%%%%%%%%%%%%%%%%%%%%%%%

\subsection{Resolution and coresolution dimensions.}

Let $M \in \C$ and $\X$ be a class of objects in $\C$. The \textbf{$\X$-resolution dimension} of $M$, denoted $\resdim_{\X}(M)$, is the smallest non-negative integer $n$ such that there is an exact sequence
\begin{align}\label{eqn:X-resolution}
0 & \to X_n \to X_{n-1} \to \cdots \to X_1 \to X_0 \to M \to 0
\end{align}
with $X_i \in \X$ for every $0 \leq i \leq n$. The sequence \emph{(\ref{eqn:X-resolution})} is said to be a \textbf{finite $\X$-resolution} of $M$. If such $n$ does not exist, we set $\resdim_{\X}(M) := \infty$\footnote{We assume $\resdim_{\X}(M) = \infty$ if either $M$ has only an $\X$-resolution of infinite length, or if $M$ has no $\X$-resolution (finite or infinite). In practice, this assumption may happen to be troublesome, but in this paper we mainly consider finite $\X$-resolutions.}. Also, we denote by $\mathcal{X}^{\wedge}$ the class of objects in $\C$ having a finite $\mathcal{X}$-resolution (or equivalently, having finite $\X$-resolution dimension).  

Dually, we have the \textbf{$\X$-coresolution dimension} of $M$, denoted $\coresdim_\X\,(M)$, and the class $\mathcal{X}^{\vee}$ of objects having a finite $\mathcal{X}$-coresolution. 

Given a class $\Y \subseteq \C$, we set  
\[
\resdim_\X\,(\Y) := \mathrm{sup}\,\{\resdim_\X\,(Y)\;:\; Y \in \Y\},
\]
and $\coresdim_\X\,(\Y)$ is defined dually.

%%%%%%%%%%%%%%%%%%%%%%%%%%%%%%%%%%%
%%%%%%%%%%%%%%%%%%%%%%%%%%%%%%%%%%%

\subsection{Approximations.}  

Let  $\X$ be a class of objects in $\C$. A morphism $f \colon X \to M$ is said to be an \textbf{$\X$-pre-cover} of $M$ if:
\begin{itemize}[itemsep=2pt,topsep=0pt]
\item[$\ai$] $X \in \X$, and 

\item[$\bi$] For every $X' \in \X$, the induced map 
\[
\Hom_\C(X', f) \colon \Hom_\C(X', X) \to \Hom_\C(X', M)
\]
is surjective; or equivalently, for every $f' \colon X' \to M$ with $X' \in \X$, there exists a morphism $h \colon X' \to X$ such that the following triangle commutes:
\begin{figure}[H]
\centering
\begin{tikzpicture}[description/.style={fill=white,inner sep=2pt}]
\matrix (m) [matrix of math nodes, row sep=3em, column sep=3em, text height=1.5ex, text depth=0.25ex]
{ X & M \\ X' \\ };
\path[->]
(m-1-1) edge node[above] {\footnotesize$f$} (m-1-2)
(m-2-1) edge node[below,sloped] {\footnotesize$f'$} (m-1-2);
\path[dotted,->]
(m-2-1) edge node[description] {\footnotesize$\exists \mbox{ } h$} (m-1-1);
\end{tikzpicture}
\caption{Precovers.}
\label{fig:precover}
\end{figure}
\end{itemize}
If in addition, in the case $X' = X$ and $f' = f$ the diagram in \emph{Figure~\ref{fig:precover}} can only be completed by automorphisms of $X$, then the $\X$-pre-cover $f$ is called an \textbf{$\X$-cover}. Furthermore, an $\X$-pre-cover $f \colon X \to M$ is \textbf{special} if $\Coker\,(f) = 0$ and $\Ker\,(f) \in \X^{\perp_1}$. The class $\X$ is said to be \textbf{pre-covering} if every object of $\C$ has an $\X$-pre-cover. Similarly, we can define \textbf{covering} and \textbf{special pre-covering} classes in $\C$, that is, $\X$ is covering (resp., special pre-covering) if every object in $\C$ has an $\X$-cover (resp., a special $\X$-pre-cover). Dually, we have the notions of \textbf{$\X$-envelopes}, \textbf{$\X$-pre-envelopes}, and \textbf{special $\X$-pre-envelopes} in $\C$, along with the corresponding notions of \textbf{enveloping}, \textbf{pre-enveloping}, and \textbf{special pre-enveloping} classes.

\begin{remark}\label{rem:approx} 
These notions of approximations are also valid for exact categories. 
\end{remark}

%%%%%%%%%%%%%%%%%%%%%%%%%%%%%%%%%%%%
%%%%%%%%%%%%%%%%%%%%%%%%%%%%%%%%%%%%

\subsection{Fundamental results in Auslander-Buchweitz Approximation Theory} 

Having the terminology and notation we have presented so far, we are ready to recall the necessary background from AB Theory. For a better understanding of the results below, and as a motivation, we present examples which come from the domain of Gorenstein homological algebra. Much of the properties we mention for Gorenstein-projective and Gorenstein-injective modules were already proven by H. Holm in his paper \emph{\cite{Holm}}, without using AB Theory. So in the following lines, we will appreciate the wide scope of the results obtained by M. Auslander and R.-O. Buchweitz in \emph{\cite{AuB}}. These examples will also be useful to motivate the theory presented in the next sections. 

Let $(\X,\omega)$ be a pair of classes of objects in $\C$. It is said that $\omega$ is \textbf{$\X$-injective} if $\id_\X(\omega) = 0$. We say that $\omega$ is a \textbf{relative cogenerator} in $\X$ if $\omega \subseteq \X$ and for any $X \in \X$ there is an exact sequence 
\[
0 \to X \to W \to X' \to 0,
\] 
with $W \in \omega$ and $X' \in \X$. \newpage

Dually, we have the notions of \textbf{$\X$-projective} and a \textbf{relative generator} in $\X$.

\begin{example}\label{ex:Gorenstein_projective} 
Recall that a left $R$-module $M$ is said to be \textbf{Gorenstein-projective} (or \textbf{G-projective} for short) if there exists an exact sequence of projective modules
\[
\bm{P} = \cdots \to P_1 \to P_0 \to P^0 \to P^1 \to \cdots
\] 
with $M = \Ker(P^0 \to P^1)$, such that $\Hom_R(\bm{P}, P)$ is an exact chain complex of abelian groups, for every $P \in \Proj(R)$. We denote by $\GProj(R)$ the class of G-projective left $R$-modules. \textbf{Gorenstein-injective} modules are defined dually, and the class of such modules will be denoted by $\GInj(R)$. 

Note that every kernel of $\bm{P}$ is a G-projective module. It follows by this fact that $\Proj(R)$ is a relative cogenerator and a relative generator in $\GProj(R)$. Moreover, the exactness of the complex $\Hom_R(\bm{P},P)$ for every $P \in \Proj(R)$ implies that $\id_{\GProj(R)}(\Proj(R)) = 0$, that is, that $\Proj(R)$ is a $\GProj(R)$-injective relative cogenerator in $\GProj(R)$. On the other hand, it is clear that $\Proj(R)$ is $\GProj(R)$-projective.   

Dually, $\Inj(R)$ is both a $\GInj(R)$-projective relative generator and a $\GInj(R)$-injective relative cogenerator in $\GInj(R)$.
\end{example}

The notions of relative generators and cogenerators provide the setting to define a sort of Frobenius category in a relative sense. Let us be more specific about this in the following lines.

\begin{definition}\label{def:Frobenius_pair}
\sloppypar{Let $(\X,\omega)$ be a pair of classes of objects in $\C$. We say that $(\X,\omega)$ is a \textbf{left Frobenius pair} in $\C$ if $\X = \Thick^{-}(\X)$, $\omega$ is an $\X$-injective relative cogenerator in $\X$, and $\omega$ is closed under direct summands in $\C$. We say that a left Frobenius pair $(\X,\omega)$ is \textbf{strong} if $\omega$ is also an $\X$-projective relative generator \mbox{in $\X$.} }

Dually, we say that a pair $(\nu,\Y)$ of classes of objects in $\C$ is a \textbf{right Frobenius pair} in $\C$ if $\Y = \Thick^{+}(\Y)$, $\nu$ is an $\Y$-projective relative generator in $\Y$, and $\nu$ is closed under direct summands in $\C$. If in addition $\nu$ is also an $\Y$-injective relative cogenerator in $\Y$, we say that $(\nu,\Y)$ is a \textbf{strong right Frobenius pair}. 
\end{definition}

\begin{example}\label{ex:G-projective_strong_Frobenius_pair} 
In \emph{\cite[Theoren 2.5]{Holm}}, it is proven that $\GProj(R)$ is a left saturated class in $\mathsf{Mod}(R)$, and so it is left thick. By the comments in \emph{Example~\ref{ex:Gorenstein_projective}}, we then have that $(\GProj(R),\Proj(R))$ is a strong left Frobenius pair in $\mathsf{Mod}(R)$. Dually, $(\Inj(R), \GInj(R))$ is a strong right Frobenius pair in $\mathsf{Mod}(R)$. 
\end{example}

\begin{remark}\label{rem:Frobenius_category} 
Let $(\X,\omega)$ be a strong left Frobenius pair in $\C$. Then $\X$ is a Frobenius category in the usual sense, that is, it is an exact category with enough projectives and injectives, where the classes of projective and injective objects in $\X$ coincide. In this case, the projective and injective objects are given by $\omega$.
\end{remark}

The proof of the following result can be found in \emph{\cite[Lemma 2.13]{MS}}. \newpage

\begin{lemma}\label{AB1}
Let $\X$ and $\Y$ be classes of objects in $\C$. Then
\[
\pd_{\Y}\,(\X^{\vee}) = \pd_{\Y}\,(\X) \mbox{ \ and \ } \id_{\X}(\Y^\wedge) = \id_{\X}(\Y).
\]
\end{lemma}

\begin{proposition}\label{AB2} 
Let $(\X,\omega)$ be a pair of classes of objects in $\C$, such that $\omega$ is $\X$-injective. Then, the following statements hold true.
\begin{itemize}[itemsep=2pt,topsep=0pt]
\item[$\ai$] $\omega^\wedge$ is $\X$-injective.

\item[$\bi$] If in addition $\omega$ is a relative cogenerator in $\X$ which is closed under direct summands in $\C$, then
\[
\omega = \{ X \in \X \mbox{ : } \id_\X(X) = 0 \} = \X \cap \omega^\wedge.
\]
\end{itemize}

Dually, if $(\nu,\Y)$ is a pair of classes of objects in $\C$ such that $\nu$ is $\Y$-projective, then the following statements hold true:
\begin{itemize}[itemsep=2pt,topsep=0pt]
\item[$\iroman$] $\nu^\vee$ is $\Y$-projective.

\item[$\iiroman$] If in addition $\nu$ is a relative generator in $\Y$ which is closed under direct summands in $\C$, then 
\[
\nu = \{ Y \in \Y \mbox{ : } \pd_{\Y}(Y) = 0 \} = \Y \cap \nu^\vee.
\]
\end{itemize}
\end{proposition}

\begin{myproof} 
Part $\ai$ follows by \emph{Lemma~\ref{AB1}}, and part $\bi$ is a consequence of \emph{\cite[Lemma 3.7]{AuB}}.
\end{myproof}

\begin{example}\label{ex:G-projective_compatible} 
By \emph{Example~\ref{ex:G-projective_strong_Frobenius_pair}}, the pair $(\GProj(R), \Proj(R))$ satisfies the hypothesis in \emph{Proposition~\ref{AB2}}. So part $\ai$ implies that $\id_{\GProj(R)}(\Proj(R)^\wedge) = 0$. Note also that $\Proj(R)^\wedge = \Proj^{< \infty}(R)$, the class of modules with finite projective dimension. It follows that if $M$ is a G-projective module, then $\Ext^i_R(M, W) = 0$ for every $W$ with finite projective dimension, and every $i > 0$. This property is also stated in \emph{\cite[Proposition 2.3]{Holm}}, and its dual is also valid for G-injective modules. On the other hand, part $\bi$ implies that 
\[
\Proj(R) = \GProj(R) \cap \Proj^{< \infty}(R),
\] 
that is, the projective dimension of a G-projective module is either $0$ or infinite. Thus, we have another proof of \emph{\cite[Proposition 10.2.3]{EJ}}, and its corresponding dual statement for G-injective modules:
\[
\Inj(R) = \GInj(R) \cap \Inj^{< \infty}(R).
\]
\end{example}

In the previous result, if in addition we assume that $\omega$ is a relative cogenerator in $\X$ and $\omega$ is closed under direct summands, we can obtain a description for the class $\X \cap \omega^\vee$. This is specified in the following result, whose proof can be found in \emph{\cite[Lemma 4.3]{AuB}}.

\begin{lemma}\label{AB3} 
Let $(\X,\omega)$ be a pair of classes of objects in $\C$, such that $\omega$ is $\X$-injective and a relative cogenerator in $\X$. If $\omega$ is closed under direct summands in $\C$, then
\[
\X \cap \omega^{\vee} = \{X \in \X\mid\id_\X(X) < \infty \}.
\]
Furthermore, we have that $\id_\X(M) = \coresdim_\omega(M)$ for every $M \in \X \cap \omega^{\vee}$. \newpage
\end{lemma}

In the following result (whose proof can be found in \emph{\cite[Theorem 1.1]{AuB}} ), the expression $\resdim_\omega(K) = -1$ just means that $K = 0$.

\begin{theorem}\label{AB4} 
\sloppypar{Let $(\X,\omega)$ be a pair of classes of objects in $\C$, such that $\X$ is closed under extensions, $0 \in \X$ and $\omega$ is a relative cogenerator in $\X$. Then the following statements hold true, for any $C \in \C$ with \mbox{$\resdim_\X(C) = n < \infty$}. }
\begin{itemize}[itemsep=2pt,topsep=0pt]
\item[$\ai$]There exist exact sequences in $\C$
\[
0 \to K\to X\xrightarrow{\varphi} C\to 0,
\]
with $\resdim_\omega(K) = n-1$ and $X \in \X$, and
\[
0 \to C \xrightarrow{\varphi'} H \to X' \to 0,
\]
with $\resdim_\omega(H) \leq n$ and $X' \in \X$.

\item[$\bi$] If $\omega$ is $\X$-injective, then
	\begin{itemize}[itemsep=2pt,topsep=0pt]
	\item[$\iroman$] $\varphi \colon X \to C$ is an $\X$-pre-cover and $K \in \X^{\perp}$.
	
	\item[$\iiroman$] $\varphi' \colon C \to H$ is an $\omega^{\wedge}$-pre-envelope and $X'\in{}^\perp(\omega^{\wedge})$. ~%
	
	\end{itemize}
 \end{itemize}
\end{theorem}

\begin{example}\label{ex:G-projective_ses} 
From now on, we denote by $\Gpd(M)$ and $\Gid(M)$ the G-projective and G-injective dimensions of a left $R$-module $M$, which is defined (see \emph{\cite[Chapter 11]{EJ}} ) as 
\[
\Gpd(M) := \resdim_{\GProj(R)}(M) \mbox{ \ and \ } \Gid(M) := \coresdim_{\GInj(R)}(M).
\] 
If we are given a Gorenstein ring $R$, then it is known that $(\GProj(R), \Proj^{< \infty}(R))$ is a cotorsion pair in $\mathsf{Mod}(R)$ (see \emph{Section~\ref{sec:relative_pairs}} to recall the definition of cotorsion pairs). Moreover, this pair is known to be complete, and one way to see this is by providing a cogenerating set for it (see \emph{\cite[Theorem 8.3]{Hovey}} for instance). This means that for every module $M$ over a Gorenstein ring $R$, there exists a short exact sequence 
\[
0 \to W \to P \to M \to 0
\]
where $P$ is G-projective and $W$ has finite projective dimension. The previous sequence is constructed directly (that is, without using the fact that cotorsion pairs with a cogenerating set are complete) by H. Holm in \emph{\cite[Theorem 2.10]{Holm}}, for every module $M$ with finite G-projective dimension and over an arbitrary associative ring $R$. Moreover, Holm also shows that if $\Gpd(M) = n$, then $\pd(W) = n-1$. The proof of this is not trivial. On the other hand, this is precisely what \emph{Theorem~\ref{AB4}} implies after setting $\X := \GProj(R)$ and $\omega := \Proj(R)$. 

It is worth recalling that if $R$ is a Gorenstein ring, then $\GProj^{< \infty}(R) = \GProj(R)^\wedge$ coincides with the whole category $\mathsf{Mod}(R)$, and so the two short exact sequences described in \emph{Theorem~\ref{AB4}} imply that $(\GProj(R), \Proj^{< \infty}(R))$ is a complete cotorsion pair. We will return to this point in \emph{Sections~\ref{sec:relative_pairs}} and \emph{\ref{sec:model_structures}}.  
\end{example}

\begin{corollary}\label{AB5} 
Let $(\X,\omega)$ be a pair of classes of objects in $\C$ such that $\X$ is closed under extensions and direct summands in $\C$, and let $\omega$ be a relative cogenerator in $\X$. If $\resdim_\X(C)\leq 1$ and $C \in {}^{\perp_1}\omega$, then $C \in \X$. 
\end{corollary}

\begin{myproof}
Since $\resdim_\X(C)\leq 1$, we get from \emph{Theorem~\ref{AB4}} an exact sequence 
\[
\varepsilon:\quad 0\to K\to X\to C\to 0,
\]
with $X \in \X$ and $K \in \omega$. Then, by the fact that $C \in {}^{\perp_1}\omega$, the exact sequence $\varepsilon$ splits and so $C \in \X$.
\end{myproof}

\begin{example} 
Set $\X := \GProj(R)$ and $\omega := \Proj(R)$. Using \emph{Corollary~\ref{AB5}}, we have that if we are given a short exact sequence
\[
0 \to A \to B \to C \to 0
\]
with $A, B \in \GProj(R)$, and $\Ext^1_R(C,P) = 0$ for every $P \in \Proj(R)$, then $C \in \GProj(R)$. This property of $\GProj(R)$ was first proven in \emph{\cite[Corollary 2.11]{Holm}}.
\end{example}

The following result (whose proof can be found in \emph{\cite[Proposition 2.1]{AuB}} ), relates the concepts of relative projective and resolution dimensions for relative injective cogenerators.

\begin{theorem}\label{AB8} 
Let $\X \subseteq \C$ be a class of objects closed under extensions and direct summands in $\C$, and $\omega$ be an $\X$-injective relative cogenerator in $\X$, which is closed under direct summands in $\C$. Then
\[
\pd_{\omega^{\wedge}}(C) = \pd_{\omega}\,(C) = \resdim_{\X}(C), \mbox{ for every } C\in \X^{\wedge}\!.
\]
\end{theorem}

Certain homological dimensions are defined as projective or injective dimensions relative to a certain class of modules (over a ring $R$), such as the \textit{FP-injective dimension}\footnote{``FP'' refers to ``finitely presented''.}. There are others, such as the G-projective dimension, which are defined as a resolution dimension relative to a class of modules. In the former case, the FP-injective dimension cannot be expressed as a coresolution dimension relative to a class, unless we assume $R$ is a coherent ring. This is not the case for the G-projective dimension, as we explain below.

\begin{example}\label{ex:G-projective_dimension}
By \emph{Theorem~\ref{AB8}}, we have that for every left $R$-module $M$ with finite G-projective dimension
\[
\Gpd(M) := \resdim_{\GProj(R)}(M) = \pd_{\Proj(R)}(M) = \pd_{\Proj^{< \infty}(R)}(M).
\]
In other words, we have that the following conditions are equivalent:
\begin{itemize}[itemsep=2pt,topsep=0pt]
\item[$\uno$] $\Gpd(M) \leq n$.

\item[$\dos$] $\Ext^i_R(M,L) = 0$ for every $i > n$ and every $L$ with finite projective dimension.

\item[$\tres$] $\Ext^i_R(M,P) = 0$ for every $i > n$ and every projective left $R$-module $P$. 
\end{itemize}
This was proven by Holm in \emph{\cite[Theorem 2.20]{Holm}}, and probably first by E. E. Enochs and O. M. G. Jenda \emph{\cite[Proposition 11.5.7]{EJ}} in the case where $R$ a Gorenstein ring. 
\end{example}

So far the reader may have already noticed that under certain conditions, the objects in the class $\X^\wedge$ have special properties. This class will turn out to be an exact sub-category of $\C$ in which we will get cotorsion pairs and construct model structures in the following sections. So the importance of $\X^\wedge$ demands to know another description for it, from which one can deduce that $\X^\wedge$ is exact indeed. This description is given in the next theorem, originally proven in \emph{\cite[Propositions 3.4 and 3.5]{AuB}}.

\begin{theorem}\label{AB9} 
Let $\X \subseteq \C$ be a pre-resolving class of objects, and $\omega$ be $\X$-injective and a relative cogenerator in $\X$. Then, the following conditions hold true.
\begin{itemize}[itemsep=2pt,topsep=0pt]
\item[$\ai$] $\X^{\wedge}$ is the smallest pre-resolving and pre-coresolving class in $\C$, containing the class $\X$.

\item[$\bi$] If $\omega$ and $\X$ are closed under direct summands, then $\add\,(\X^\wedge) = \X^\wedge$. 
\end{itemize}
In particular, if $(\X,\omega)$ is a left Frobenius pair, then $\X^\wedge = \Thick\,(\X)$. 
\end{theorem}

\begin{proposition}\label{AB10} 
Let $(\X,\omega)$ be a left Frobenius pair in $\C$. Then, for any $C \in \X^\wedge$ and $n \geq 0$, the following conditions are equivalent.
\begin{itemize}[itemsep=2pt,topsep=0pt]
\item[(a)] $\resdim_\X(C) \leq n$.

\item[(b)] If 
\[
0\to K_n\to X_{n-1}\to\cdots \to X_1\to X_0\to C\to 0
\] 
is an exact sequence, with $X_i\in\X,$ then $K_n\in\X$.
\end{itemize}
\end{proposition}

The previous result, proven in \emph{\cite[Proposition 3.3]{AuB}}, has also applications to Gorenstein Homological Algebra. The particular case where $\X := \GProj(R)$ and $\omega := \Proj(R)$ was first proven in \emph{\cite[Theorem 2.20]{Holm}}, and extends the equivalence mentioned in \emph{Example~\ref{ex:G-projective_dimension}}.

The cotorsion pairs we construct in the following sections, and the model structures associated to them on $\X^\wedge$, involve the classes $\omega$ and $\omega^\wedge$. So we devote the rest of this section to present properties of these classes. We begin with the following description of $\omega^\wedge$, whose proof can be found to \emph{\cite[Proposition 3.6]{AuB}}.

\begin{proposition}\label{AB11} 
Let $\X\subseteq\C$ be a class of objects closed under extensions, and $\omega$ be a class closed under direct summands in $\C$, which is $\X$-injective and a relative cogenerator in $\X$. Then
\[
\omega^\wedge = \X^\perp \cap \X^\wedge.
\]
\end{proposition}

Given a left Frobenius pair $(\X,\omega)$, the class $\omega^\wedge$ is not necessarily thick. The following three results sort of measure how far is $\omega^\wedge$ from being thick.

\begin{proposition}\label{AB12} 
Let $(\X,\omega)$ be a left Frobenius pair in $\C$. Then
\[
\X^{\wedge} \cap {}^{\perp}\omega = \X = \X^{\wedge} \cap {}^{\perp}(\omega^{\wedge}).
\]
\end{proposition}

\begin{myproof} 
By \emph{Proposition~\ref{AB2}} $\ai$, we know that $\X \subseteq {}^\perp\omega$ and $\X \subseteq {}^\perp(\omega^\wedge)$.

We assert that $\X^{\wedge} \cap {}^{\perp}(\omega^{\wedge}) \subseteq \X$. Indeed, let $C \in \X^{\wedge} \cap {}^{\perp}(\omega^{\wedge})$. Then, from \emph{Theorem~\ref{AB4}}, there is an exact sequence 
\[
\varepsilon : \;0\to C\to Y\to X\to 0,
\] 
where $X \in \X \subseteq {}^\perp(\omega^{\wedge})$ and $Y \in  \omega^{\wedge}$. Since $C \in {}^\perp(\omega^{\wedge})$, it follows that $Y \in {}^\perp(\omega^{\wedge})$, and therefore $\Ext^ i_\C(Y, \omega^\wedge) = 0$, for any $i \geq 1$. But now, using \emph{Theorem~\ref{AB8}}, we get 
\[
0 = \pd_{\omega^{\wedge}}(Y) = \resdim_{\X}(Y)
\]
and thus $Y \in \X$. After that, we have the exact sequence $\varepsilon$, with $Y, X \in \X$. Hence $C \in \X$ since $\X$ is pre-resolving.

Finally, the inclusion $\X^{\wedge} \cap {}^{\perp}\omega \subseteq \X$ follows as in the preceding proof. 
\end{myproof}

In order to get a description of the class $\Thick(\omega),$ we introduce the following definition.

\begin{definition} 
For any class $\X \subseteq \C$ of objects in $\C$, we set 
\[
\mathsf{Inj}_\X^{<\infty}(\C) := \{ C \in \C \;|\;\id_\X(C) < \infty \}.
\]
\end{definition}

Note that $\mathsf{Inj}_\X^{<\infty}(\C)$ is a thick sub-category of the abelian category $\C$. 

For every full sub-category $\Y \subseteq \C$, we set $\mathsf{Inj}_{\X}^{< \infty}(\Y) := \mathsf{Inj}_{\X}^{< \infty}(\C) \cap \Y$.

\begin{theorem}\label{AB14} 
Let $\X \subseteq \C$ be a pre-resolving class of objects in $\C$, and $\omega$ be an $\X$-injective relative cogenerator in $\X$. Then, the following statements hold.
\begin{itemize}[itemsep=2pt,topsep=0pt]
\item[$\ai$] $(\omega^\wedge)^\vee = \mathsf{Inj}_\X^{<\infty}(\X^\wedge)$.

\item[$\bi$] If $\X$ and $\omega$ are closed under direct summands in $\C$ (and so $(\X,\omega)$ is a left Frobenius pair in $\C$), then
\[
(\omega^\wedge)^\vee = \Thick\,(\omega).
\]
\end{itemize}
\end{theorem}

\begin{myproof} Part $\ai$ follows by \emph{\cite[Proposition 4.2]{AuB}}. 

We now prove part $\bi$. Let $\X$ and $\omega$ be closed under direct summands in $\C$. Then, by \emph{Theorem~\ref{AB9}}, we know that $\X^\wedge = \Thick\,(\X)$, and thus $(\omega^\wedge)^\vee$ is a thick sub-category in $\C$ (see $\ai$). Assume that $\B$ is a thick sub-category of $\C$ containing $\omega$. Since $\B$ is closed under cokernel of monomorphisms, it follows that $\omega^\wedge \subseteq \B$, and using now that $\B$ is closed under kernel of epimorphisms we get that $(\omega^\wedge)^\vee \subseteq \B$.
\end{myproof}

In the next theorem we obtain another equality involving the class $\mathsf{Inj}^{< \infty}_{\X}(\X^\wedge)$. It is a result due to M. Auslander, R.-O. Buchweitz, and I. Reiten \emph{\cite{AuB,AR}}. The statement given below is a simplification of the one given in \emph{\cite[Theorem 1.12.10]{Hashimoto}}, and adapted to our terminology and notation.

\begin{theorem}\label{theo:ponja} 
Let $\C$ be an abelian category, $\X$ a left thick class, and $\Y$ a right thick class contained in $\X^\wedge$, such that $\omega := \X \cap \Y$ is an $\X$-injective relative cogenerator in $\X$. Then
\[
\Y = \omega^\wedge = \X^\wedge \cap \X^\perp = \X^\wedge \cap \X^{\perp_1} = \mathsf{Inj}^{< \infty}_{\X}(\X^\wedge) \cap \omega^\perp.
\]
\end{theorem}

\begin{myproof}
By \emph{Proposition~\ref{AB11}}, we have $\omega^\wedge = \X^\wedge \cap \X^\perp$. We now show that $\omega^\wedge = \Y$. Indeed, since $\omega \subseteq \Y$, we have $\omega^\wedge \subseteq \Y^\wedge$, and since $\Y$ is right thick, we obtain $\Y^\wedge = \Y$. The inclusion $\omega^\wedge \subseteq \Y$ follows. Now let $Y \in \Y$. Knowing that $\Y \subseteq \X^\wedge$, we can apply \emph{Theorem~\ref{AB4}} to get a short exact sequence 
\[
0 \to K \to X \to Y \to 0
\] 
with $X \in \X$ and $K \in \omega^\wedge \subseteq \Y$. Since $\Y$ is closed under extensions, we have that $X \in \X \cap \Y =: \omega$. It follows that $Y \in \omega^\wedge$, that is, $\Y \subseteq \omega^\wedge$. So far, we have proven the equalities $\Y = \omega^\wedge = \X^\wedge \cap \X^\perp$. 

For the third equality, it is clear that $\mathcal{X}^\wedge \cap \mathcal{X}^\perp \subseteq \mathcal{X}^\wedge \cap \mathcal{X}^{\perp_1}$. Now let $Z \in \mathcal{X}^\wedge \cap \mathcal{X}^{\perp_1}$. By \emph{Theorem~\ref{AB4}} again, there exists a short exact sequence
\[
0 \to Z \to W \to X \to 0
\] 
with $W \in \omega^\wedge$ and $X \in \mathcal{X}$. Since $\Ext^1_{\C}(X,Z) = 0$, the previous sequence splits, $Z$ is a direct summand of $W \in \omega^\wedge = \Y$, and so $Z \in \Y = \mathcal{X}^\wedge \cap \mathcal{X}^\perp$. Hence, the third equality follows.

Finally, we prove $\mathcal{Y} = \mathsf{Inj}^{< \infty}_{\mathcal{X}}(\mathcal{X}^\wedge) \cap \omega^\perp$. Note that $\Y \subseteq \X^\perp \subseteq \omega^\perp$, and by \emph{Theorem~\ref{AB14}}, we know 
\[
\Y = \omega^\wedge \subseteq (\omega^\wedge)^\vee = \mathsf{Inj}^{< \infty}_{\X}(\X^\wedge).
\] 
So $\Y \subseteq \mathsf{Inj}^{< \infty}_{\X}(\X^\wedge) \cap \omega^\perp$. 

Now let $Z \in \mathsf{Inj}^{< \infty}_{\mathcal{X}}(\mathcal{X}^\wedge) \cap \omega^\perp$, with $\mathrm{id}_{\X}(Z) = k < \infty$. We use induction on $k$ to prove $Z \in \Y$. Suppose $k = 1$. Given $X \in \X$, there is an exact sequence 
\[
0 \to X \to W \to X' \to 0
\]
with $W \in \omega$ and $X' \in \X$. Then there is an exact sequence
\[
\Ext^1_{\C}(W,Z) \to \Ext^1_{\C}(X,Z) \to \Ext^2_{\C}(X',Z)
\]
where $\Ext^1_{\C}(W,Z) = 0$ since $Z \in \omega^\perp$, and $\Ext^2_{\C}(X',Z) = 0$ since $\id_{\X}(Z) = 1$. It follows $Z \in \X^\wedge \cap \X^{\perp_1} = \Y$. Now if $\mathrm{id}_{\mathcal{X}}(Z) = k$, then in a similar way we can show that $\Ext^{k-1}_{\C}(X,Z) = 0$ for every $X \in \mathcal{X}$. Repeating this procedure, we finally get that $Z \in \mathcal{X}^{\perp_1}$, and hence the result follows.  
\end{myproof}

%%%%%%%%%%%%%%%%%%%%%%%%%%%%%%%%%%%%%%%%%%%%%%%%%%%%%%%%%%%%%
%%%%%%%%%%%%%%%%%%%%%%%%%%%%%%%%%%%%%%%%%%%%%%%%%%%%%%%%%%%%%

\section{Relative cotorsion pairs}\label{sec:relative_pairs}

The section is devoted to present the notion of cotorsion pairs relative to a thick sub-category of an abelian category. We begin this section recalling the concept of cotorsion pairs in exact categories, and then introduce the relative $\Sb$-cotorsion pairs as complete cotorsion pairs in a thick sub-category $\Sb$ of an abelian category $\C$. Later, we provide a characterization for this concept which, along with the results presented in \emph{Section~\ref{sec:ABtheory}}, will allow us to construct relative cotorsion pairs from Frobenius pairs. Throughout this section, the term ``sub-category'' means ``full sub-category''.

%%%%%%%%%%%%%%%%%%%%%%%%%%%%%%%%%%%
%%%%%%%%%%%%%%%%%%%%%%%%%%%%%%%%%%%

\subsection{Cotorsion pairs in exact categories}

The notion of cotorsion pair was first introduced by L. Salce in \emph{\cite{S}}. It is the analog of a torsion pair where the functor $\Hom_\C(-,-)$ is replaced by $\Ext^1_\C(-,-)$. Roughly speaking, two classes $\A$ and $\B$ of objects in $\C$ form a cotorsion pair if they are orthogonal to each other with respect to the functor $\Ext^1_{\C}(-,-)$. 

Most of the homological algebra done for abelian categories carries over to exact categories (This can be appreciated in the detailed survey \emph{\cite{Buhler}} of exact categories written by T. B\"uhler). So it is not surprising the existence of the notion of cotorsion pair for exact categories.

The concept of cotorsion pairs in exact categories given below is due to J. Gillespie (See \emph{\cite[Definition 2.1]{GiExact}}). Due to the purposes of this paper, we split this concept into \emph{left} and \emph{right} cotorsion pairs. 

Let $(\mathcal{E}, \tau)$ be an exact category, where $\tau$ is a class of ``short exact sequences'' in $\mathcal{E}$. The axioms of exact categories allows us to construct the extension functors $\Ext^i_{\mathcal{E}}(-,-)$ as in abelian categories, that is, using the Baer-Yoneda description. For example, in the case $i = 1$, $\Ext^1_{\mathcal{E}}(A,B)$ is the abelian group of equivalence classes of short exact sequences
\begin{align}\label{eqn:exact}
0 & \to B \to C \to A \to 0 
\end{align}
in $\tau$, where the zero element is given by the class of the split sequence 
\[
0 \to B \to A \oplus B \to A \to 0.
\]
Every morphism $B \to C$ appearing in a short exact sequence \emph{(\ref{eqn:exact})} in $\tau$ is called an \textbf{admissible monomorphism}. Dually, we have \textbf{admissible epimorphisms}. An object $I \in \mathcal{E}$ is \textbf{$\tau$-injective} if any admissible monomorphism $I \to C$ splits, or equivalently, if $\Ext^1_{\mathcal{E}}(X,I) = 0$ for every $X \in \mathcal{E}$. \textbf{$\tau$-Projective}\footnote{If there is no need to specify the class $\tau$ of short exact sequences in $\mathcal{E}$, $\tau$-injective and $\tau$-projective objects will be simply referred as projective and injective objects in $\mathcal{E}$.} objects in $\mathcal{E}$ have a dual description.

\begin{example}\label{ex:thick_are_exact}
Let $\C$ be an abelian category. Then, $\C$ is exact with $\tau$ formed by the family of all short exact sequences in $\C$. Thus, every result or definition presented below in the context of exact categories will be also valid for abelian categories.

Now let $\Sb$ be a thick sub-category of $\C$. Then $\Sb$ is also an exact category, where the family $\tau$ is formed by the short exact sequences \emph{(\ref{eqn:exact})} such that $A, B, C \in \Sb$. This exact category $\Sb$ is not necessarily abelian. In fact, $\Sb$ is abelian if, and only if, $\Sb$ is an \emph{admissible sub-category} of $\C$ (See \emph{\cite[Proposition 2.3]{MMSS}} for details).
\end{example}

\begin{definition} 
Let $(\mathcal{E}, \tau)$ be an exact category, and $\F$ and $\G$ be classes of objects in $\mathcal{E}$. The pair $(\F, \G)$ is said to be a \textbf{left cotorsion pair} if $\F = {}^{\perp_1}\G$\footnote{Orthogonal classes in exact categories are defined as in abelian categories.}. We say that a left cotorsion pair $(\F,\G)$ in $\mathcal{E}$ is \textbf{complete} if for every $X \in \mathcal{E}$ there exists a short exact sequence
\begin{align}\label{eqn:precover}
0 \to G \to F \to X \to 0
\end{align}
with $F \in \F$ and $G \in \G$. Note that in this case, the morphism $F \to X$ in the sequence \emph{(\ref{eqn:precover})} is a special $\F$-pre-cover\footnote{(Pre-)cover and (pre-)envelopes in exact categories are defined as in abelian categories.}.

Dually, we have the concept of (\textbf{complete}) \textbf{right cotorsion pair}. 

Finally, we say that $(\F,\G)$ is a \textbf{cotorsion pair} if it is both a left and a right cotorsion pair. A cotorsion pair $(\F,\G)$ is \textbf{complete} if it is complete as a left and as a right cotorsion pair. 
\end{definition}

An exact category $(\mathcal{E},\tau)$ is said to have  \textbf{enough $\tau$-injectives} if for every object $X \in \mathcal{E}$ there exists an admissible monomorphisms $X \to I$, where $I$ is a $\tau$-injective object of $\mathcal{E}$. If $\mathcal{E}$ satisfies the dual property, we say that $\mathcal{E}$ have \textbf{enough $\tau$-projectives}. 

The following result states that the dual notions of pre-covers and pre-envelopes are equivalent when considered with respect to a cotorsion pair. The reader can see the original proof in \emph{\cite[Corollary 2.4]{S}}, whose arguments carry over to exact categories.

\begin{lemma} 
For a cotorsion pair $(\F,\G)$ in an exact category $\mathcal{E}$, with enough projectives and injectives, the following conditions are equivalent.
\begin{itemize}[itemsep=2pt,topsep=0pt]
\item[$\ai$] Every object in $\mathcal{E}$ has a special $\F$-pre-cover.

\item[$\bi$] Every object in $\mathcal{E}$ has a special $\G$-pre-envelope.
\end{itemize}
\end{lemma}

\begin{definition}\label{def:hereditary} 
Let $(\mathcal{F,G})$ be a cotorsion pair in an exact category $\mathcal{E}$. We say that $(\mathcal{F,G})$ is \textbf{left hereditary} if $\F$ is resolving (in $\mathcal{E}$). Dually, we have the notion of \textbf{right hereditary} cotorsion pair in $\mathcal{E}$. A \textbf{hereditary} cotorsion pair in $\mathcal{E}$ is a cotorsion pair which is both left and right hereditary. 
\end{definition}

The dual notions of resolving and coresolving classes are equivalent when considered with respect to a cotorsion pair, as specified in the following result (See \emph{\cite[Theorem 1.2.10]{GR}}).

\begin{lemma}\label{lem:hereditary_exact} 
For a cotorsion pair $(\F,\G)$ in an exact category $(\mathcal{E},\tau)$, with enough $\tau$-projectives and $\tau$-injectives, the following conditions are equivalent.
\begin{itemize}[itemsep=2pt,topsep=0pt]
\item[$\ai$] $(\mathcal{F,G})$ is left hereditary.

\item[$\bi$] $(\mathcal{F,G})$ is right hereditary. 

\item[$\ci$] $\id_\F\,(\G) = 0$. 
\end{itemize}
\end{lemma}

Note that, in an hereditary cotorsion pair $(\F,\G)$ in an exact category $\mathcal{E}$, we have that $\F$ is left saturated and $\G$ is right saturated (in $\mathcal{E}$). \newpage

%%%%%%%%%%%%%%%%%%%%%%%%%%%%%%%%%%
%%%%%%%%%%%%%%%%%%%%%%%%%%%%%%%%%%

\subsection{Cotorsion pairs relative to thick sub-categories}

From now on, we focus on a special type of complete cotorsion pair in a thick sub-category $\Sb$ of an abelian category $\C$. We provide in \emph{Proposition~\ref{CP1}} a characterization for these pairs, more suitable to study properties and construct examples related to AB Theory.

\begin{definition}\label{def:left_cotorsion_pair}  
\sloppypar{Let $\Sb$ be a thick sub-category of an abelian category $\C$ (and then an exact sub-category). Let $\F$ and $\G$ be two classes of objects in $\C$. We say that $(\F,\G)$ is a \textbf{left $\Sb$-cotorsion pair} in $\C$ if $(\F,\G)$ is a complete left cotorsion pair in the exact category $\Sb$. Dually, we have the definition of \textbf{right $\Sb$-cotorsion pairs} in $\C$. Finally, we say that $(\F,\G)$ is an \textbf{$\Sb$-cotorsion pair} if it is both a left and a right $\Sb$-cotorsion pair \mbox{in $\C$.}}
\end{definition}

\begin{remark}\label{rem:left_pairs_properties} Note the following facts about $\Sb$-cotorsion pairs and orthogonal classes:
\begin{itemize}[itemsep=2pt,topsep=0pt]
\item[$\ai$] If $(\F,\G)$ is a left $\Sb$-cotorsion pair, then $\F$ and $\G$ are sub-classes of  $\Sb$. 

\item[$\bi$] For any $\X$ and $\Y$ classes of objects in $\C,$ we consider the relative perpendicular classes 
\[
{}^{\perp_{i,\Y}}\X := {}^{\perp_i}\X \cap \Y \mbox{ \ and \ } \X^{\perp_{i,\Y}} := \X^{\perp_i} \cap \Y.
\]
It follows that if $\Sb$ is a thick sub-category of $\C$, the relative perpendicular classes in $\C$ with respect to $\Sb$ are set according to the previous definition. Hence, if $(\F,\G)$ is a left (resp., right) $\Sb$-cotorsion pair, then 
\[
\F = {}^{\perp_{1,\Sb}}\G = {}^{\perp_1}\G \cap \Sb \mbox{ (resp., $\G = \F^{\perp_{1,\Sb}} = \F^{\perp_1} \cap \Sb$)}. 
\]

\item[$\ci$] If $(\F,\G)$ is a left $\Sb$-cotorsion pair in $\C$, then $\F$ is closed under extensions, since $\F = {}^{\perp_{1,\Sb}}\G$. Dually, $\G$ is closed under extension for every right $\Sb$-cotorsion pair $(\F,\G)$ in $\C$.
\end{itemize}
\end{remark}

\begin{proposition}\label{CP1} 
Let $(\F,\G)$ be a pair of classes of objects in $\C$ and $\Sb$ a thick sub-category of $\C$. Consider the following conditions: 
\begin{itemize}[itemsep=2pt,topsep=0pt]
\item[$\scpuno$] $\F, \G \subseteq \Sb$, and $\F$ is closed under direct summands in $\C$. 

\item[$\scpdos$] $\F, \G \subseteq \Sb$, and $\G$ is closed under direct summands in $\C$. 

\item[$\scptres$] $\Ext^1_\C(\F,\G) = 0$.

\item[$\scpcuatro$] For every object $S \in \Sb$, there exists an epic $\F$-pre-cover $\varphi \colon F \to S$, with $\Ker(\varphi) \in \G$.\footnote{Note that if $\G \subseteq \Sb$ and condition $\scptres$ holds, then these $\F$-pre-covers are special since $\G \subseteq \F^{\perp_{1,\Sb}}$. However, the inclusion $\G \supseteq \F^{\perp_{1,\Sb}}$ is not true in general, and so not every special $\F$-pre-cover of $S$ has kernel in $\G$.}

\item[$\scpcinco$] For every object $S \in \Sb$, there exists a monic $\G$-pre-envelope $\psi \colon S \to G$, with $\Coker(\psi) \in \F$.
\end{itemize}
Then the following statements hold true:
\begin{itemize}[itemsep=2pt,topsep=0pt]
\item[$\ai$] $(\F,\G)$ is a left $\Sb$-cotorsion pair in $\C$ if, and only if, $\F$ and $\G$ satisfy $\scpuno$, $\scptres$ and $\scpcuatro$.

\item[$\bi$] $(\F,\G)$ is a right $\Sb$-cotorsion pair in $\C$ if, and only if, $\F$ and $\G$ satisfy $\scpdos$, $\scptres$ and $\scpcinco$. 

\item[$\ci$] $(\F,\G)$ is an $\Sb$-cotorsion pair in $\C$ if, and only if, $\F$ and $\G$ satisfy conditions from $\scpuno$ to $\scpcinco$. \newpage
\end{itemize}
\end{proposition}

\begin{myproof} 
We only prove part $\ai$, since $\bi$ is dual, and $\ci$ follows from $\ai$ and $\bi$.

Suppose $(\F,\G)$ is a left $\Sb$-cotorsion pair. Then by \emph{Remark~\ref{rem:left_pairs_properties}}, we have $\F$ and $\G$ are sub-classes of $\Sb$. The same remark tell us that $\F = {}^{\perp_{1}}\G \cap \Sb$, and so $\F$ is closed under direct summands since both ${}^{\perp_{1}}\G$ and $\Sb$ are. Hence, condition $\scpuno$ follows. Condition $\scptres$ follows by the equality $\F = {}^{\perp_{1,\Sb}}\G$, and by the fact that $\F$ and $\G$ are sub-classes of $\Sb$. Finally, $\scpcuatro$ follows by the fact that $(\F,\G)$ is a complete left cotorsion pair in $\Sb$. The right-hand morphisms in the exact sequence \emph{(\ref{eqn:precover})} are epic $\F$-pre-covers with kernel in $\G$. 

Now suppose that $\F$ and $\G$ are two classes of objects in $\C$ satisfying conditions $\scpuno$, $\scptres$ and $\scpcuatro$. Then by $\scpuno$, we have that $\F$ and $\G$ are sub-classes of $\Sb$. Let us check the equality $\F = {}^{\perp_{1,\Sb}}\G$. The inclusion $\F \subseteq {}^{\perp_{1,\Sb}}\G$ follows by the equality $\Ext^1_{\C}(\F,\G) = 0$. Now let $S \in {}^{\perp_{1,\Sb}}\G$. By $\scpcuatro$, there exists a short exact sequence 
\[
0 \to \Ker(\varphi) \to F \xrightarrow{\varphi} S \to 0
\]
with $F \in \F$ and $\Ker(\varphi) \in \G$. Then $\Ext^1_{\C}(S,\Ker(\varphi)) = 0$, and so the previous short exact sequence splits. It follows that $S$ is a direct summand of $F \in \F$, and since $\F$ is closed under direct summands by $\scpuno$, we have $S \in \F$, and hence the inclusion ${}^{\perp_{1,\Sb}}\G \subseteq \F$ follows. So far we have proven that $(\F,\G)$ is a left cotorsion pair in $\Sb$, and the completeness of $(\F,\G)$ follows by $\scpcuatro$. 
\end{myproof}

%%%%%%%%%%%%%%%%%%%%%%%%%%%%%%%%%%
%%%%%%%%%%%%%%%%%%%%%%%%%%%%%%%%%%

\subsection{Relative cotorsion pairs from Frobenius pairs}

The characterization of (left and right) $\Sb$-cotorsion pairs given in \emph{Proposition~\ref{CP1}} allows us to construct easily cotorsion pairs from Frobenius pairs (See \emph{Definition~\ref{def:Frobenius_pair}}). Later on we will study correspondences between these two notions.

\begin{theorem}\label{CP2} 
Let $(\X,\omega)$ be a left Frobenius pair in $\C$. Then, the following statements hold:
\begin{itemize}[itemsep=2pt,topsep=0pt]
\item[$\ai$] $(\X,\omega^\wedge)$ is an $\X^\wedge$-cotorsion pair in $\C$.

\item[$\bi$] The following equalities hold 
\[
\omega^\wedge = \X^\perp \cap \X^\wedge, \mbox{ \ } \omega = \X \cap \omega^\wedge, \mbox{ \ and \ } \X = \X^\wedge \cap {}^{\perp}\omega = \X^\wedge \cap {}^\perp(\omega^\wedge).
\]
\end{itemize}
Dually, if $(\nu,\Y)$ is a right Frobenius pair in $\C$, then the following statements hold:
\begin{itemize}[itemsep=2pt,topsep=0pt]
\item[$\iroman$] $(\nu^\vee, \Y)$ is a $\Y^\vee$-cotorsion pair in $\C$.

\item[$\iiroman$] The following equalities hold:
\[
\nu^\vee = {}^\perp\Y \cap \Y^\vee, \mbox{ \ } \nu = \Y \cap \nu^\vee, \mbox{ \ and \ } \Y = \Y^\vee \cap \nu^\perp = \Y^\vee \cap (\nu^\vee)^\perp.
\]
\end{itemize}
\end{theorem}

\begin{myproof} 
By \emph{Theorem~\ref{AB9}}, we know that $\X^\wedge$ is a thick sub-category in $\C$. On the other hand, by \emph{Lemma~\ref{AB1}} we have $\id_\X(\omega^\wedge) = \id_\X(\omega) = 0$, and so $\scptres$ follows. Furthermore, \emph{Proposition~\ref{AB11}} gives us that $\omega^{\wedge} = \X^{\perp} \cap \X^{\wedge}$, and hence $\omega^{\wedge}$ is closed under direct summands. Condition $\scpdos$ then follows; while $\scpuno$ is true since $\X$ is left thick. Note that $\scpcuatro$ and $\scpcinco$ hold by \emph{Theorem~\ref{AB4}}, and hence we conclude by \emph{Proposition~\ref{CP1}} that $(\X,\omega^\wedge)$ is an $\X^\wedge$-cotorsion pair in $\C$. Furthermore, $\X^{\wedge} \cap {}^{\perp}\omega = \X = \X^{\wedge} \cap {}^{\perp}(\omega^{\wedge})$ follow from \emph{Proposition~\ref{AB12}}. Finally, the equality $\omega = \X \cap \omega^\wedge$ follows from \emph{Proposition~\ref{AB2}}.
\end{myproof}

If we impose an exact condition on $(\X,\omega)$ in the previous theorem, then it is possible to construct another $\X^\wedge$-cotorsion pair in $\C$.

\begin{theorem}\label{theorem:projective_pair}
Let $(\X,\omega)$ be a strong left Frobenius pair in $\C$. Then the following conditions hold:
\begin{itemize}[itemsep=2pt,topsep=0pt]
\item[$\ai$] $(\omega,\X^\wedge)$ is a $\X^\wedge$-cotorsion pair in $\C$. 

\item[$\bi$] $\omega^\wedge = \Thick(\omega)$. 
\end{itemize}

Dually, if $(\nu,\Y)$ is a strong right Frobenius pair in $\C$, then:
\begin{itemize}[itemsep=2pt,topsep=0pt]
\item[$\iroman$] $(\Y^\vee, \nu)$ is a $\Y^\vee$-cotorsion pair in $\C$.

\item[$\iiroman$] $\nu^\vee = \Thick(\nu)$.
\end{itemize}
\end{theorem}

\begin{myproof} 
\begin{itemize}[itemsep=2pt,topsep=0pt]
\item[$\ai$] By \emph{Theorem~\ref{AB9}}, we get that $\X^\wedge$ is a thick sub-category of $\C$. Now note that $\omega$ is closed under direct summands by hypothesis, and $\X^\wedge$ satisfies the same property since it is thick. Moreover, $\omega \subseteq \X^\wedge$. Then conditions $\scpuno$ and $\scpdos$ follows. 

To show $\mathsf{(scp3)}$, it suffices to use \emph{Lemma~\ref{AB1}} and the condition $\pd_{\X}(\omega) = 0$ as follows:
\[
0 = \pd_{\X}(\omega) = \id_{\omega}(\X) = \id_{\omega}(\X^\wedge). 
\]

Note that $\scpcinco$ is trivial since $0 \in \omega$ (which follows using the property that $\omega$ is closed under direct summands). 

Finally, we show $\scpcuatro$. Let $Y \in \X^\wedge$. By \emph{Theorem~\ref{CP2}}, there exists a short exact sequence 
\[
0 \to K \xrightarrow{\alpha} X \xrightarrow{\beta} Y \to 0
\] 
with $K \in \omega^\wedge$ and $X \in \X$. On the other hand, there exists a short exact sequence 
\[
0 \to X' \xrightarrow{i} W \xrightarrow{p} X \to 0
\]
with $X' \in \X$ and $W \in \omega$, since $\omega$ is a relative generator in $\X$. Now consider the composition $q = \beta \circ p$. By Snake Lemma and the fact that $p$ is epic, we obtain a short exact sequence 
\[
0 \to \Ker(p) \to \Ker(q) \to \Ker(\beta) \to 0
\]
where $\Ker(p) = X' \in \X$ and $\Ker(\beta) = K \in \omega^\wedge \subseteq \X^\wedge$. Then $\Ker(q) \in \X^\wedge$ since $\X^\wedge$ is thick. Hence, we have a short exact sequence
\[
0 \to \Ker(q) \to W \xrightarrow{q} Y \to 0
\]
with $W \in \omega$ and $\Ker(q) \in \X^\wedge$, proving $\scpcuatro$.

\item[$\bi$] Indeed, by \emph{Theorem~\ref{AB14}}, we have 
\[
\X \cap \Thick(\omega) = \X \cap \mathsf{Inj}^{< \infty}_{\X}(\X^\wedge) = \mathsf{Inj}^{< \infty}_{\X}(\X).
\]
By \emph{Lemma~\ref{AB3}}, we have $\mathsf{Inj}^{< \infty}_{\X}(\X) = \X \cap \omega^\vee$, and by \emph{Proposition~\ref{AB2}} we have $\X \cap \omega^\vee = \omega$. Hence, $\X \cap \Thick(\omega) = \omega$ holds true. Setting $\Y := \Thick(\omega)$ in \emph{Theorem~\ref{theo:ponja}}, it follows that $\omega^\wedge = \Thick(\omega)$. 
\end{itemize}
\end{myproof}

\begin{corollary}
If $(\X,\omega)$ is a strong left Frobenius pair in $\C$, then $\omega$ is the class of projective objects in the exact sub-category $\X^\wedge \subseteq \C$. Moreover, $\X^\wedge$ has enough projectives. 
\end{corollary}

\begin{myproof}
It follows from \emph{Theorem~\ref{theorem:projective_pair}}, since $\omega = {}^{\perp_{1,\X^\wedge}}(\X^\wedge) = \Proj(\X^\wedge)$.
\end{myproof}

So far we know that the concept of $\Sb$-cotorsion pair is a description of the notion of (left and right) completeness of cotorsion pairs in $\Sb$. In the next section, we introduce the property of ``being hereditary'' for $\Sb$-cotorsion pairs in $\C$, and compare it with the standard definition of hereditary cotorsion pairs in $\Sb$ (that is, \emph{Definition~\ref{def:hereditary}}).

%%%%%%%%%%%%%%%%%%%%%%%%%%%%%%%%
%%%%%%%%%%%%%%%%%%%%%%%%%%%%%%%%

\subsection{Hereditary relative cotorsion pairs}

We now study the analogous notion of hereditary cotorsion pairs in the relative context (See \emph{Definition~\ref{def:hereditary_relative}}). However, this notion will not be equivalent to that of hereditary cotorsion pairs in exact categories. We later in this section present some conditions under which the $\X^\wedge$-cotorsion pairs obtained in \emph{Theorem~\ref{CP2}} and \emph{Theorem~\ref{theorem:projective_pair}} are hereditary in the sense of the following definition.

\begin{definition}\label{def:hereditary_relative} 
We say that an $\Sb$-cotorsion pair $(\F,\G)$ in $\C$ is \textbf{left hereditary} if $\F$ is a resolving class in $\Sb$. Furthermore, we say that $(\F,\G)$ is \textbf{left strong hereditary} if $\F$  is a resolving class in $\C$. 

The notions of \textbf{right} (\textbf{strong}) \textbf{hereditary} and (\textbf{strong}) \textbf{hereditary} $\Sb$-cotorsion pairs in $\C$ are defined similarly. 
\end{definition}

 Is it possible to establish some condition under which the projective objects of $\C$ and $\Sb$ coincide? This question is settled in the following result.

\begin{proposition}\label{prop:same_projectives}
Let $\Sb$ be a thick sub-category of $\C$. If $\Proj(\C) \subseteq \Sb$ and $\C$ has enough projectives, then $\Proj(\Sb) = \Proj(\C)$. 
\end{proposition}

\begin{myproof}
First, note that the inclusion $\Proj(\C) \subseteq \Proj(\Sb)$ is clear since $\Proj(\C) \subseteq \Sb$. Now, let $Q \in \Proj(\Sb)$. Since $\C$ has enough projectives, there exists a short exact sequence
\begin{align}\label{eqn:res}
0 & \to K \to P \to Q \to 0
\end{align}
with $P \in \Proj(\C) \subseteq \Proj(\Sb)$. We also know that $\Sb$ is a thick sub-category, and so $K \in \Sb$. It follows \emph{(\ref{eqn:res})} is a short exact sequence in $\Sb$, and since $Q \in \Proj(\Sb)$, we have \emph{(\ref{eqn:res})} splits, and so $Q$ is a direct summand of $P$, which in turn implies $Q \in \Proj(\C)$. 
\end{myproof}

As a consequence of the previous result, we have the following remark.

\begin{remark}\label{rem:resolving_thick} 
Let $(\F,\G)$ be a cotorsion pair in $\Sb.$ If $\C$ has enough projectives, then $\F$ is resolving in $\C$ if, and only if, $\F$ is resolving in $\Sb$ and $\Proj(\C) \subseteq \F$. 
\end{remark}

\begin{theorem}\label{CP3} 
Let $(\F,\G)$ be a left strong hereditary $\Sb$-cotorsion pair in $\C$ and $\omega := \F \cap \G$. If $\C$  has enough projectives, then the following statements hold:
\begin{itemize}[itemsep=2pt,topsep=0pt]
\item[$\ai$] $\Ext^i_\C(\F,\G) = 0$ for $i\geq 1.$

\item[$\bi$] $(\F,\omega)$ is a left Frobenius pair in $\C$.

\item[$\ci$] $(\omega,\G)$ is a right Frobenius pair in $\C$. 

\item[$\di$] If $\G \subseteq \F^{\wedge},$ then $(\F,\G)$ is an $\F^{\wedge}$-cotorsion pair and
\[
\G = \omega^{\wedge} = \F^{\wedge} \cap \F^{\perp} \mbox{ \ and \ } \F = \F^{\wedge} \cap {}^{\perp}\omega = \F^{\wedge} \cap {}^{\perp}(\omega^{\wedge}).
\]
\end{itemize}
\end{theorem}

\begin{myproof} 
\begin{itemize}[itemsep=2pt,topsep=0pt]
\item[$\ai$] Let $F \in \F$ and $G \in \G$. Since $\F$ is resolving (in $\C$) and $\C$ has enough projectives, we have that $\Omega^{i-1}(F) \subseteq \F$, for any $i \geq 1$. Then we have that $\Ext^i_\C(F,G) \simeq \Ext^1_\C(F', G) = 0$, for every $F' \in \Omega^{i-1}(F)$.

\item[$\bi$] First, we need to check that $\F$ is left thick. Since $\F$ is resolving, it is only left to check that $\F$ is closed under direct summands in $\C$, which follows from the equality $\F = {}^{\perp_1}\G \cap \Sb$. 

Note that $\omega$ is closed under direct summands since both $\F$ and $\G$ are. So, in order to show that $(\F,\omega)$ is a left Frobenius pair in $\C$, it remains to check that $\omega$ is an $\F$-injective relative cogenerator in $\F$. Part $\ai$ implies that $\omega$ is $\F$-injective, since $\id_\F(\omega) \leq \id_\F(\G) = 0$. Now let $F \in \F$. By $\scpcinco$, there exists a short exact sequence
\[
0 \to F \to W \to F' \to 0
\]
in $\C$ with $F' \in \F$ and $W \in \G$. Since $\F$ is closed under extensions, we obtain $W \in \F \cap \G =: \omega$, proving that $\omega$ is a relative cogenerator in $\F$. Hence, $\bi$ follows.

\item[$\ci$] We first check that $\G$ is right thick. The equality $\G = \F^{\perp_1} \cap \Sb$ implies that $\G$ is closed under extensions and direct summands in $\C$. Now consider a short exact sequence 
\[
0 \to A \to B \to C \to 0
\] 
in $\Sb$, with $A, B \in \G$. Applying the functor $\Hom_\C(F,-)$ to the above exact sequence, with $F$ running over $\F$, we get an exact sequence
\[
\Ext^1_\C(F,B) \to \Ext^1_\C(F,C) \to \Ext^2_\C(F,A)
\]
where $\Ext^1_{\C}(F,B) = 0$ and $\Ext^2_{\C}(F,A) = 0$ since $\id_{\F}(\G) = 0$ by $\ai$. Hence, it follows that $C \in \F^{\perp_1} \cap \Sb = \G$. Thus, we have that $\G$ is right thick. The rest of the proof follows as in $\bi$.

\item[$\di$] Since $(\F,\omega)$ is a left Frobenius pair in $\C,$ we get from \emph{Theorem~\ref{CP2}} that
\[
\omega^{\wedge} = \F^{\wedge} \cap \F^{\perp}\;\text{ and }\;\F^{\wedge}\cap{}^{\perp}\omega = \F = \F^{\wedge} \cap {}^{\perp}(\omega^{\wedge}).
\]
We assert that $\G = \omega^{\wedge}$. Indeed, from $\ci$ we get that $\omega^\wedge \subseteq \G$. In order to prove that $\G \subseteq \omega^{\wedge}$, it suffices to see that $\G \subseteq \F^{\wedge} \cap \F^{\perp}$. But this follows from $\ai$, since $\G \subseteq \F^{\wedge}$.
\end{itemize}
\end{myproof}

The following result shows us how to obtain strong hereditary relative cotorsion pairs from hereditary cotorsion pairs in abelian categories.

\begin{corollary}\label{CP4} 
Let $(\F,\G)$ be a left hereditary complete cotorsion pair in $\C$ and $\omega := \F \cap \G$. If $\C$  has enough projectives, then  
$(\F, \G \cap \F^{\wedge})$ is a left strong hereditary $\F^\wedge$-cotorsion pair in $\C$, and
\[
\omega^\wedge = \G \cap \F^\wedge = \F^\wedge \cap \F^\perp \mbox{ \ and \ } \F = \F^\wedge \cap {}^\perp\omega = \F^\wedge \cap {}^\perp(\omega^\wedge).
\]
\end{corollary}

\begin{myproof} 
Note that $(\F,\G)$ is a left strong hereditary $\C$-cotorsion pair in $\C$. Then by \emph{Theorem~\ref{CP3}}, we get that $(\F,\omega)$ is a left Frobenius pair. Hence by $\emph{Theorem~\ref{CP2}}$, it follows that $(\F, \omega^\wedge)$ is a left strong hereditary $\F^\wedge$-cotorsion pair in $\C$, and furthermore, we have the equalities 
\[
\omega^\wedge = \F^\perp \cap \F^\wedge \mbox{ \ and \ } \F = \F^\wedge \cap {}^\perp\omega = \F^\wedge \cap {}^\perp(\omega^\wedge).
\]
where $\G = \F^{\perp_1} = \F^\perp$, since $\F$ is resolving and $\C$ has enough projectives. 
\end{myproof}

We close this section presented the conditions for which the $\X^\wedge$-cotorsion pairs $(\X, \omega^\wedge)$ and $(\omega, \X^\wedge)$ are strong hereditary.

\begin{theorem}\label{theo:induced_pairs_hereditary}
Let $(\X,\omega)$ be a strong left Frobenius pair in $\C$. Then the following conditions are equivalent:
\begin{itemize}[itemsep=2pt,topsep=0pt]
\item[$\ai$] $\Proj(\C) \subseteq \omega$ and $\Inj(\C) \subseteq \X^\wedge$.

\item[$\bi$] $(\omega,\X^\wedge)$ is a strong hereditary $\X^\wedge$-cotorsion pair in $\C$.

\item[$\ci$] $(\X, \omega^\wedge)$ is a strong hereditary $\X^\wedge$-cotorsion pair in $\C$.
\end{itemize} 
\end{theorem}

\begin{myproof} First, note by \emph{Theorems~\ref{CP2} and \ref{theorem:projective_pair}} that $(\X,\omega^\wedge)$ and $(\omega,\X^\wedge)$ are $\X^\wedge$-cotorsion pairs in $\C$, and $\omega^\wedge$ is thick. 
\begin{itemize}[itemsep=2pt,topsep=0pt]
\item[$\bullet$] $\ai \Longrightarrow \bi$: It suffices to show that $\omega$ is pre-resolving in $\C$ and $\X^\wedge$ is pre-coresolving in $\C$. The latter fact holds since $\X^\wedge$ is thick. The equality $\omega = {}^{\perp_1}(\X^\wedge) \cap \X^\wedge$ implies that $\omega$ is closed under extensions. To show that $\omega$ is also closed under taking kernel of epimorphisms in $\omega$, suppose we are given a short exact sequence
\begin{align*}
0 & \to A \to B \to C \to 0
\end{align*}
with $B, C \in \omega$. Since $\omega \subseteq \X$, we have $B, C \in \X$, and so $A \in \X$ since $\X$ is left thick. On the other hand, the inclusion $\omega \subseteq \omega^\wedge$ and the fact that $\omega^\wedge$ is thick implies that $A \in \X \cap \omega^\wedge$, where $\X \cap \omega^\wedge = \omega$ by \emph{Theorem~\ref{CP2}}. Hence $A \in \omega$, and $\bi$ follows.

\item[$\bullet$] $\bi \Longrightarrow \ci$: Suppose the $\X^\wedge$-cotorsion pair $(\omega,\X^\wedge)$ in $\C$ is strong hereditary. Then $\Proj(\C) \subseteq \omega$ and $\Inj(\C) \subseteq \X^\wedge$. The class $\omega^\wedge$ is pre-coresolving, since it is thick. Note also that $\Inj(\C) \subseteq \X^\perp$. It follows $\Inj(\C) \subseteq \X^\perp \cap \X^\wedge = \omega^\wedge$, where the last equality comes from \emph{Theorem~\ref{CP2}}. Hence, $\omega^\wedge$ is coresolving in $\C$. Concerning the class $\X$, we have on the one hand that $\X$ is pre-resolving since $\X = \Thick^{-}(\X)$. On the other hand, we know by hypothesis that $\Proj(\C) \subseteq \omega \subseteq \X$. Then $\X$ is a resolving class in $\C$. Hence, $(\X,\omega^\wedge)$ is a strong hereditary $\X^\wedge$-cotorsion pair in $\C$, and $\ci$ follows.

\item[$\bullet$] $\ci \Longrightarrow \ai$: Now suppose $(\X,\omega^\wedge)$ is a strong hereditary $\X^\wedge$-cotorsion pair in $\C$. Since $\omega^\wedge$ is coresolving in $\C$, we have $\Inj(\C) \subseteq \omega^\wedge \subseteq \X^\wedge$. On the other hand, the inclusion $\Proj(\C) \subseteq \X$ and the equality $\omega = {}^{\perp_1}(\X^\wedge) \cap \X^\wedge$ imply that $\Proj(\C) \subseteq \omega$. Therefore, $\ai$ holds. 
\end{itemize}
\end{myproof}

%%%%%%%%%%%%%%%%%%%%%%%%%%%%%%%%
%%%%%%%%%%%%%%%%%%%%%%%%%%%%%%%%

\subsection{Some examples of relative cotorsion pairs}

We devote this section to study several examples of $\Sb$-cotorsion pairs. We will present Gorenstein and Ding-Chen Homological Algebra encoded as particular instances of the setting presented before by Frobenius pairs.

%%%%%%%%%%%%%%%%%%%%%%%%%%%%%%%%%%%%%
%%%%%%%%%%%%%%%%%%%%%%%%%%%%%%%%%%%%%

\subsubsection*{Pairs from Gorenstein-projective and Gorenstein-injective modules:}

Consider $\C := \mathsf{Mod}(R)$, $\X := \GProj(R)$ and $\omega := \Proj(R)$. The class $\X$ is left thick by \emph{\cite[Theorem 2.5]{Holm}} and $\omega$ is $\X$-injective by \emph{\cite[Proposition 2.3]{Holm}}. By \emph{Example~\ref{ex:G-projective_strong_Frobenius_pair}}, we know that $(\GProj(R), \Proj(R))$ is a strong left Frobenius pair in $\mathsf{Mod}(R)$. Then by \emph{Theorems~\ref{CP2}} and \emph{\ref{theorem:projective_pair}} we have that $(\GProj(R), \Proj^{< \infty}(R))$ and $(\Proj(R), \GProj^{< \infty}(R))$ are $\GProj^{< \infty}(R)$-cotorsion pairs in $\mathsf{Mod}(R)$. This pairs are not in general left strong hereditary, since the inclusions $\Inj(R) \subseteq \GProj^{< \infty}(R)$ and $\Inj(R) \subseteq \Proj^{< \infty}(R)$ are not necessarily true. Note that these pairs are (left and right) hereditary as cotorsion pairs in the exact category $\GProj^{< \infty}(R)$. However, this two notions of hereditary cotorsion pairs coincide in the case $R$ is a Gorenstein ring, where the equalities $\GProj^{< \infty}(R) = \mathsf{Mod}(R)$ and $\Proj^{< \infty}(R) = \Inj^{< \infty}(R)$ hold.

Recall from \emph{Example~\ref{ex:G-projective_compatible}} the equality $\Proj(R) = \GProj(R) \cap \Proj^{< \infty}(R)$. Furthermore, the following equalities also hold:
\begin{align*}
\Proj^{< \infty}(R) & = \GProj(R)^\perp \cap \GProj^{< \infty}(R), \\
\GProj(R) & = \GProj^{< \infty}(R) \cap {}^\perp \Proj(R) = \GProj^{< \infty}(R) \cap {}^\perp \Proj^{< \infty}(R).
\end{align*}
These equalities are generalizations of the orthogonality relations $\Proj^{< \infty}(R) = \GProj(R)^\perp$ and $\GProj(R) = {}^\perp\Proj^{< \infty}(R)$ which hold in the case $R$ is a Gorenstein ring, since in that case we have a cotorsion pair $(\GProj(R), \Proj^{< \infty}(R))$ and the equality $\GProj^{< \infty}(R) = \mathsf{Mod}(R)$ (see \emph{\cite[Remark 11.5.10]{EJ}}). Then the following result follows.

\begin{corollary} 
If $\GProj^{< \infty}(R) = \mathsf{Mod}(R)$, then $(\GProj(R), \Proj^{< \infty}(R))$ is a hereditary complete cotorsion pair in $\mathsf{Mod}(R)$.
\end{corollary}

We also have dual conclusions if we set $\Y := \GInj(R)$ and $\nu := \Inj(R)$ in \emph{Theorems~\ref{CP2}} and \emph{\ref{theorem:projective_pair}}.

%%%%%%%%%%%%%%%%%%%%%%%%%%%%%%%%%%%%%
%%%%%%%%%%%%%%%%%%%%%%%%%%%%%%%%%%%%%

\subsubsection*{Pairs from Ding-projective and Ding-injective modules:}

The results by H. Holm that we cited previously from \emph{\cite{Holm}}, along with the arguments proving them, carry over to the notions of Ding-projective and Ding-injective modules. 

Recall that a left $R$-module $M$ is said to be \textbf{Ding-projective} (or \textbf{D-projective} for short) if there exists an exact sequence of projective modules
\[
\bm{P} = \cdots \to P_1 \to P_0 \to P^0 \to P^1 \to \cdots
\] 
with $M = \Ker(P^0 \to P^1)$, such that $\Hom_R(\bm{P}, F)$ is an exact chain complex of abelian groups, for every flat module $F \in \Flat(R)$. We denote by $\DProj(R)$ the class of D-projective left $R$-modules. Dually, a left $R$-module $N$ is \textbf{Ding-injective} if there exists an exact sequence of injective modules
\[
\bm{I} = \cdots \to I_1 \to I_0 \to I^0 \to I^1 \to \cdots
\]
with $N = \Ker(I_0 \to I^0)$, such that $\Hom_R(J, \bm{I})$ is an exact chain complex of abelian groups, for every FP-injective module $J$ (that is, a left $R$-module in $\mathcal{FP}^{\perp_1}$, where $\mathcal{FP}$ denotes the class of finitely presented left $R$-modules).

Using a similar reasoning as in the example of G-projective modules, if we set $\C := \mathsf{Mod}(R)$, $\X := \DProj(R)$ as the class of D-projective modules, and $\omega := \Proj(R)$, then we have that the pair $(\DProj(R), \Proj(R))$ is strong left Frobenius. Hence  we have left strong hereditary $\DProj^{< \infty}(R)$-cotorsion pairs $(\DProj(R), \Proj^{< \infty}(R))$ and $(\Proj(R), \DProj^{< \infty}(R))$ in $\mathsf{Mod}(R)$, where $\DProj^{< \infty}(R)$) denotes the class of modules with finite D-projective dimension. The first of these cotorsion pairs is a generalization of the cotorsion pair $(\DProj(R), \Flat^{< \infty}(R))$ found by J. Gillespie in \emph{\cite{Gi}}. Moreover, we have the equalities:
\begin{align*}
\Proj(R) & = \DProj(R) \cap \Proj^{< \infty}(R), \\
\Proj^{< \infty}(R) & = \DProj(R)^\perp \cap \DProj^{< \infty}(R), \\
\DProj^{< \infty}(R) \cap {}^\perp \Proj(R) & = \DProj(R) = \DProj^{< \infty}(R) \cap {}^\perp(\Proj^{< \infty}(R)).
\end{align*}
In \emph{\cite[Proposition 3.8]{Gi}}, it is proven that a D-projective module is either projective or has flat dimension $\infty$. This can be written as the equality $\Proj(R) = \DProj(R) \cap \Flat^{< \infty}(R)$. It follows that non-projective flat modules are not Ding-projective. For this reason, we set $\omega := \Proj(R)$ instead of $\omega := \Flat(R)$. Thus, the first equality above can be extended to
\[
\DProj(R) \cap \Flat^{< \infty}(R) = \Proj(R) = \DProj(R) \cap \Proj^{< \infty}(R).
\]
Dually, if $\DInj(R)$ denotes the class of Ding-injective modules, we have a strong right Frobenius pair $(\Inj(R), \DInj(R))$, and the corresponding $\DInj^{< \infty}(R)$-cotorsion pairs obtained from it.

%%%%%%%%%%%%%%%%%%%%%%%%%%%%%%%%%%%%%
%%%%%%%%%%%%%%%%%%%%%%%%%%%%%%%%%%%%%

\subsubsection*{Pairs from Gorenstein-AC-projective and Gorenstein-AC-injective modules:}

Gorenstein AC-projective and Gorenstein AC-injective modules where defined by D. Bravo, J. Gillespie and M. Hovey in \emph{\cite[Sections 5 and 8]{BHG}}. We recall these notions in the following lines. 
\begin{itemize}[itemsep=2pt,topsep=0pt]
\item[$\uno$] A left $R$-module $Q$ is \textbf{of type ${\rm FP}_\infty$} if it has an infinite presentation, that is, if there exists a long exact sequence 
\[
\cdots \to F_1 \to F_0 \to Q \to 0
\]
where $F_k$ is finitely generated and free, for every $k \geq 0$. 

\item[$\dos$] A left $R$-module $E$ is \textbf{absolutely clean} if $\Ext^1_R(Q,E) = 0$ for every left $R$-module $Q$ of type ${\rm FP}_\infty$.

\item[$\tres$] A left $R$-module $L$ is \textbf{level} if $\Tor^R_1(Q,L) = 0$ for every right $R$-module $Q$ of type ${\rm FP}_\infty$.

\item[$\cuatro$] A left $R$-module $M$ is said to be \textbf{Gorenstein AC-projective} if there exists an exact chain complex of projective modules
\[
\bm{P} = \cdots \to P_1 \to P_0 \to P^0 \to P^1 \to \cdots
\]
with $M = \Ker(P^0 \to P^1)$ and such that the induced chain complex $\Hom_R(\bm{P}, L)$ is exact for every level left $R$-module $L$.

\item[$\cinco$] A left $R$-module $N$ is said to be \textbf{Gorenstein AC-injective} if there exists an exact chain complex of injective modules
\[
\bm{I} = \cdots \to I_1 \to I_0 \to I^0 \to I^1 \to \cdots
\]
with $N = \Ker(I^0 \to I^1)$ and such that the induced chain complex $\Hom_R(E,\bm{I})$ is exact for every absolutely clean left $R$-module $E$. 
\end{itemize}

Now set $\nu := \Inj(R)$ and $\Y := \mathsf{GInj}_{\rm AC}(R)$ as the class of Gorenstein AC-injective modules. Then, note that $\Inj(R) \subseteq \mathsf{GInj}_{\rm AC}(R)$. On the other hand, $\mathsf{GInj}_{\rm AC}(R)$ is right thick since it is the right half of a hereditary cotorsion pair (see \emph{\cite[Lemma 5.6]{BHG}}). Now the fact that $\Inj(R)$ is a $\mathsf{GInj}_{\rm AC}(R)$-projective relative generator in $\mathsf{GInj}(R)$ follows by the definition of Gorenstein AC-injective modules. On the other hand, $\Ext^i_R(V,Y) = 0$ for every $V \in \Inj(R)$ and every $Y \in \mathsf{GInj}_{\rm AC}(R)$, since we can compute $\Ext^i_R(V,Y)$ by using an injective coresolution of $Y$ which is $\Hom_R(E,-)$-exact for every absolutely clean module $E$. Finally, we know that $\Inj(R)$ is closed under direct summands. Then, $(\Inj(R), \GInj_{\rm AC}(R))$ is a strong right Frobenius pair, and so by the duals of \emph{Theorems~\ref{CP2}} and \emph{\ref{theorem:projective_pair}} we obtain $\GInj^{< \infty}_{\rm AC}(R)$-cotorsion pairs $(\Inj(R)^\vee, \mathcal{GI}_{{\rm Ac}}(R))$ and $(\GInj^{< \infty}_{\rm AC}(R), \Inj(R))$ in $\mathsf{Mod}(R)$, which are right hereditary.   

In a similar way, from the definition of Gorenstein AC-projective modules and from the corresponding results in \emph{\cite{BHG}}, we have a strong left Frobenius pair $(\GProj_{\rm AC}(R), \Proj(R))$, where $\GProj_{\rm AC}(R)$ denotes the class of Gorenstein AC-projective modules. Moreover, we have left  hereditary $\GProj^{< \infty}_{\rm AC}(R)$-cotorsion pairs $(\GProj_{\rm AC}(R), \Proj^{< \infty}(R))$ and $(\Proj(R), \GProj^{< \infty}_{\rm AC}(R))$ in $\mathsf{Mod}(R)$.

%%%%%%%%%%%%%%%%%%%%%%%%%%%%%%%%%%%%%%%
%%%%%%%%%%%%%%%%%%%%%%%%%%%%%%%%%%%%%%%

\section{Exact model structures from Frobenius pairs}\label{sec:model_structures}

In this section, given a strong left Frobenius pair $(\X, \omega)$ in $\C$, we will obtain a model structure on the exact category $\X^\wedge$ (see \emph{Example~\ref{ex:thick_are_exact}}), whose homotopy category represents, in some sense, a generalization of the stable module category of a ring. 

\subsection{Hovey-Gillespie Correspondence}

We first recall the notion of exact model structures, along with the statement of \emph{Hovey-Gillespie Correspondence}, in order to obtain such model structures from the cotorsion pairs in \emph{Theorems~\ref{CP2}} and \emph{\ref{theorem:projective_pair}}. Our goal in this section is to show the following result.

\begin{theorem}\label{theo:Model_category_Frobenius} 
Let $(\X,\omega)$ be a strong left Frobenius pair in $\C$. Then there exists a unique exact model structure on $\X^\wedge$, denoted $\mathcal{M}^{\rm proj}_{\rm AB}(\X,\omega)$ and referred as the \textbf{projective Auslander-Buchweitz model structure}, such that $\X$ is the class of cofibrant objects, $\X^\wedge$ is the class of fibrant objects, and $\omega^\wedge$ is the class of trivial objects. 

Dually, if $(\nu,\Y)$ is a strong right Frobenius pair in $\C$, then there exists a unique exact model structure on $\Y^\vee$, denoted $\mathcal{M}^{\rm inj}_{\rm AB}(\nu,\Y)$ and referred as the \textbf{injective Auslander-Buchweitz model structure}, such that $\Y^\vee$ is the class of cofibrant objects, $\Y$ is the class of fibrant objects, and $\nu^\vee$ is the class of trivial objects.
\end{theorem}

The concept of model categories was introduced by D. Quillen in 1967 (see \emph{\cite{Quillen}}). Roughly speaking, a \textbf{model category} is given by a category $\mathcal{E}$ (normally assumed bicomplete) along with three classes of morphisms, called cofibrations, fibrations and weak equivalences, which provide a natural setting to do homotopy theory in $\mathcal{E}$. The triple of cofibrations, fibrations and weak equivalences is referred as a \textbf{model structure} on $\mathcal{E}$. By ``doing homotopy theory'', we mean that every model category has an associated \textbf{homotopy category}, denoted $\textrm{Ho}(\mathcal{E})$, obtained after localizing $\C$ at the class of weak equivalences. 

Let $\mathcal{E}$ be an exact category. From now on, if $\mathcal{E}$ is equipped with a model structure, we will consider the associated triple $(\mathcal{Q,R,T})$, where:
\begin{itemize}[itemsep=2pt,topsep=0pt]
\item[$\bullet$] $\mathcal{Q}$ is the class of \textbf{cofibrant objects} in $\mathcal{E}$, that is, objects $Q$ of $\mathcal{E}$ such that the only morphism $0 \to Q$ is a cofibration.

\item[$\bullet$] $\mathcal{R}$ is the class of \textbf{fibrant objects} in $\mathcal{E}$, that is, objects $R$ of $\mathcal{E}$ such that the only morphism $R \to 0$ is a fibration.

\item[$\bullet$] $\mathcal{T}$ is the class of \textbf{trivial objects} in $\mathcal{E}$, that is, objects $T$ of $\mathcal{E}$ such that the only morphism $0 \to T$ is a weak equivalence. 
\end{itemize}
We do not go into the details of model categories and homotopy theory, but we recommend M. Hovey's book \emph{\cite{HoveyBook}} as a very modern approach to these topics. Another recommended source with a nice homological approach is \emph{\cite[Chapter VIII]{BR}}, by A. Beligiannis and I. Reiten. 

In the case where $\mathcal{E}$ is an exact category, there is an appealing interplay between a special type of model structures on $\mathcal{E}$ and the concept of cotorsion pairs. These model structures are called \textit{exact}, and were defined by J. Gillespie in \emph{\cite{GiExact}}. 

A model structure on an exact category $\mathcal{E}$ is \textbf{exact} if:
\begin{itemize}[itemsep=2pt,topsep=0pt]
\item[$\bullet$] Cofibrations are given by the admissible monomorphisms with cokernel in $\mathcal{Q}$.

\item[$\bullet$] Fibrations are given by the admissible epimorphisms with kernel in $\mathcal{R}$.
\end{itemize}
On the other hand, given three classes $\F$, $\G$ and $\mathcal{W}$ of objects in $\mathcal{E}$, we say that $\F$, $\G$ and $\mathcal{W}$ form a \textbf{Hovey triple} $(\mathcal{F,G,W})$ is $\mathcal{W}$ is thick, and if $(\F \cap \mathcal{W},\G)$ and $(\F,\G \cap \mathcal{W})$ are complete cotorsion pairs in $\mathcal{E}$. Gillespie proved in \emph{\cite[Theorem 3.3]{GiExact}} the following one-to-one correspondence between exact model structures on $\mathcal{E}$ and Hovey triples in $\mathcal{E}$, in the case where the exact category $\mathcal{E}$ is \textbf{weakly idempotent complete}, that is, if every split monomorphism has a cokernel and every split epimorphism has a kernel. (See \emph{\cite[Definition 2.2]{GiExact}}).

\begin{theorem}[Hovey-Gillespie Correspondence]\label{theo:Hovey-Gillespie} 
Let $\mathcal{E}$ be an exact category with an exact model structure. Then $(\mathcal{Q,R,T})$ is a Hovey triple in $\mathcal{E}$. If in addition $\mathcal{E}$ is weakly idempotent complete, the converse also holds. That is, if $(\mathcal{F,G,W})$ is a Hovey triple in $\mathcal{E}$,  then there is a unique exact model structure on $\mathcal{E}$ such that $\F$, $\G$ and $\mathcal{W}$ are the classes of cofibrant, fibrant, and trivial objects, respectively.  
\end{theorem}

The previous theorem was originally proved by M. Hovey in \emph{\cite{Hovey}} in the case where $\mathcal{E}$ is an abelian category, where the bijective correspondence was established between \textbf{abelian} model structures and Hovey triples. Abelian model structures are defined in the same way as the exact ones: cofibrations are given by the monomorphisms with cokernel in $\mathcal{Q}$, and fibrations by the epimorphisms with kernel in $\mathcal{R}$. Notice that every abelian category $\C$ is exact, where the admissible monomorphisms and the admissible epimorphisms are given by the monomorphisms and epimorphisms in $\C$, respectively.

In this paper, we only use the converse of \emph{Hovey-Gillespie Correspondence}. Thick categories are examples of exact categories which are weakly idempotent complete, and so the converse of \emph{Hovey-Gillespie Correspondece} holds in the thick sub-categories $\X^\wedge, \Y^\vee \subseteq \C$ obtained from a left Frobenius pair $(\X,\omega)$ and a right Frobenius pair $(\nu,\Y)$. If in addition we assume that $(\X,\omega)$ is strong, then by \emph{Theorems~\ref{CP2}} and \emph{\ref{theorem:projective_pair}} we have two $\X^\wedge$-cotorsion pairs $(\X, \omega^\wedge)$ and $(\omega, \X^\wedge)$, where $\omega = \X \cap \omega^\wedge$. These are complete cotorsion pairs in the exact category $\X^\wedge$, and after setting $\F := \X$, $\G := \X^\wedge$, and $\mathcal{W} := \omega^\wedge$ in \emph{Theorem~\ref{theo:Hovey-Gillespie}}, we obtain the exact model structure on $\X^\wedge$ described in \emph{Theorem~\ref{theo:Model_category_Frobenius}}.

\begin{definition}\label{def:projective_model}
An exact model structure on an exact category $\mathcal{E}$ is \textbf{projective} if every object in $\mathcal{E}$ is fibrant. Dually, we have the concept of \textbf{injective} exact model structures. 
\end{definition}

\begin{example}
If $(\X,\omega)$ is a strong left Frobenius pair and $(\nu,\Y)$ a strong right Frobenius pair in $\C$, then the model structures $\mathcal{M}^{\rm proj}_{\rm AB}(\X,\omega)$ and $\mathcal{M}^{\rm inj}_{\rm Ab}(\nu,\Y)$ described in \emph{Theorem~\ref{theo:Model_category_Frobenius}} are projective and injective, respectively. 
\end{example}
 
\newpage

The following result provides an easy method to get exact model structures from one $\Sb$-cotorsion pair.

\begin{corollary}\label{coro:model_hereditary} 
Let $(\F,\G)$ be an $\Sb$-cotorsion pair in $\C$, where $\C$ has enough projectives and injectives, $\F$ is resolving in $\C$, and $\G$ coresolving in $\C$. If $\G \subseteq \F^\wedge$ then $(\F,\G)$ and $(\F \cap \G, \F^\wedge)$ are $\F^\wedge$-cotorsion pairs in $\C$. Dually, if $\F \subseteq \G^\vee$ then $(\F,\G)$ and $(\G^\vee, \F \cap \G)$ are $\G^\vee$-cotorsion pairs in $\C$.
\end{corollary}

We close this section presenting some examples of Auslander-Buchweitz model structures, which turn out to be generalizations of already known abelian model structures in Gorenstein and Ding-Chen Homological Algebra.

\begin{example}\label{ex:AB_model_structures} We generalize, in the sense explained below, some abelian model structures found by D. Bravo, M. Hovey and J. Gillespie. 
\begin{itemize}[itemsep=2pt,topsep=0pt]
\item[$\uno$] From the strong left Frobenius pair $(\GProj(R),\Proj(R))$ in $\mathsf{Mod}(R)$, we obtain the projective AB model structure 
\[
\mathcal{M}^{\rm proj}_{\rm AB}(\GProj(R), \Proj(R)) = (\GProj(R), \GProj^{< \infty}(R), \Proj^{< \infty}(R))
\] 
on $\GProj^{< \infty}(R)$. This is the only exact model structure on $\GProj^{< \infty}(R)$ with $\GProj(R)$, $\GProj^{< \infty}(R)$ and $\Proj^{< \infty}(R)$ as the classes of cofibrant, fibrant and trivial objects, respectively. Dually, there is a unique injective exact model structure on $\GInj^{< \infty}(R)$ such that $\GInj^{< \infty}(R)$, $\GInj(R)$ and $\Inj^{< \infty}(R)$ are the classes of cofibrant, fibrant, and trivial objects, respectively. 

These model structures generalize Hovey's projective and injective abelian model structures on $\mathsf{Mod}(R)$, with $R$ a Gorenstein ring (see \emph{\cite[Theorem 8.6]{Hovey}}). Notice we do not impose any condition on the ground ring $R$. However, we need to pay a price for this. We do not get abelian but exact model structures on $\GProj^{< \infty}(R)$ and $\GInj^{< \infty}(R)$, exact sub-categories of $\mathsf{Mod}(R)$. On the other hand, we have already mentioned that if $R$ is a Gorenstein ring, the previous two sub-categories coincide with $\mathsf{Mod}(R)$, and in this case, the model structures $\mathcal{M}^{\rm proj}_{\rm AB}(\GProj(R), \Proj(R))$ and $\mathcal{M}^{\rm inj}_{\rm AB}(\Inj(R), \GInj(R))$ are the abelian model structures described in \emph{\cite[Theorem 8.6]{Hovey}}.

\item[$\dos$] From the strong left and right Frobenius pairs 
\[
(\DProj(R), \Proj(R)) \mbox{ \ and \ } (\Inj(R), \GInj(R))
\] 
in $\mathsf{Mod}(R)$, we obtain the projective and injective AB model structures
\begin{align*}
\mathcal{M}^{\rm proj}_{\rm AB}(\DProj(R), \Proj(R)) & = (\DProj(R), \DProj^{< \infty}(R), \Proj^{< \infty}(R)), \\
\mathcal{M}^{\rm inj}_{\rm AB}(\Inj(R), \DInj(R)) & = (\DInj^{< \infty}(R), \DInj(R), \Inj^{< \infty}(R))
\end{align*}
on $\DProj^{< \infty}(R)$ and $\DInj^{< \infty}(R)$, respectively. 

These model structures are generalizations, in the sense specified in the previous example, of the abelian model structures on $\mathsf{Mod}(R)$ (with $R$ a Ding-Chen ring) found by J Gillespie in \emph{\cite[Theorem 4.7]{Gi}}. However,  the authors are not aware if $\mathsf{Mod}(R) = \DProj^{< \infty}(R)$ in the case where $R$ is a Ding-Chen ring.

\item[$\tres$]	Finally, in the more general setting of Gorenstein AC-projective and Gorenstein AC-injective modules, the strong left and right Frobenius pairs 
\[
(\GProj_{\rm AC}(R), \Proj(R)) \mbox{ \ and \ } (\Inj(R), \GInj_{\rm AC}(R))
\] 
give rise to the projective and injective AB model structures
\begin{align*}
\mathcal{M}^{\rm proj}_{\rm AB}(\GProj_{\rm AC}(R), \Proj(R)) & = (\GProj_{\rm AC}(R), \GProj^{< \infty}_{\rm AC}(R), \Proj^{< \infty}(R)), \\
\mathcal{M}^{\rm inj}_{\rm AB}(\Inj(R), \GInj_{\rm AC}(R)) & = (\GInj^{< \infty}_{\rm AC}(R), \GInj_{\rm AC}(R), \Inj^{< \infty}(R)).
\end{align*} 

These model structures are related to the abelian Gorenstein AC-projective and Gorenstein AC-injective model structures on $\mathsf{Mod}(R)$ (with $R$ an arbitrary ring) described in \emph{\cite[Theorems 5.5 and 8.5]{BHG}} by D. Bravo, M. Hovey and J. Gillespie. For those model structures, the classes of trivial objects have a description different to the one given above (see \emph{\cite[Lemmas 5.4 and 8.4]{BHG}}). On the other hand, the authors are not aware if the classes $\GProj^{< \infty}_{\rm AC}(R)$ and $\GInj^{< \infty}_{\rm AC}(R)$ coincide with the whole category $\mathsf{Mod}(R)$. If this turns out to be true, we would know another way to obtain the Gorenstein AC-projective and Gorenstein AC-injective modules structures.

\item[$\cuatro$] Examples $\mathsf{(1})$, $\dos$ and $\tres$ can also be obtained from \emph{Corollary~\ref{coro:model_hereditary}}. If $(\F,\G)$ is a $\Sb$-cotorsion pair in an abelian category with enough projectives and injectives, then:
	\begin{itemize}[itemsep=2pt,topsep=0pt]
	\item[$\ai$] If $\G \subseteq \F^\wedge$, then there exists a unique exact model structure on $\F^\wedge$ where $\F$ is the class of cofibrant objects, $\F^\wedge$ is the class of fibrant objects, and $\G$ is the class of trivial objects.
	
	\item[$\bi$] If $\F \subseteq \G^\vee$, then there exists a unique exact model structure on $\G^\vee$ where $\G^\vee$ is the class of cofibrant objects, $\G$ is the class of fibrant objects, and $\F$ is the class of trivial objects. ~%
	\end{itemize}
\end{itemize}
\end{example}

\subsection{Some remarks on sub-model structures}

We begin this section with the following definition.

\begin{definition}\label{def:hereditary_Hovey_triple}
We say that a Hovey triple $(\mathcal{F,G,W})$ in an exact category $\mathcal{E}$ is \textbf{left hereditary} if the cotorsion pairs $(\F, \G \cap \mathcal{W})$ and $(\F \cap \mathcal{W}, \G)$ are both left hereditary, that is, the classes $\F$ and $\F \cap \mathcal{W}$ are both resolving in $\mathcal{E}$. Dually, we have the definition of \textbf{right hereditary Hovey triples}. Finally, a Hovey triple is \textbf{hereditary} if it is both left and right hereditary. 

In the case $\mathcal{E}$ is an exact sub-category of an abelian category $\C$, the Hovey triple $(\mathcal{F,G,W})$ (in $\mathcal{E}$)  is said to be \textbf{left strong hereditary}, \textbf{right strong hereditary}, or \textbf{strong hereditray} if the cotorsion pairs $(\F, \G \cap \mathcal{W})$ and $(\F \cap \mathcal{W}, \G)$ are both left strong hereditary, right strong hereditary, or strong hereditary, respectively, in $\mathcal{E}.$ (See \emph{Definition~\ref{def:hereditary_relative}}). 
\end{definition}

In his paper \emph{\cite{GiExact}}, J. Gillespie constructs from a hereditary Hovey triple $(\mathcal{F,G,W})$ in $\mathcal{E}$, sub-model structures (of the unique model structure on $\mathcal{E}$ resulting from $(\mathcal{F,G,W})$) on the full sub-categories $\mathcal{Q}$, $\mathcal{R}$ and $\mathcal{Q} \cap \mathcal{R}$ of cofibrant, fibrant, and cofibrant-fibrant objects, respectively. The first property obtained after assuming that a Hovey triple $(\mathcal{F,G,W})$ is hereditary, is that the resulting sub-categories $\mathcal{Q}$, $\mathcal{R}$ and $\mathcal{T}$ are exact and weakly idempotent complete (see \emph{\cite[Lemma 5.1 and Proposition 5.2]{GiExact}}). The resulting sub-model structures on $\mathcal{Q}$, $\mathcal{R}$ and $\mathcal{Q} \cap \mathcal{R}$ are described in \emph{\cite[Proposition 5.2]{GiExact}}. In this section, we apply this result to check which are the sub-model structures obtained from a strong left Frobenius pair. 

The following result is a consequence of \emph{Theorem~\ref{theo:induced_pairs_hereditary}}.

\begin{corollary}
Let $(\X,\omega)$ be a strong left Frobenius pair in $\C$. Then $(\X,\X^\wedge,\omega^\wedge)$ is a strong hereditary Hovey triple in $\X^\wedge$ if, and only if, $\Proj(\C) \subseteq \omega$ and $\Inj(\C) \subseteq \X^\wedge$. 
\end{corollary}

Applying the previous corollary, along with \emph{\cite[Proposition 5.2]{GiExact}}, we obtain the following sub-model structures from a strong left Frobenius pair $(\X,\omega)$.

\begin{proposition}\label{prop:sub-model}
Let $(\X,\omega)$ be a strong left Frobenius pair in $\C$. If $\Proj(\C) \subseteq \omega$\footnote{If $\C$ has enough projectives, this inclusion can be replaced by $\Proj(\C) \subseteq \X^\wedge$, since \emph{Proposition~\ref{prop:same_projectives}} asserts $\Proj(\C) = \Proj(\X^\wedge)$, and so $\Proj(\C) \subseteq \omega$ (notice $\Proj(\X^\wedge) \subseteq \omega$ since $(\omega,\X^\wedge)$ is a cotorsion pair in $\X^\wedge$).} 
and $\Inj(\C) \subseteq \X^\wedge$, then there exist the following sub-model structures of $\mathcal{M}^{\rm proj}_{\rm AB}(\X, \omega)$:
\begin{itemize}[itemsep=2pt,topsep=0pt]
\item[$\ai$] An exact model structure on $\X$, where the cofibrant and fibrant objects are given by $\X$, and the trivial objects by $\omega$.

\item[$\bi$] An exact model structure on $\omega^\wedge$, where $\omega$ is the class of cofibrant objects, and $\omega^\wedge$ is the class of fibrant objects. 
\end{itemize}
\end{proposition}

In \emph{\cite[Proposition 5.2]{GiExact}}, the resulting sub-model structure on $\X^\wedge$ coincides with $\mathcal{M}^{\rm proj}_{\rm AB}(\X, \omega)$. This happens because $\mathcal{M}^{\rm proj}_{\rm AB}(\X, \omega)$ is a projective model structure.

\begin{example}
The model structure from $\ai$ is an example of what Gillespie calls a \textbf{Frobenius model structure}, that is, every object in the exact category is cofibrant and fibrant.

Recall that a Frobenius category is an exact category in which the notions of projective and injective objects coincide, and there are enough projectives and injectives. In every Frobenius category $\mathcal{E}$, $(\mathcal{E,E}, \Proj(\mathcal{E}))$ is a Hovey triple, which gives rise to a Frobenius model structure with $\Proj(\mathcal{E})$ as the class of trivial objects. 

One example of a Frobenius category is given by $\mathsf{Mod}(R)$, with $R$ a quasi-Frobenius ring. Another example is given by $\GProj(R)$ (with $R$ an arbitrary ring). Note that $\GProj(R)$ and $\mathsf{Mod}(R)$ coincide in the case where $R$ is quasi-Frobenius. It follows there is a unique exact Frobenius model structure on $\GProj(R)$ with $\Proj(R)$ as the class of trivial objects. This structure can also be obtained by setting $\C := \mathsf{Mod}(R)$, $\X := \GProj(R)$ and $\omega := \Proj(R)$ in \emph{Proposition~\ref{prop:sub-model}}.
\end{example}

\subsection{The homotopy category of an Auslander-Buchweitz model structure}

Let $(\X,\omega)$ be a strong left Frobenius pair in $\C$ and $\mathcal{M}^{\rm proj}_{\rm AB}(\X,\omega)$ be the unique projective AB model structure corresponding to $(\X,\omega)$. Let $\textrm{Ho}(\X^\wedge)$ be the homotopy category of $\mathcal{M}^{\rm proj}_{\rm AB}(\X,\omega)$. \newpage

In this section, we present an explicit description of $\textrm{Ho}(\X^\wedge)$, and notice how this homotopy category is equivalent to the stable category of a Frobenius sub-category of $\C$.  

Let us recall some results by J. Gillespie on homotopy theory for exact model structures. Given an exact model structure $(\mathcal{Q,R,T})$ on an exact category $\mathcal{E}$, two morphisms $f, g \colon X \to Y$ in $\mathcal{E}$ are \textbf{right homotopic} if, and only if, $g - f$ factors through an object in $\mathcal{Q} \cap \mathcal{T}$. Dually, $f$ and $g$ are \textbf{left homotopic} if, and only if, $g - f$ factors through an object in $\mathcal{R} \cap \mathcal{T}$. This was proven by Gillespie in \emph{\cite[Proposition 4.4]{GiExact}}. 

If the Hovey triple $(\mathcal{Q,R,T})$ is hereditary, then the sub-model structures on $\mathcal{Q}$, $\mathcal{R}$, and $\mathcal{Q} \cap \mathcal{R}$ are \textbf{full equivalent sub-model structures} of $\mathcal{E}$, that is, the inclusions $i_{\mathcal{Q}} \colon \mathcal{Q} \hookrightarrow \mathcal{E}$, $i_{\mathcal{R}} \colon \mathcal{R} \hookrightarrow \mathcal{E}$, and $i_{\mathcal{Q} \cap \mathcal{R}} \colon \mathcal{Q} \cap \mathcal{R} \hookrightarrow \mathcal{E}$ preserve the corresponding model structures, and the induced functors $\textrm{Ho}(i_{\mathcal{Q}}) \colon \textrm{Ho}(\mathcal{Q}) \to \textrm{Ho}(\mathcal{E})$, $\textrm{Ho}(i_{\mathcal{R}}) \colon \textrm{Ho}(\mathcal{R}) \to \textrm{Ho}(\mathcal{E})$, and $\textrm{Ho}(i_{\mathcal{Q} \cap \mathcal{R}}) \colon \textrm{Ho}(\mathcal{Q} \cap \mathcal{R}) \to \textrm{Ho}(\mathcal{E})$ are equivalence of categories. This is proven in \emph{\cite[Proposition 5.2 and Corollary 5.4]{GiExact}}. 

The Hovey triple associated to the projective AB model structure $\mathcal{M}^{\rm proj}_{\rm AB}(\X, \omega)$ is hereditary, as proven below.

\begin{proposition}\label{prop:AB_Hovey_triples_hereditaries}
Let $(\X,\omega)$ be a left Frobenius pair in $\C$. Then the class $\omega^\wedge$ is right thick. If in addition $(\X,\omega)$ is strong, then $(\X,\X^\wedge,\omega^\wedge)$ is a hereditary Hovey triple in $\X^\wedge$. 
\end{proposition}

\begin{myproof}
By \emph{Theorem~\ref{CP2}}, we have the equality $\omega^\wedge = \X^\perp \cap \X^\wedge$. The classes $\X^\perp$ and $\X^\wedge$ are clearly right thick, and hence $\omega^\wedge$ is right thick. 

Now suppose $(\X,\omega)$ is strong. Then we have two cotorsion pairs $(\X,\omega^\wedge)$ and $(\omega,\X^\wedge)$ in the exact sub-category $\X^\wedge \subseteq \C$. First, note that $\omega$ is the class of projective objects in $\X^\wedge$. It follows that $\X$ and $\omega$ are resolving classes in $\X^\wedge$. On the other hand, the classes $\omega^\wedge$ and $\X^\wedge$ are right thick, and so they are coresolving in $\X^\wedge$. 
\end{myproof}

It follows by the previous proposition and the comments above that, if $\X^\wedge$ is equipped with the projective AB model structure $\mathcal{M}^{\rm proj}_{\rm AB}(\X,\omega)$, then $\mathcal{Q} \cap \mathcal{R} := \X$ is equipped with the sub-model structure $(\X,\X,\omega)$ of $\mathcal{M}^{\rm proj}_{\rm AB}(\X,\omega)$, and the corresponding homotopy categories $\textrm{Ho}(\X^\wedge)$ and $\textrm{Ho}(\X)$ are equivalent. The homotopy category $\textrm{Ho}(\X)$ is easier to describe than $\textrm{Ho}(\X^\wedge)$. First of all, on the sub-category $\mathcal{Q} \cap \mathcal{R} := \X \subseteq \X^\wedge$, the relations of being left and right homotopic coincide (denoted $\sim$). On the other hand, by \emph{\cite[Theorem 1.2.10 $\mathsf{(i)}$]{HoveyBook}}, there is an equivalence of categories: 
\[
\textrm{Ho}(\X) \simeq ({\mathcal{Q} \cap \mathcal{R}) / \sim} = \X / \sim.
\]
By \emph{Remark~\ref{rem:Frobenius_category}}, $\X$ is a Frobenius category, and for any two morphisms $f, g \colon X \to Y$, $f \sim g$ if, and only if, $g - f$ factors through a projective object in $\X$ (that is, an object in $\omega$). Hence, the quotient $\X / \sim$ is the \textbf{stable category} of $\X$. Now with respect to $\textrm{Ho}(\X^\wedge)$, we have by \emph{\cite[Theorem 1.2.10 (ii)]{HoveyBook}} that for every $X, Y \in \X^\wedge$ there is a natural isomorphism
\begin{align*}
\Hom_{\textrm{Ho}(\X^\wedge)}(X,Y) & \simeq \Hom_{\X^\wedge}(QX, RY) / \sim,
\end{align*}
where $QX$ denotes the cofibrant replacement of $X$, and $RY$ the fibrant replacement of $Y$. Since every object in $\X^\wedge$ is fibrant, we have $RY = Y$. We summarize these facts in the following result.

\begin{proposition}\label{prop:homotopy_AB_structures}
Let $(\X,\omega)$ be a strong left Frobenius pair in $\C$. For the projective AB model structure $\mathcal{M}^{\rm proj}_{\rm AB}(\X,\omega) = (\X,\X^\wedge,\omega^\wedge)$ on $\X^\wedge$, there is a natural isomorphism
\[
\Hom_{{\rm Ho}(\X^\wedge)}(X,Y) \cong \Hom_{\X^\wedge}(QX, Y) / \sim
\]
for every $X, Y \in \X^\wedge$, where $f \sim g$ if, and only if, $g - f$ factors through an object in $\omega$. Moreover, $\textrm{Ho}(\X^\wedge)$ is equivalent to the stable category $\X / \sim$. 

Dually, if $(\nu, \Y)$ is a strong right Frobenius pair in $\C$, then for the injective AB model structure $\mathcal{M}^{\rm inj}_{\rm AB}(\nu, \Y) = (\Y^\vee, \Y, \nu^\vee)$ on $\Y^\vee$ there is a natural isomorphism
\[
\Hom_{{\rm Ho}(\Y^\vee)}(X,Y) \cong \Hom_{\Y^\vee}(X, RY) / \sim
\]
for every $X, Y \in \Y^\vee$, where $f \sim g$ if, and only if, $g - f$ factors through an object in $\nu$. Moreover, $\textrm{Ho}(\Y^\vee)$ is equivalent to the stable category $\Y / \sim$. 
\end{proposition}

The previous proposition encodes the stable module category $\mathsf{Stmod}(R)$ of a quasi-Frobenius ring $R$ as a particular example, by setting $\C := \mathsf{Mod}(R)$ (with $R$ a quasi-Frobenius ring), $\X := \GProj(R) = \mathsf{Mod}(R)$, and $\omega := \Proj(R)$. This is explained in more detail below.

\begin{example} 
We recover some examples of homotopy categories previously constructed by D. Bravo, M. Hovey and J. Gillespie. 
\begin{itemize}[itemsep=2pt,topsep=0pt]
\item[$\uno$] Consider the homotopy category $\mathsf{Ho}(\GProj^{< \infty}(R))$ of the projective AB model structure 
$$\mathcal{M}^{\rm proj}_{\rm AB}(\GProj(R), \Proj(R)).$$ By \emph{Proposition~\ref{prop:homotopy_AB_structures}}, we have that two morphisms $f, g \colon X \to Y$ in $\GProj^{< \infty}(R)$ are homotopic if their difference $g - f$ factors through a projective module. The homotopy category of this model structure is the projective stable module category $\GProj(R) / \sim$, which is also the homotopy category of the Hovey's projective abelian model structure $(\GProj(R), \mathsf{Mod}(R), \Proj^{< \infty}(R))$ on $\mathsf{Mod}(R)$, when $R$ is a Gorenstein ring (see \emph{\cite[Section 9]{Hovey}}). This stable module category coincides with the usual stable module category $\mathsf{Stmod}(R)$ in the case where $R$ is a quasi-Frobenius ring, that is, a $0$-Gorenstein ring. Using  \emph{Proposition~\ref{prop:homotopy_AB_structures}}, we can obtain similar conclusions in the injective case, and the same stable module category from $\mathsf{Ho}(\GInj^{< \infty}(R))$ in the case $R$ is a quasi-Frobenius ring.

\item[$\dos$] Recalling \emph{Example~\ref{ex:AB_model_structures}} $\tres$, the homotopy categories of the AB  model structures $$\mathcal{M}^{\rm proj}_{\rm AB}(\GProj_{\rm AC}(R), \Proj(R))\quad\text{ and }\quad\mathcal{M}^{\rm inj}_{\rm AB}(\Inj(R), \GInj_{\rm AC}(R))$$  are exactly the homotopy categories obtained in \emph{\cite[Theorems 5.7 and 8.7]{BHG}}. 
\end{itemize}
\end{example}

%%%%%%%%%%%%%%%%%%%%%%%%%%%%%%%%%%%%%%%%%%%%%%%%%%%%%%%%%%%%%%%
%%%%%%%%%%%%%%%%%%%%%%%%%%%%%%%%%%%%%%%%%%%%%%%%%%%%%%%%%%%%%%%

\section{Auslander-Buchweitz contexts: correspondences with relative cotorsion pairs, Frobenius pairs and model structures}\label{sec:AB_contexts}

We conclude this paper presenting one-to-one correspondences between the objects we have been studying so far: relative cotorsion pairs and AB model structures. Auslander-Buchweitz contexts will play an important role in this section, and will also appear in this correspondence.

\subsection{AB contexts vs. Frobenius pairs vs. relative cotorsion pairs}

The following definition of Auslander-Buchweitz contexts is due to \emph{\cite[Theorem 1.1.2.10]{Hashimoto}}, but it is written according to the terminology we have been using so far.

\begin{definition}\label{def:AB-context}
Let $(\A,\B)$ be a pair of classes of objects in $\C$ and $\omega := \A \cap \B$. We say that $(\A,\B)$ is a \textbf{left weak Auslander-Buchweitz pre-context} (\textbf{left weak AB pre-context} for short) in $\C$ if:
\begin{itemize}[itemsep=2pt,topsep=0pt]
\item[$\ai$] The pair $(\A, \omega)$ is a left Frobenius pair.

\item[$\bi$] $\B = \Thick^{+}(\B)$.  
\end{itemize}
If in addition:
\begin{itemize}[itemsep=2pt,topsep=0pt] 
\item[$\bullet$] $(\A,\B)$ satisfies $\B \subseteq \A^\wedge$, then we say that $(\A,\B)$ is a \textbf{left weak AB context} in $\C$; 

\item[$\bullet$] $(\A,\B)$ satisfies $\A^\wedge = \C$, then we say that $(\A,\B)$ is a \textbf{left AB context} in $\C$.
\end{itemize}

The notion of \text{right weak AB pre-context} is defined dually, that is, a pair $(\A,\B)$ of classes of objects in $\C$, with $\nu := \A \cap \B$, such that:
\begin{itemize}[itemsep=2pt,topsep=0pt]
\item[$\iroman$] The pair $(\nu,\B)$ is a right Frobenius pair.

\item[$\iiroman$] $\A = \Thick^{-}(\A)$.
\end{itemize} 
If in addition:
\begin{itemize}[itemsep=2pt,topsep=0pt] 
\item[$\bullet$] $(\A,\B)$ satisfies $\A \subseteq \B^\vee$, then we say that $(\A,\B)$ is a \textbf{right weak AB context} in $\C$; 

\item[$\bullet$] $(\A,\B)$ satisfies $\B^\vee = \C$, then we say that $(\A,\B)$ is a \textbf{right AB context} in $\C$.
\end{itemize}
\end{definition}

\begin{example}\label{ex:ABponja}
The classes $\X$ and $\Y$ in \emph{Theorem~\ref{theo:ponja}} form a left weak AB context $(\X,\Y)$. 
\end{example}

\begin{theorem}\label{theo:from_AB-context_to_hereditary_pairs}
Let $(\A,\B)$ be a left weak AB context in $\C$, and $\omega := \A \cap \B$. Then:
\begin{itemize}[itemsep=2pt,topsep=0pt]
\item[$\ai$] $\omega = \A \cap \A^\perp$ and $\omega^\wedge = \B$.

\item[$\bi$] $(\A,\B)$ is a $\A^\wedge$-cotorsion pair in $\C$ with $\id_{\A}(\B) = 0$.
\end{itemize}
\end{theorem}

\begin{myproof} 
On the one hand, $(\A,\omega)$ is a left Frobenius pair in $\C$, and by \emph{Theorem~\ref{CP2}} we have an $\A^\wedge$-cotorsion pair $(\A, \omega^\wedge)$ with $\omega = \A \cap \omega^\wedge = \A \cap (\A^\perp \cap \A^\wedge) = \A \cap \A^\perp$. On the other hand, $\omega \subseteq \B$ and $\B = \Thick^{+}(\B)$ imply $\omega^\wedge \subseteq \B$. We now show the remaining inclusion. Let $N \in \B$. Then $N \in \A^\wedge$ since $\B \subseteq \A^\wedge$, and so we can apply \emph{Theorem~\ref{AB4}} in order to obtain a short exact sequence
\[
0 \to K \to A \to N \to 0
\]
with $A \in \A$ and $K \in \omega^\wedge$. Since $\omega^\wedge \subseteq \B$, we have that $K \in \B$. Then $A \in \A \cap \B =: \omega$ since $\B$ is closed under extensions. It follows $N \in \omega^\wedge$. Hence, $\B = \omega^\wedge$, and so $(\A,\B)$ is an $\A^\wedge$-cotorsion pair in $\C$. 

It is only left to show that $\id_{\A}(\B) = 0$. Let  $B \in \B$. Consider a short exact sequence
\[
0 \to K \to A \to B \to 0
\]
with $A \in \A$ and $K \in \omega^\wedge = \B$. By the dual of \emph{\cite[Lemma 1.1]{AM}}, we have
\[
\id_{\A}(B) \leq \max\{ \id_{\A}(A), \id_{\A}(K) - 1 \},
\]
where $\id_{\A}(K) \leq \id_{\A}(\omega^\wedge)$, and $\id_{\A}(\omega^\wedge) = \id_{\A}(\omega)$ by \emph{Lemma~\ref{AB1}}. On the one hand, $\id_{\A}(\omega) = 0$. Then $\id_{\A}(K) = 0$, and hence 
\[
\id_{\A}(B) \leq \id_{\A}(A).
\]
On the other hand, $A \in \A \cap \B = \omega$, and so $\id_{\A}(A) = 0$. Therefore, $\id_{\A}(B) = 0$ for every $B \in \B$. 
\end{myproof}

\begin{lemma}\label{lem:relative_pair_to_A_thick}
Let $(\F,\G)$ be a left $\Sb$-cotorsion pair in $\C$ such that $\id_{\F}(\G) = 0$. Then $\F = \Thick^{-}(\F)$. Dually, if $(\F,\G)$ is a right $\Sb$-cotorsion pair in $\C$ such that $\pd_{\G}(\F) = 0$, then $\G = \Thick^{+}(\G)$. 
\end{lemma}

\begin{myproof}
First, we know that $\F = {}^{\perp_{1,\Sb}}\G := {}^{\perp_1}\G \cap \Sb$. Thus, we have that $\F$ is closed under direct summands and extensions. It suffices to show that $\F$ is closed under kernels of epimorphisms in $\F$. So suppose we are given a short exact sequence
\[
0 \to A \to B \to C \to 0
\]
with $B, C \in \F$. Let $G \in \G$. Then we have an exact sequence 
\[
\Ext^1_{\C}(B,G) \to \Ext^1_{\C}(A,G) \to \Ext^2_{\C}(C,G)
\]
where $\Ext^1_{\C}(B,G) = 0$ and $\Ext^2_{\C}(C,G) = 0$, since $\id_{\F}(\G) = 0$. Hence $\Ext^1_{\C}(A,G) = 0$ for every $G \in \G$, and so $A \in {}^{\perp_1}\G$. On the other hand, $A \in \Sb$ since $\Sb$ is thick and $B, C \in \F \subseteq \Sb$. Therefore, $A \in {}^{\perp_1}\G \cap \Sb =: {}^{\perp_{1,\Sb}}\G = \F$. 
\end{myproof}

\begin{proposition}\label{prop:cotorsion_to_precontext}
Let $(\F,\G)$ be an $\Sb$-cotorsion pair in $\C$ with $\id_{\F}(\G) = 0$. Then $(\F,\G)$ is a left and right weak AB pre-context in $\C$. 
\end{proposition}

\begin{myproof}
The equalities $\F = \Thick^{-}(\F)$ and $\G = \Thick^{+}(\G)$ follow by \emph{Lemma~\ref{lem:relative_pair_to_A_thick}} and its dual. On the other hand, $\id_{\F}(\G) = 0$ implies $\id_{\F}(\omega) = 0$, where $\omega := \F \cap \G$. Now let $F \in \F \subseteq \Sb$. Then there exists a short exact sequence 
\[
0 \to F \to G \to F' \to 0
\]
with $F' \in \F$ and $G \in \G$. Since $\F$ is closed under extensions, we get $G \in \F \cap \G =: \omega$, and so $\omega$ is a relative cogenerator in $\F$. Finally, $\pd_{\G}(\F) = \id_{\F}(\G) = 0$ and $\omega \subseteq \F$ imply $\pd_{\G}(\omega) = 0$. Dually, for every $G \in \G$ there exists a short exact sequence
\[
0 \to G' \to F \to G \to 0
\]
with $F \in \F$ and $G' \in \G$, and so $F \in \F \cap \G =: \omega$. Therefore, $\omega$ is a relative generator in $\G$. 
\end{myproof}

\begin{theorem}\label{theo:ThickF-cotorsion_pair}
Let $(\F,\G)$ be a $\Thick(\F)$-cotorsion pair in $\C$ with $\id_{\F}(\G) = 0$. Then:
\begin{itemize}[itemsep=2pt,topsep=0pt]
\item[$\ai$] $\Thick(\F) = \F^\wedge$.

\item[$\bi$] $(\F,\G)$ is a left weak AB context in $\C$.
\end{itemize}
\end{theorem}

\begin{myproof} 
 By \emph{Proposition~\ref{prop:cotorsion_to_precontext}}, $(\F,\G)$ is a left weak AB pre-context. Then $(\F,\omega)$ is a left Frobenius pair, and so  \emph{Theorem~\ref{AB9}} implies $\F^\wedge = \Thick(\F)$. Finally, we only need to show that $\G \subseteq \F^\wedge$. But $\G \subseteq \Thick(\F)=\F^\wedge$ and thus the result follows. 
\end{myproof}

From now on, $\C^2$ will denote the product category $\C \times \C$.

\begin{theorem}\label{theo:correspondence_with_AB-context}
Let $\C$ be an abelian category. For the classes
\begin{align*}
\mathfrak{F} & := \{ (\X, \omega) \subseteq \C^2 \mbox{ \emph{: $(\X,\omega)$ is a left Frobenius pair in $\C$}}\}, \\
\mathfrak{C} & := \{ (\A,\B) \subseteq \C^2 \mbox{ \emph{: $(\A,\B)$ is a left weak AB context in $\C$}} \}, \\
\mathfrak{P} & := \{ (\F,\G) \subseteq \C^2 \mbox{ \emph{: $(\F,\G)$ is a $\Thick(\F)$-cotorsion pair in $\C$ with $\id_{\F}(\G) = 0$}} \},
\end{align*}
\sloppypar{the equality $\mathfrak{C} = \mathfrak{P}$ holds. Moreover, there exists a one-to-one correspondence between the classes \mbox{$\mathfrak{F}$ and $\mathfrak{C}$.}}

Dually, for the classes
\emph{
\begin{align*}
\mathfrak{F}\op & := \{ (\nu, \Y) \subseteq \C^2 \mbox{ : $(\nu,\Y)$ is a right Frobenius pair in $\C$}\}, \\
\mathfrak{C}\op & := \{ (\A,\B) \subseteq \C^2 \mbox{ : $(\A,\B)$ is a right weak AB context in $\C$} \}, \\
\mathfrak{P}\op & := \{ (\F,\G) \subseteq \C^2 \mbox{ : $(\F,\G)$ is a $\Thick(\G)$-cotorsion pair in $\C$ with $\pd_{\G}(\F) = 0$} \},
\end{align*}
}
\sloppypar{the equality $\mathfrak{C}\op = \mathfrak{P}\op$ holds, and there exists a one-to-one correspondence between the classes \mbox{$\mathfrak{F}\op$ and $\mathfrak{C}\op$.}}
\end{theorem}

\begin{myproof}
Define the mapping
\begin{align*}
\Phi \colon \mathfrak{F} & \too \mathfrak{C} \\
(\X,\omega) & \mapsto (\A := \X, \B := \omega^\wedge).
\end{align*}
We first show that $\Phi$ is well defined, that is, that if $(\X,\omega)$ is a left Frobenius pair, then $(\A := \X, \B := \omega^\wedge)$ is a left weak AB context. By \emph{Proposition~\ref{AB2}} we have 
\[
(\A, \A \cap \B) = (\X, \X \cap \omega^\wedge) = (\X,\omega),
\]
and so $(\A, \A \cap \B)$ is a left Frobenius pair. Then it is only left to show that $\B = \Thick^{+}(\B)$ and that $\B \subseteq \A^\wedge$. By \emph{Theorem~\ref{CP2}}, $\B = \omega^\wedge = \X^\perp \cap \X^\wedge$, and also $\X^\perp$ and $\X^\wedge$ are right thick. Then, $\B$ is right thick and hence $\B = \Thick^{+}(\B)$. Finally, the inclusion $\B = \omega^\wedge \subseteq \X^\wedge = \A^\wedge$ is clear. 

We have that $\Phi$ is a well defined mapping. To show that it is one-to-one, we construct an inverse for it. Consider the mapping
\begin{align*}
\Psi \colon \mathfrak{C} & \too \mathfrak{F} \\
(\A,\B) & \mapsto (\X := \A, \omega := \A \cap \B)
\end{align*}
which is well defined by the definition of left weak AB context. We show $\Psi \circ \Phi = \id_{\mathfrak{F}}$ and $\Phi \circ \Psi = \id_{\mathfrak{C}}$.

Let $(\X,\omega)$ be a left Frobenius pair. Then
\[
\Psi \circ \Phi(\X,\omega) = \Psi(\X, \omega^\wedge) = (\X, \X \cap \omega^\wedge) = (\X, \omega)
\]
where $\X \cap \omega^\wedge = \omega$ by \emph{Proposition~\ref{AB2}}. On the other hand, if $(\A,\B)$ is a left weak AB context, we have
\[
\Phi \circ \Psi(\A,\B) = \Phi(\A,\A \cap \B) = (\A, (\A \cap \B)^\wedge) = (\A, \B)
\]
where $(\A \cap \B)^\wedge = \B$ by \emph{Theorem~\ref{theo:from_AB-context_to_hereditary_pairs}}. 

Now we focus on showing $\mathfrak{C} = \mathfrak{P}$. Let $(\A,\B) \in \mathfrak{C}$ be a left weak AB context. Then $(\A,\B)$ is a $\A^\wedge$-cotorsion pair with $\id_{\A}(\B) = 0$, by \emph{Theorem~\ref{theo:from_AB-context_to_hereditary_pairs}}, where $\A^\wedge = \Thick(\A)$ by \emph{Theorem~\ref{AB9}}. Hence $(\A,\B) \in \mathfrak{P}$, and $\mathfrak{C} \subseteq \mathfrak{P}$. The remaining inclusion $\mathfrak{C} \supseteq \mathfrak{P}$ follows by \emph{Theorem~\ref{theo:ThickF-cotorsion_pair}}.
\end{myproof}

\begin{figure}[H]
\centering 
\scriptsize
\begin{tikzpicture}[description/.style={fill=white,inner sep=2pt}]
\matrix (m) [matrix of math nodes, row sep=0.5em, column sep=0em, text height=3.5ex, text depth=1.5ex]
{ {} & {} & {} & {} & {\color{red}{(\A,\A \cap \B)}}& {} & {} & {} & {} & {} \\ 
  {} & {} & {} & {\color{white}{(\X,\omega)}} & {\Large\bm{\mathfrak{F}}} & {\color{red}{(\X,\omega)}} & {} & {} & {} \\
  {} & {} & {} & {} & {} & {} & {} & {} & {} & {} \\
  {\color{red}{(\A,\B)}} & {\Large\bm{\mathfrak{C}}} & {} & {} & {} & {} & {} & {\Large\bm{\mathfrak{P}}} & {\color{red}{(\X,\omega^\wedge)}} \\ };
\path[red,|->]
($(m-4-1.east)+(-0.25,0.4)$) edge node[above,sloped] {$\Psi$} ($(m-1-5.west)-(-0.45,+0.4)$)
($(m-2-6.east)+(-0.5,-0.4)$) edge node[above,sloped] {$\Phi$} ($(m-4-8.west)-(-0.85,-0.2)$)
;
\path[->] 
(m-4-2) edge [thick] node[above,sloped] {$\Psi$} (m-2-5)
(m-2-5) edge [thick] node[above,sloped] {$\Phi$} (m-4-8)
;
\path[-,font=\scriptsize]
(m-4-2) edge [double, thick, double distance=2pt] (m-4-8)
;
\end{tikzpicture}
\caption[Correspondences between left Frobenius pairs, left weak AB contexts, and relative cotorsion pairs in abelian categories]{Correspondences between left Frobenius pairs, left weak AB contexts, and relative cotorsion pairs in abelian categories.}
\label{fig:correspondence_abelian}
\end{figure}

\begin{figure}[H]
\centering 
\scriptsize
\begin{tikzpicture}[description/.style={fill=white,inner sep=2pt}]
\matrix (m) [matrix of math nodes, row sep=0.5em, column sep=0em, text height=3.5ex, text depth=1.5ex]
{ {} & {} & {} & {} & {\color{blue}{(\A \cap \B,\B)}}& {} & {} & {} & {} & {} \\ 
  {} & {} & {} & {\color{white}{(\nu,\Y)}} & {\Large\bm{\mathfrak{F}\op}} & {\color{blue}{(\nu,\Y)}} & {} & {} & {} \\
  {} & {} & {} & {} & {} & {} & {} & {} & {} & {} \\
  {\color{blue}{(\A,\B)}} & {\Large\bm{\mathfrak{C}\op}} & {} & {} & {} & {} & {} & {\Large\bm{\mathfrak{P}\op}} & {\color{blue}{(\nu^\vee,\Y)}} \\ };
\path[blue,|->]
($(m-4-1.east)+(-0.2,0.3)$) edge node[above,sloped] {$\Psi\op$} ($(m-1-5.west)-(-0.45,+0.4)$)
($(m-2-6.east)+(-0.25,-0.4)$) edge node[above,sloped] {$\Phi\op$} ($(m-4-8.west)-(-1.2,-0.2)$)
;
\path[->] 
(m-4-2) edge [thick] node[above,sloped] {$\Psi\op$} (m-2-5)
(m-2-5) edge [thick] node[above,sloped] {$\Phi\op$} (m-4-8)
;
\path[-,font=\scriptsize]
(m-4-2) edge [double, thick, double distance=2pt] (m-4-8)
;
\end{tikzpicture}
\caption[Correspondences between right Frobenius pairs, right weak AB contexts, and relative cotorsion pairs in abelian categories]{Correspondences between right Frobenius pairs, right weak AB contexts, and relative cotorsion pairs in abelian categories.}
\label{fig:correspondence_abelian_dual}
\end{figure}

%%%%%%%%%%%%%%%%%%%%%%%%%%%%%%%%%%%
%%%%%%%%%%%%%%%%%%%%%%%%%%%%%%%%%%%

\subsection{Relative cotorsion pairs vs. covering classes}

The class $\mathfrak{G}$ defined by
\[
\mathfrak{G} := \left\{ \F \subseteq \C \mbox{ : $\F$ is a left saturated class in $\C$, and pre-covering in $\Thick(\F)$} \right\}
\]
plays a role in \emph{Theorem~\ref{theo:correspondence_with_AB-context}} when we impose some extra conditions on the category $\C$. In order to be more specific about this, we need the following definition.

\begin{definition}\label{def:perfect}
Let $\Sb$ be a thick sub-category of $\C$. A left $\Sb$-cotorsion pair $(\F,\G)$ is said to be \textbf{perfect} if every object in $\Sb$ has a $\F$-cover. \textbf{Perfect right $\Sb$-cotorsion pairs} are defined dually. Finally, an $\Sb$-cotorsion pair is \textbf{perfect} if every object in $\Sb$ has an $\F$-cover and a $\G$-envelope. 
\end{definition}

It is possible to study some interplay between the classes $\mathfrak{G}$ and:
\[
\mathfrak{P}' := \left\{ (\F,\G) \subseteq \C^2 \mbox{ : } \begin{array}{ll} \mbox{$(\F,\G)$ is a perfect left $\Thick(\F)$-cotorsion pair in $\C$,} \\ \mbox{with $\id_{\F}(\G) = 0$ and $\Proj(\C) \subseteq \F$} \end{array} \right\}. 
\]
We first see how to map elements in $\mathfrak{P}'$ to $\mathfrak{G}$. Let $\Theta$ be the following mapping:
\begin{align*}
\Theta \colon \mathfrak{P}' & \to \mathfrak{G} \\
(\F,\G) & \mapsto \F.
\end{align*}
We assert that $\Theta$ is well defined. This is a consequence of the following result.

\begin{proposition}\label{prop:correspondence_with_G}
Let $(\F,\G)$ be a left $\Thick(\F)$-cotorsion pair in $\C$ such that $\id_{\F}(\G) = 0$. Then $\F$ is left thick and a pre-covering class in $\Thick(\F)$.
\end{proposition}

\begin{myproof}
By \emph{$\mathsf{(scp4)}$} in \emph{Proposition~\ref{CP1}}, every object in $\Thick(\F)$ has an epic $\F$-pre-cover. It is only left to show that $\F = \Thick^{-}(\F)$,  but this follows by \emph{Lemma~\ref{lem:relative_pair_to_A_thick}}. 
\end{myproof}

It follows for every $(\F,\G) \in \mathfrak{P}'$, that $\F$ is left thick and a pre-covering class in $\Thick(\F)$. This, plus the fact that $\Proj(\C) \subseteq \F$, imply that $\F$ is left saturated in $\C$, that is, $\Theta$ is well defined.

Under certain conditions, it is possible to construct a left inverse for the previous mapping. Namely, we need the ground category $\C$ to be Krull-Schmidt. Recall that a category $\C$ is \textbf{Krull-Schmidt} if it is additive and if every object decomposes into a finite direct sum of objects having local endomorphism rings.

\begin{theorem}\label{theo:correspondence_with_G}
Let $\C$ be an abelian Krull-Schmidt category with enough projectives. Then the mapping $\Theta \colon \mathfrak{P}' \to \mathfrak{G}$ defines an injection.
\end{theorem}

Before presenting a proof for the previous theorem, we need the following lemma.

\begin{lemma}\label{lem:correspondence_with_G}
Let $\C$ be a Krull-Schmidt abelian category with enough projectives, and $\Sb$ be a thick sub-category of $\C$. If $\F \subseteq \Sb$ is a left saturated class in $\C$ which is also pre-covering in $\Sb$, then the following conditions hold:
\begin{itemize}[itemsep=2pt,topsep=0pt]
\item[$\ai$] For each $S \in \Sb$, there exists an exact sequence
\[
0 \to \Ker(\varphi) \to F \xrightarrow{\varphi} S \to 0
\]
where $\varphi$ is an $\F$-cover, and $\Ker(\varphi) \in \F^{\perp_1}$. 

\item[$\bi$] $\F^\perp = \F^{\perp_1}$.
 
\item[$\ci$] Setting $\G := \F^\perp \cap \Sb$, the pair $(\F,\G)$ is a perfect left $\Sb$-cotorsion pair in $\C$, with $\id_{\F}(\G) = 0$. 
\end{itemize}

Dually, if $\C$ is a Krull-Schmidt abelian category with enough injectives, and $\G \subseteq \Sb$ is a right saturated class in $\C$ which is also pre-enveloping in $\Sb$, then $(\F,\G)$ is a perfect right $\Sb$-cotorsion pair in $\C$, with $\pd_{\G}(\F) = 0$, where $\F := {}^\perp \G \cap \Sb$ and ${}^\perp \G = {}^{\perp_1} \G$.
\end{lemma}

\begin{myproof} 
We only prove $\ai$, $\bi$, and $\ci$. 
\begin{itemize}[itemsep=2pt,topsep=0pt]
\item[$\ai$] Let $S \in \Sb$. On the one hand, $S$ has a $\F$-pre-cover. Since $\C$ is a Krull-Schmidt category, the endomorphism ring ${\rm End}_{\C}(S)$ is semiperfect (see \emph{\cite[Corollary 4.4]{Krause}}). On the other hand, by \emph{\cite[Corollary 2.1.10 (b)]{GT}}\footnote{The statement and proof in the given reference is for the category of left $R$-modules. However, this result is also valid in abelian categories.}, $S$ has an $\F$-cover, which can be taken epic since $\Proj(\C) \subseteq \F$. It follows there exists a short exact sequence
\[
0 \to \Ker(\varphi) \to F \xrightarrow{\varphi} S \to 0
\]
where $\varphi \colon \F \to S$ is an $\F$-cover. Then, using the \emph{Wakamatsu Lemma} (see \emph{\cite[Lemma 5.12]{GT}}\footnote{Although the statement and proof appearing in the given reference are written for the category of right $R$-modules, the arguments carry over to any abelian category.}), we have that $\Ker(\varphi) \in \F^{\perp_1}$.

\item[$\bi$] It suffices to show $\F^{\perp_1} \subseteq \F^\perp$. So let $Y \in \F^{\perp_1}$ and $F \in \F$.  Since $\C$ has enough projectives, for each $i \geq 0$
 we can find an exact sequence 
\[
0 \to K_i \to P_{i-2} \to \cdots \to P_1 \to P_0 \to F \to 0
\]
where $K_i \in \Omega^{i-1}(F)$. Note that $K_i \in \F$ since $\F$ is left saturated in $\C$. Then by \emph{Dimension Shifting}, we have the equality 
\[
\Ext^i_{\C}(X,Y) \cong \Ext^1_{\C}(K_i, Y) = 0.
\]
Hence, $\F^{\perp_1} = \F^\perp$ follows.

\item[$\ci$] We check conditions $\scpuno$, $\scptres$ and $\scpcuatro$ in \emph{Proposition~\ref{CP1}} for the pair $(\F,\G)$. 
	\begin{itemize}[itemsep=2pt,topsep=0pt]
	\item[$\scpuno$] The inclusions $\F, \G \subseteq \Sb$ are immediate. On the other hand, $\F$ is closed under direct summands in $\C$ since $\F$ is left thick.

	\item[$\scptres$] $\Ext^1_{\C}(\F,\G) = 0$ follows since $\G \subseteq \F^\perp$. \newpage

	\item[$\scpcuatro$] By $\ai$ and $\bi$, for every $S \in \Sb$ there exists a short exact sequence 
	\[
	0 \to \Ker(\varphi) \to F \xrightarrow{\varphi} S \to 0	
	\]
	where $\varphi \colon \F \to S$ is an $\F$-cover and $\Ker(\varphi) \in \F^\perp$. Moreover, $\Ker(\varphi) \in \Sb$ since $F, S \in \Sb$ and $\Sb$ is thick. Hence, $\Ker(\varphi) \in \F^\perp \cap \Sb = \G$. 
	\end{itemize}
So far we have proven that $(\F,\G)$ is a perfect left $\Sb$-cotorsion pair in $\C$, and the condition $\id_{\F}(\G) = 0$ follows by the definition of $\G$. 
\end{itemize}
\end{myproof}

\begin{proof}[\textbf{Proof of Theorem~\ref{theo:correspondence_with_G}}] \

\vspace{-0.1cm}
Consider the mapping
\begin{align*}
\Omega \colon \mathfrak{G} & \to \mathfrak{P}' \\
\F & \mapsto (\F, \G := \F^\perp \cap \Thick(\F)). 
\end{align*}
We show $\Omega$ is well defined and a left inverse of $\Theta$. Let $\F \in \mathfrak{G}$, that is, $\F$ is a left saturated class in $\C$ which is pre-covering in $\Thick(\F)$. Then $\F = \Thick^{-}(\F)$ and $\Proj(\C) \subseteq \F$. Since $\F$ is pre-covering in $\Thick(\F)$ with $\Proj(\C) \subseteq \F$ and $\C$ has enough projectives, then $\F$ is epic pre-covering. By \emph{Lemma~\ref{lem:correspondence_with_G}}, we have $(\F,\G) \in \mathfrak{P}'$. Finally, the equality $\Theta \circ \Omega = {\rm id}_{\mathfrak{G}}$ is easy to verify, and so $\Omega$ defines an injective embedding of $\mathfrak{G}$ into $\mathfrak{P}'$. 
\end{proof}

Is it possible to restrict $\Theta$ on a sub-class of $\mathfrak{P}'$, say $\tilde{\mathfrak{P}}$, in such a way that $\Theta|_{\tilde{\mathfrak{P}}}$ defines a one-to-one correspondence to a sub-class of $\mathfrak{G}$, say $\tilde{\mathfrak{G}}$? This question is settled in the following result.

\begin{theorem}\label{theo:correspondence_G_P_one_to_one}
Let $\C$ be a Krull-Schmidt abelian category with enough projectives and injectives. Then, there is a one-to-one correspondence between the following classes: 
\begin{align*}
\tilde{\mathfrak{P}} & := \left\{ (\F,\G) \subseteq \C^2 \mbox{ \emph{:} } \begin{array}{ll} \mbox{\emph{$(\F,\G)$ is a perfect $\Thick(\F)$-cotorsion pair in $\C$ with $\id_{\F}(\G) = 0$,}} \\ \mbox{\emph{$\Proj(\C) \subseteq \F$, and $\Inj(\C) \subseteq \Thick(\F)$}} \end{array} \right\}, \\
\tilde{\mathfrak{G}} & := \left\{ \F \subseteq \C \mbox{ {\rm :} } \begin{array}{ll} \mbox{\emph{$\F$ is a left saturated class in $\C$, and pre-covering in $\Thick(\F)$,}} \\ \mbox{\emph{such that $\Inj(\C) \subseteq \Thick(\F)$}} \end{array} \right\}.
\end{align*}
\end{theorem}

\begin{myproof}
First, we check that every pair $(\F,\G)$ in $\tilde{\mathfrak{P}}$ is mapped to $\tilde{\mathfrak{G}}$ via $\Theta$. By \emph{Proposition~\ref{prop:correspondence_with_G}}, we have that $\F$ is a left saturated class in $\C$, pre-covering in $\Thick(\F)$. The validity of the inclusion $\Inj(\C) \subseteq \Thick(\F)$ implies $\F \in \tilde{\mathfrak{G}}$. It follows that the restriction of $\Theta$ on $\tilde{\mathfrak{P}}$ defines a mapping $\tilde{\Theta} \colon \tilde{\mathfrak{P}} \to \tilde{\mathfrak{G}}$. 

Now we construct an inverse $\tilde{\Omega} \colon \tilde{\mathfrak{G}} \to \tilde{\mathfrak{P}}$. Set $\tilde{\Omega} := \Omega|_{\tilde{\mathfrak{G}}}$. We check $\tilde{\Omega}$ is well defined. Let $\F \in \tilde{\mathfrak{G}}$. By \emph{Lemma~\ref{lem:correspondence_with_G}}, we have that $(\F,\G := \F^\perp \cap \Thick(\F))$ is a perfect left $\Thick(\F)$-cotorsion pair in $\C$, such that $\id_{\F}(\G) = 0$. The inclusion $\Inj(\C) \subseteq \Thick(\F)$ is true by hypothesis, and $\Proj(\C) \subseteq \F$ holds since $\F$ is left saturated in $\C$. It is only left to show that $(\F,\G)$ is a right $\Thick(\F)$-cotorsion pair. We only check $\scpcinco$, since we already know $\scptres$, and $\scpdos$ is easy to verify. Let $S \in \Thick(\F)$. We show there exists a short exact sequence 
\[
0 \to S \to G \to F \to 0
\]
where $F \in \F$ and $G \in \G$. Since $\C$ has enough injectives, and $\Inj(\C) \subseteq \Thick(\F)$, there exists a short exact sequence
\[ 
0 \to S \to I \to C \to 0
\]
where $I \in \Inj(\C)$ and $C \in \Thick(\F)$. On the other hand, since $(\F,\G)$ is a left $\Thick(\F)$-cotorsion pair, there exists a short exact sequence
\[
0 \to G' \to F \to C \to 0
\]
where $F \in \F$ and $G' \in \G$. Taking the pullback of $I \to C$ and $F \to C$, we have the following commutative diagram with exact rows and columns, where the bottom right square is a pullback square: 
\begin{figure}[H]
\centering
\begin{tikzpicture}[description/.style={fill=white,inner sep=2pt}]
\matrix (m) [matrix of math nodes, row sep=2.5em, column sep=2.5em, text height=1.25ex, text depth=0.25ex]
{ & & 0 & 0 \\ 
& & G' & G' \\ 
0 & S & G & F & 0 \\ 
0 & S & I & C & 0 \\ 
& & 0 & 0 \\ };
\path[->]
(m-3-3)-- node[pos=0.5] {\footnotesize$\mbox{\bf pb}$} (m-4-4)
(m-1-3) edge (m-2-3) (m-1-4) edge (m-2-4)
(m-2-3) edge (m-3-3) (m-2-4) edge (m-3-4)
(m-3-1) edge (m-3-2) (m-3-2) edge (m-3-3) (m-3-3) edge (m-3-4) edge (m-4-3) (m-3-4) edge (m-3-5) edge (m-4-4)
(m-4-1) edge (m-4-2) (m-4-2) edge (m-4-3) (m-4-3) edge (m-4-4) edge (m-5-3) (m-4-4) edge (m-4-5) edge (m-5-4);
\path[-,font=\scriptsize]
(m-2-3) edge [double, thick, double distance=2pt] (m-2-4)
(m-3-2) edge [double, thick, double distance=2pt] (m-4-2);
\end{tikzpicture}
\caption[An application of pullbacks]{Pullbacks along epimorphisms preserves short exact sequences.}
\label{fig:pullback}
\end{figure}
The central row satisfies $\scpcinco$, that is, $\G$ is a special pre-enveloping in $\Thick(\F)$. Moreover, $\G$ is right saturated in $\C$. To show this, note that it is clear that $\G$ is closed under extensions and direct summands in $\C$. On the other hand $\G$ is closed under cokernels of monomorphisms in $\G$, since $\G := \F^\perp \cap \Thick(\F)$, where $\F^\perp$ and $\Thick(\F)$ are both closed under taking cokernels of monomorphisms. Finally, the inclusion $\Inj(\C) \subseteq \F^\perp$ is clear, and $\Inj(\C) \subseteq \Thick(\F)$ by hypothesis, and hence $\Inj(\C) \subseteq \G$. We have that $\G$ is a right saturated class in $\C$ which is pre-enveloping in $\Thick(\F)$. By \emph{Lemma~\ref{lem:correspondence_with_G}}, we have that $\G$ is the right half of a perfect right $\Thick(\F)$-cotorsion pair in $\C$, and so $\G$ is enveloping in $\Thick(\F)$. Therefore, $(\F,\G)$ is a perfect $\Thick(\F)$-cotorsion pair in $\C$, and thus we can define a mapping $\tilde{\Omega} \colon \tilde{\mathfrak{G}} \to \tilde{\mathfrak{P}}$ as $\tilde{\Omega} := \Omega|_{\tilde{\mathfrak{G}}}$. 

Now, we show that $\tilde{\Theta}$ and $\tilde{\Omega}$ are inverse of each other. The equality $\tilde{\Theta} \circ \tilde{\Omega} = {\rm id}_{\tilde{\mathfrak{G}}}$ is easy to verify. Now let $(\F,\G) \in \tilde{\mathfrak{P}}$. We have:
\begin{align*}
\tilde{\Omega} \circ \tilde{\Theta}(\F,\G) & = \tilde{\Omega}(\F) = (\F, \F^\perp \cap \Thick(\F)).
\end{align*}
Since $(\F,\G)$ is a right $\Thick(\F)$-cotorsion pair and $\F^\perp = \F^{\perp_1}$, the equality $\G = \F^{\perp_1} \cap \Thick(\F) = \F^\perp \cap \Thick(\F)$ follows, and so $\tilde{\Omega} \circ \tilde{\Theta}(\F,\G) = (\F,\G)$. Therefore, $\tilde{\Theta}$ defines a bijection between $\tilde{\mathfrak{P}}$ and $\tilde{\mathfrak{G}}$. 
\end{myproof}

Using \emph{Theorem~\ref{theo:correspondence_G_P_one_to_one}}, the correspondence from \emph{Theorem~\ref{theo:correspondence_with_AB-context}} can be extended in the case where $\C$ is a Krull-Schmidt abelian category with enough projectives and injectives. Specifically, we have the following result.

\begin{corollary}\label{coro:correspondence_with_AB-context}
Let $\C$ be a Krull-Schmidt abelian category with enough projectives and injectives. Then there exists a one-to-one correspondence between the following classes:
\begin{align*}
\tilde{\mathfrak{F}} & := \left\{ (\X, \omega) \subseteq \C^2 \mbox{ \emph{: $(\X,\omega)$ is a left Frobenius pair in $\C$, with $\Proj(\C) \subseteq \X$, and $\Inj(\C) \subseteq \X^\wedge$}} \right\}, \\
\tilde{\mathfrak{P}} & := \left\{ (\F,\G) \subseteq \C^2 \mbox{ \emph{:} } \begin{array}{ll} \mbox{\emph{$(\F,\G)$ is a perfect $\Thick(\F)$-cotorsion pair in $\C$ with $\id_{\F}(\G) = 0$,}} \\ \mbox{\emph{$\Proj(\C) \subseteq \F$, and $\Inj(\C) \subseteq \Thick(\F)$}} \end{array} \right\}, \\
\tilde{\mathfrak{C}} & := \left\{ (\A,\B) \subseteq \C^2 \mbox{\emph{ : }} \begin{array}{ll} \mbox{\emph{$(\A,\B)$ is a left weak AB context in $\C$, with $\Proj(\C) \subseteq \A$, and}} \\ \mbox{\emph{$\Inj(\C) \subseteq \Thick(\A)$}} \end{array} \right\}, \\
\tilde{\mathfrak{G}} & := \left\{ \F \subseteq \C \mbox{ \emph{:}} \begin{array}{ll} \mbox{\emph{$\F$ is a left saturated class in $\C$, and pre-covering in $\Thick(\F)$,}} \\ \mbox{\emph{such that $\Inj(\C) \subseteq \Thick(\F)$}} \end{array} \right\},
\end{align*}
where $\tilde{\mathfrak{P}} = \tilde{\mathfrak{C}}$, and the corresponding mappings are denoted by $\tilde{\Theta}$, $\tilde{\Omega}$, $\tilde{\Psi}$, and $\tilde{\Phi}$.

Dually, there is a one-to-one correspondence between the classes:
\emph{
\begin{align*}
\tilde{\mathfrak{F}}\op & := \left\{ (\nu, \Y) \subseteq \C^2 \mbox{ : $(\nu,\Y)$ is a right Frobenius pair in $\C$, with $\Inj(\C) \subseteq \Y$, and $\Proj(\C) \subseteq \Y^\vee$} \right\}, \\
\tilde{\mathfrak{P}}\op & := \left\{ (\F,\G) \subseteq \C^2 \mbox{ \emph{:} } \begin{array}{ll} \mbox{$(\F,\G)$ is a perfect $\Thick(\G)$-cotorsion pair in $\C$ with $\pd_{\G}(\F) = 0$,} \\ \mbox{$\Inj(\C) \subseteq \G$, and $\Proj(\C) \subseteq \Thick(\G)$} \end{array} \right\}, \\
\tilde{\mathfrak{C}}\op & := \left\{ (\A,\B) \subseteq \C^2 \mbox{ \emph{:} } \begin{array}{ll} \mbox{$(\A,\B)$ is a right weak AB context in $\C$, with $\Inj(\C) \subseteq \B$, and} \\ \mbox{$\Proj(\C) \subseteq \Thick(\B)$} \end{array} \right\}, \\
\tilde{\mathfrak{G}}\op & := \left\{ \G \subseteq \C \mbox{ \emph{:} } \begin{array}{ll} \mbox{$\G$ is a right saturated class in $\C$, and pre-enveloping in $\Thick(\G)$,} \\ \mbox{such that $\Proj(\C) \subseteq \Thick(\G)$} \end{array} \right\},
\end{align*}
}
where $\tilde{\mathfrak{P}}\op = \tilde{\mathfrak{C}}\op$, and the corresponding mappings are denoted by $\tilde{\Theta}\op$, $\tilde{\Omega}\op$, $\tilde{\Psi}\op$, and $\tilde{\Phi}\op$.
\end{corollary}

\begin{figure}[H]
\centering 
\scriptsize
\begin{tikzpicture}[description/.style={fill=white,inner sep=2pt}]
\matrix (m) [matrix of math nodes, row sep=0.5em, column sep=0.5em, text height=3.5ex, text depth=1.5ex]
{ {\color{red}{(\A,\A \cap \B)}} & {\Large\bm{\tilde{\mathfrak{F}}}} & {} & {} & {} & {\Large\bm{\tilde{\mathfrak{G}}}} & {\color{red}{\F}} \\
  {} & {} & {\color{red}{(\X,\omega)}} & {} & {\color{red}{\F}} & {} & {} \\
  {} & {} & {} & {} & {} & {} \\
  {} & {} & {\color{red}{(\X,\omega^\wedge)}} & {} & {\color{red}{(\F, \F^\perp \cap \Thick(\F))}} & {} & {} \\
  {\color{red}{(\A,\B)}} & {\Large\bm{\tilde{\mathfrak{C}}}} & {} & {} & {} & {\Large\bm{\tilde{\mathfrak{P}}}} & {\color{red}{(\F,\G)}} \\ };
\path[red,|->]
(m-5-1) edge node[above,sloped] {$\tilde{\Psi}$} (m-1-1)
(m-2-3) edge node[above,sloped] {$\tilde{\Phi}$} (m-4-3)
(m-2-5) edge node[above,sloped] {$\tilde{\Omega}$} (m-4-5)
(m-5-7) edge node[below,sloped] {$\tilde{\Theta}$} (m-1-7)
;
\path[->] 
($(m-1-2.east)+(0,-0.2)$) edge [thick] ($(m-1-6.west)-(0,0.2)$)
($(m-1-2.east)+(-0.2,-0.4)$) edge [thick] node[right] {$\tilde{\Phi}$} ($(m-5-2.west)-(-0.4,-0.4)$)
($(m-1-6.east)+(-0.45,-0.4)$) edge [thick] node[left] {$\tilde{\Omega}$} ($(m-5-6.west)-(-0.25,-0.4)$)
;
\path[<-] 
($(m-1-2.east)+(0,0)$) edge [thick] ($(m-1-6.west)-(0,-0)$)
($(m-1-2.east)+(-0.35,-0.4)$) edge [thick] node[left] {$\tilde{\Psi}$} ($(m-5-2.west)-(-0.25,-0.4)$)
($(m-1-6.east)+(-0.3,-0.4)$) edge [thick] node[right] {$\tilde{\Theta}$} ($(m-5-6.west)-(-0.4,-0.4)$)
;
\path[-,font=\scriptsize]
(m-5-2) edge [double, thick, double distance=2pt] (m-5-6)
;
\end{tikzpicture}
\caption[Correspondences between left Frobenius pairs, left weak AB contexts, and relative cotorsion pairs in abelian categories, and pre-covering classes in Krull-Schmidt abelian categories]{Correspondences between left Frobenius pairs, left weak AB contexts, and relative cotorsion pairs, and pre-covering classes in Krull-Schmidt abelian categories.}
\label{fig:correspondence_KS}
\end{figure}

\begin{figure}[H]
\centering 
\scriptsize
\begin{tikzpicture}[description/.style={fill=white,inner sep=2pt}]
\matrix (m) [matrix of math nodes, row sep=0.5em, column sep=0.5em, text height=3.5ex, text depth=1.5ex]
{ {\color{blue}{(\A \cap \B,\B)}} & {\Large\bm{\tilde{\mathfrak{F}}\op}} & {} & {} & {} & {\Large\bm{\tilde{\mathfrak{G}}\op}} & {\color{blue}{\G}} \\
  {} & {} & {\color{blue}{(\nu,\Y)}} & {} & {\color{blue}{\G}} & {} & {} \\
  {} & {} & {} & {} & {} & {} \\
  {} & {} & {\color{blue}{(\nu^\vee,\Y)}} & {} & {\color{blue}{({}^\perp\G \cap \Thick(\G),\G)}} & {} & {} \\
  {\color{blue}{(\A,\B)}} & {\Large\bm{\tilde{\mathfrak{C}}\op}} & {} & {} & {} & {\Large\bm{\tilde{\mathfrak{P}}\op}} & {\color{blue}{(\F,\G)}} \\ };
\path[blue,|->]
(m-5-1) edge node[above,sloped] {$\tilde{\Psi}\op$} (m-1-1)
(m-2-3) edge node[above,sloped] {$\tilde{\Phi}\op$} (m-4-3)
(m-2-5) edge node[above,sloped] {$\tilde{\Omega}\op$} (m-4-5)
(m-5-7) edge node[below,sloped] {$\tilde{\Theta}\op$} (m-1-7)
;
\path[->] 
($(m-1-2.east)+(0,-0.2)$) edge [thick] ($(m-1-6.west)-(0,0.2)$)
($(m-1-2.east)+(-0.6,-0.4)$) edge [thick] node[right] {$\tilde{\Phi}\op$} ($(m-5-2.west)-(-0.4,-0.4)$)
($(m-1-6.east)+(-0.9,-0.4)$) edge [thick] node[left] {$\tilde{\Omega}\op$} ($(m-5-6.west)-(-0.25,-0.4)$)
;
\path[<-] 
($(m-1-2.east)+(0,0)$) edge [thick] ($(m-1-6.west)-(0,-0)$)
($(m-1-2.east)+(-0.75,-0.4)$) edge [thick] node[left] {$\tilde{\Psi}\op$} ($(m-5-2.west)-(-0.25,-0.4)$)
($(m-1-6.east)+(-0.75,-0.4)$) edge [thick] node[right] {$\tilde{\Theta}\op$} ($(m-5-6.west)-(-0.4,-0.4)$)
;
\path[-,font=\scriptsize]
(m-5-2) edge [double, thick, double distance=2pt] (m-5-6)
;
\end{tikzpicture}
\caption[Correspondences between right Frobenius pairs, right weak AB contexts, and relative cotorsion pairs, and pre-enveloping classes in Krull-Schmidt abelian categories]{Correspondences between right Frobenius pairs, right weak AB contexts, and relative cotorsion pairs, and pre-enveloping classes in Krull-Schmidt abelian categories.}
\label{fig:correspondence_KS_dual}
\end{figure}
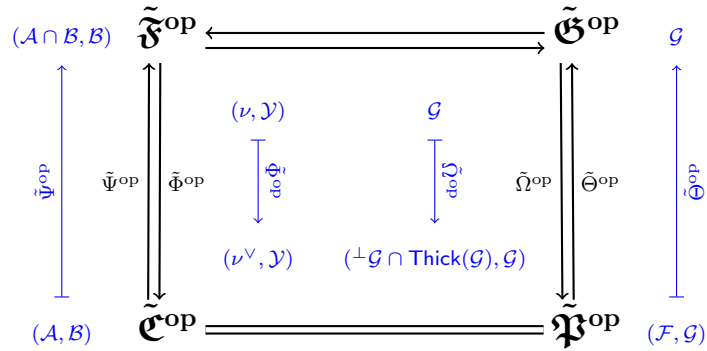

%%%%%%%%%%%%%%%%%%%%%%%%%%%%%%%%%%
%%%%%%%%%%%%%%%%%%%%%%%%%%%%%%%%%%

\subsection{Some remarks on perfect cotorsion pairs}

The problem of obtaining covers by classes of modules has had an important interest recently in Homological Algebra and Representation Theory of Algebras. This in part has been motivated by the Flat Cover Conjecture\footnote{Every left $R$-module has a flat cover.}, settled by L. Bican, R. El Bashir, and E. E. Enochs in \emph{\cite{FlatCover}}. Several authors have studied conditions under which it is possible to obtain covers. For example, it is known that every left $R$-module over a perfect ring has a projective cover. This result is also valid in the category $\mathsf{mod}(\Lambda)$ of finitely generated modules over an Artin $R$-algebra $\Lambda$, where $R$ is a commutative Artinian ring with identity. In a more general context, H. Holm and P. J\o rgensen have established in \emph{\cite[Theorem 3.4]{HolmJorgensen}} certain conditions under which a class $\F$ of left $R$-modules is covering. Namely, if $\F$ contains the ground ring $R$ and is closed under extensions, direct sums, pure sub-modules, and pure quotient of modules, then $(\F,\F^{\perp_1})$ is a perfect cotorsion pair, and hence $\F$ is covering. In the following result, we provide other conditions under which a class of object in an abelian category is covering.

\begin{corollary}\label{coro:Holm_Jorgensen}
Let $\C$ be a Krull-Schmidt abelian category with enough projectives and injectives. Then, there is a one-to-one correspondence between the following classes: 
\emph{
\begin{align*}
\tilde{\mathfrak{P}}(\C) & := \left\{ (\F,\G) \subseteq \C^2 \mbox{ : $(\F,\G)$ is a perfect hereditary cotorsion pair in $\C$ with $\Thick(\F) = \C$} \right\}, \\
\tilde{\mathfrak{G}}(\C) & := \left\{ \F \subseteq \C \mbox{ : $\F$ is a left saturated and pre-covering class in $\C$ such that $\Thick(\F) = \C$} \right\}.
\end{align*}
}
\end{corollary}

\begin{myproof}
\sloppypar{By \emph{Lemma~\ref{lem:hereditary_exact}} and \emph{Theorem~\ref{theo:correspondence_G_P_one_to_one}}, the bijection is given by $(\F,\G) \mapsto \F$ with inverse \mbox{$\F \mapsto (\F,\F^\perp)$.}} 
\end{myproof}

\begin{remark}
For the case $\C := \mathsf{mod}(\Lambda)$, where $\Lambda$ is an Artin algebra, M. Auslander and I. Reiten proved in \emph{\cite[Proposition 3.3]{AR}} that $(\F,\F^\perp)$ is a perfect and hereditary cotorsion pair in $\mathsf{mod}(\Lambda)$ with $\F^\wedge = \mathsf{mod}(\Lambda)$, whenever $\F$ is a left saturated and pre-covering class in $\mathsf{mod}(\Lambda)$ such that $\F^\wedge = \mathsf{mod}(\Lambda)$. The previous corollary states that these two assertions are in fact equivalent. 
\end{remark}

%%%%%%%%%%%%%%%%%%%%%%%%%%%%%%%%%%%
%%%%%%%%%%%%%%%%%%%%%%%%%%%%%%%%%%%

\subsection{Relationship with some Auslander-Reiten correspondence theorems}

In this section, we give some remarks on the relation between perfect cotorsion pairs, pre-enveloping classes, and cotilting modules, within the framework of the correspondences we have studied so far. In their seminal paper \textit{Applications of Contravariantly Finite Subcategories} \emph{\cite{AR}}, M. Auslander and I. Reiten proved the following correspondence theorem:

\begin{theorem}\label{theo:Auslander_Reiten} 
Let $\Lambda$ be an Artin $R$-algebra.  
\begin{itemize}[itemsep=2pt,topsep=0pt]
\item[$\ai$] There exists a one-to-one correspondence between the following classes:
	\begin{itemize}[itemsep=2pt,topsep=0pt] 
	\item[$\ai\text{-}\iroman$] The class of isomorphisms classes of basic cotilting modules.
	
	\item[$\ai\text{-}\iiroman$] The class of contravariantly finite resolving sub-categories $\F \subseteq \mathsf{mod}(\Lambda)$ satisfying $\F^\wedge = \mathsf{mod}(\Lambda)$. 
	
	\item[$\ai\text{-}\iiiroman$] The class of complete cotorsion pairs $(\F,\G)$ with $\F$ resolving and $\F^\wedge = \mathsf{mod}(\Lambda)$. 
	\end{itemize}

\item[$\bi$] There exists a one-to-one correspondence between the classes:
	\begin{itemize}[itemsep=2pt,topsep=0pt]
	\item[$\bi\text{-}\iroman$] The class of isomorphism classes of basic cotilting modules.
 	
	\item[$\bi\text{-}\iiroman$] The class of covariantly finite coresolving sub-categories of $\mathsf{inj}^{< \infty}(\Lambda)$. 
	\end{itemize}
\end{itemize}
\end{theorem}

The correspondence $\ai\text{-}\iroman \leftrightarrow \ai\text{-}\iiroman$ is given by $[C] \mapsto {}^\perp C$, where $[C]$ denotes the isomorphism class of $C$, that is, the class of finitely generated left $\Lambda$-modules isomorphic to $C$. The inverse of the previous correspondence is given by $\F \mapsto [C_\F]$, where $C_\F$ is the basic cotilting $\Lambda$-module defined as the direct sum of indecomposable $\X$-injective modules. The proof of this can be found in \emph{\cite[Chapter 8 by I. Reiten, Theorem 2.2 (c)]{handbook}} or in \emph{\cite[Theorem 5.5 (a)]{AR}}. 

On the other hand, the correspondence $\ai\text{-}\iroman \leftrightarrow \ai\text{-}\iiiroman$ is given by $[C] \mapsto ({}^\perp C, \add(\C)^\wedge)$, whose inverse is defined as $(\F,\G) \mapsto [C_{\F\cap\G}]$, where $C_{\F\cap\G}$ is the direct sum of pairwise non-isomorphic indecomposable finitely generated left $\Lambda$-modules in $\F \cap \G$. This fact is proven in \emph{\cite[Chapter 8 by I. Reiten, Corollary 2.3 (b)]{handbook}}. 

Finally, in \emph{\cite[Theorem 5.5 (b)]{AR}} or \emph{\cite[Chapter 8 by I. Reiten, Theorem 2.2 (d)]{handbook}}, one can check that the correspondence $\bi\text{-}\iroman \leftrightarrow \bi\text{-}\iiroman$ is given by $[C] \mapsto \add(C)^\wedge$, whose inverse is given by $\G \mapsto [C_\G]$, where $C_{\G}$ is the direct sum of pairwise non-isomorphic indecomposable $\G$-projective finitely generated left $\Lambda$-modules. 

The previous theorem is one of the motivations of \emph{Corollary~\ref{coro:Holm_Jorgensen}}. This can be more appreciated if we explain how to connect \emph{Corollary~\ref{coro:Holm_Jorgensen}} and \emph{Theorem~\ref{theo:Auslander_Reiten}}. Before that, and for a better understanding, we recall some definitions. 

For the rest of this section, we explain how the statements $\ai$ and $\bi$ are actually equivalent. Before doing that, we recall the terminology used in the previous theorem: 
\begin{itemize}[itemsep=2pt,topsep=0pt]
\item[$\bullet$] \emph{\cite[Definition 8.1.1]{GT}}: A left $R$-module $C$ is \textbf{cotilting} provided that:
\begin{itemize}[itemsep=2pt,topsep=0pt]
\item[$\ai$] $\id(C) < \infty$.

\item[$\bi$] $\Ext^i_R(C^I, C) = 0$ for every $i > 0$ and every index set $I$, where $C^I$ denotes the direct product of copies of $C$ indexed by $I$. 

\item[$\ci$] $\resdim_{\mathsf{Prod}(C)}(Q) < \infty$, where $\mathsf{Prod}(C)$ denotes the class of direct summands of arbitrary direct products of copies of $C$ and $Q$ is an injective cogenerator.
\end{itemize}

\item[$\bullet$] \emph{\cite{AR}}: A left $R$-module $M$ is \textbf{basic} if in a direct sum decomposition no indecomposable module appears more than once.

\item[$\bullet$] A left $R$-module $M$ is \textbf{self-orthogonal} if $\Ext^i_R(M,M) = 0$ for every $i > 0$. 

\item[$\bullet$] In \emph{\cite{AR}}, pre-covering classes are called \textbf{contravariantly finite}, and pre-enveloping classes are called \textbf{covariantly finite}. 

\item[$\bullet$] In \emph{\cite{AR}}, sub-categories of $\mathsf{mod}(\Lambda)$ are considered to be closed under isomorphisms and direct summands, so it follows that resolving and coresolving classes in \emph{\cite{AR}} are left saturated and right saturated in $\mathsf{mod}(\Lambda)$, respectively, according to our terminology. 
\end{itemize}

For the rest of this section, our goal will be to present the previous theorem in a more general context. For this, we will define a categorical setting called \emph{Auslander-Reiten category} (see \emph{Definition~\ref{def:AR_category}}).

\begin{definition}\label{def:AR_category}
Let $\C$ be an abelian category. We say that a full sub-category $\X$ of $\C$ satisfies the \textbf{Auslander-Reiten condition} (\textbf{AR condition} for short) if the following two conditions are equivalent:
\begin{itemize}[itemsep=2pt,topsep=0pt]
\item[$\mathsf{(ar1)}$] $\X \subseteq \mathsf{Inj}^{< \infty}(\C)$.

\item[$\mathsf{(ar2)}$] $({}^\perp \X)^\wedge = \C$. 
\end{itemize}

If in addition $\C$ has enough projectives and injectives, we say that $\C$ is an \textbf{Auslander-Reiten category} if every right saturated and special pre-enveloping sub-category of $\C$ satisfies the AR condition.
\end{definition}

\begin{example}\label{ex:AR_category} 
Consider the category $\mathsf{mod}(\Lambda)$ of finitely generated left modules over an Artin algebra $\Lambda$. In \cite[Propositions 5.3 and 5.5]{AR}, Auslander and Reiten proved that if $\G$ is a right saturated pre-enveloping sub-category of $\mathsf{mod}(\Lambda)$, and ${}^\perp\G$ is the associated left saturated pre-covering sub-category of $\mathsf{mod}(\Lambda)$, then $({}^\perp\G)^\wedge = \mathsf{mod}(\Lambda)$ if, and only if, $\G \subseteq \mathsf{inj}^{< \infty}(\Lambda)$, where $\mathsf{inj}^{< \infty}(\Lambda)$ denotes the class of finitely generated left $\Lambda$-modules with finite injective dimension. \emph{Definition~\ref{def:AR_category}} is motivated on this result. 
\end{example}

\begin{proposition}\label{prop:AR_condition}
Let $\C$ be an abelian category with enough projectives and injectives. Then there exists a one-to-one correspondence between the following classes:
\begin{align*}
\tilde{\mathfrak{G}}\op_{\mathsf{AR}}(\C) & := \left\{ \G \subseteq \C \mbox{ \emph{:} } \begin{array}{ll} \mbox{\emph{$\G$ is a right saturated and special pre-enveloping sub-category of $\C$,}} \\ \mbox{\emph{such that $\G$ satisfies the AR condition}} \end{array} \right\}, \\
\tilde{\mathfrak{P}}_{\mathsf{AR}}(\C) & := \left\{ (\F,\G) \subseteq \C^2 \mbox{ \emph{:} } \begin{array}{ll} \mbox{\emph{$(\F,\G)$ is a complete and hereditary cotorsion pair in $\C$, such that}} \\ \mbox{\emph{$\G$ satisfies the AR condition}} \end{array} \right\}.
\end{align*}
Consider the classes:
\emph{
\begin{align*}
\tilde{\mathfrak{C}}_{\mathsf{AR}}(\mathcal{C}) & := \{ (\A,\B) \subseteq \C^2 \mbox{ : $(\A,\B)$ is a left AB context in $\C$} \}, \\
\tilde{\mathfrak{F}}_{\mathsf{AR}}(\mathcal{C}) & := \{ (\X,\omega) \subseteq \C^2 \mbox{ : $(\X,\omega)$ is a left Frobenius pair with $\X^\wedge = \C$} \}.
\end{align*}
}
If in addition, $\C$ is an Auslander-Reiten category, then the following equalities hold:
\begin{itemize}[itemsep=2pt,topsep=0pt]
\item[$\bullet$] $\tilde{\mathfrak{G}}\op_{\mathsf{AR}}(\C) = \left\{ \G \subseteq \C \mbox{ \emph{:} } \begin{array}{ll} \mbox{\emph{$\G$ is a right saturated and special pre-enveloping}} \\ \mbox{\emph{sub-category of $\Inj^{< \infty}(\C)$}} \end{array} \right\}$.

\item[$\bullet$] $\tilde{\mathfrak{P}}_{\mathsf{AR}}(\C) = \left\{ (\F,\G) \subseteq \C^2 \mbox{ \emph{:} } \begin{array}{ll} \mbox{\emph{$(\F,\G)$ is a complete hereditary cotorsion pair in $\C$,}} \\ \mbox{\emph{with $\Thick(\F) = \C$}} \end{array} \right\}$. 

\item[$\bullet$] $\tilde{\mathfrak{P}}_{\mathsf{AR}}(\C) = \tilde{\mathfrak{C}}_{\mathsf{AR}}(\C)$.
\end{itemize}
Moreover, there exists a one-to-one correspondence between the classes $\tilde{\mathfrak{C}}_{\mathsf{AR}}(\mathcal{C})$ and $\tilde{\mathfrak{F}}_{\mathsf{AR}}(\C)$.
\end{proposition}

\begin{myproof}
Consider the mapping $\tilde{\Omega}_{\mathsf{AR}} \colon \tilde{\mathfrak{G}}\op_{\mathsf{AR}}(\C) \to \tilde{\mathfrak{P}}_{\mathsf{AR}}(\C)$ given by $\G \mapsto ({}^\perp\G,\G)$ for every $\G \in \tilde{\mathfrak{G}}\op_{\mathsf{AR}}(\C)$. Let us check that it is well defined. It suffices to show that $({}^\perp\G,\G)$ is indeed a cotorsion pair. The fact that $({}^\perp\G,\G)$ is complete and hereditary will follow from the hypothesis that $\C$ has enough projectives and injectives, and the properties of $\G$. First note that class ${}^\perp\G$ is left thick, and $\G$ is right thick. Consider $\omega := {}^\perp\G \cap \G$. Since $\G$ is special pre-enveloping and right saturated in $\C$, for every $X \in {}^\perp\G$ there exists a short exact sequence 
\[
0 \to X \to G \to K \to 0
\]
where $G \in \G$ and $K \in {}^{\perp_1}\G = {}^\perp\G$. Now since ${}^\perp\G$ is closed under extensions, we have $G \in {}^\perp\G \cap \G = \omega$. On the other hand, it is clear that $\id_{{}^\perp\G}(\omega) = 0$. Then, $\omega$ is a ${}^\perp\G$-injective relative cogenerator in ${}^\perp\G$. Finally, note that $\G \subseteq ({}^\perp\G)^\wedge$ since $\G$ satisfies the Auslander-Reiten condition. Thus, by \emph{Theorem~\ref{theo:ponja}}, we have $\G = ({}^\perp\G)^\wedge \cap ({}^\perp\G)^{\perp_1} = ({}^\perp\G)^{\perp_1}$. Hence, the mapping $\tilde{\Omega}_{\mathsf{AR}}$ is well defined, and its inverse is given by the projection $\tilde{\Omega}_{\mathsf{AR}} \colon (\F,\G) \mapsto \G$. 

Now suppose that $\C$ is an Auslander-Reiten category. The equalities 
\begin{align*}
\tilde{\mathfrak{G}}\op_{\mathsf{AR}}(\C) & = \{ \G \subseteq \C \mbox{ : $\G$ is a right saturated and special pre-enveloping sub-category of $\Inj^{< \infty}(\C)$} \}, \\
\tilde{\mathfrak{P}}_{\mathsf{AR}}(\C) & = \{ (\F,\G) \subseteq \C^2 \mbox{ : $(\F,\G)$ is a complete hereditary cotorsion pair in $\C$, with $\Thick(\F) = \C$} \},
\end{align*}
are clear. The inclusion $\tilde{\mathfrak{P}}_{\mathsf{AR}}(\C) \supseteq \tilde{\mathfrak{C}}_{\mathsf{AR}}(\C)$ follows by \emph{Theorem~\ref{theo:from_AB-context_to_hereditary_pairs}}, and $\tilde{\mathfrak{P}}_{\mathsf{AR}}(\C) \subseteq \tilde{\mathfrak{C}}_{\mathsf{AR}}(\C)$ by \emph{Theorem~\ref{theo:ThickF-cotorsion_pair}}.

Finally, it is not hard to see that the mappings
\[
\begin{tikzpicture}[description/.style={fill=white,inner sep=2pt}]
\matrix (m) [matrix of math nodes, row sep=1em, column sep=2em, text height=1.25ex, text depth=0.25ex]
{ (\X,\omega^\wedge) & (\X,\omega) \\
  \tilde{\mathfrak{C}}_{\mathsf{AR}}(\C) & \tilde{\mathfrak{F}}_{\mathsf{AR}}(\C) \\
  (\A,\B) & (\A,\A \cap \B) \\ };
\path[->]
($(m-2-1.east)+(0,-0.1)$) edge node[below] {\footnotesize$\tilde{\Psi}_{\mathsf{AR}}$} ($(m-2-2.west)+(0,-0.1)$)
;
\path[<-]
($(m-2-1.east)+(0,+0.1)$) edge node[above] {\footnotesize$\tilde{\Phi}_{}\mathsf{AR}$} ($(m-2-2.west)+(0,+0.1)$)
;
\path[|->]
(m-1-2) edge (m-1-1)
(m-3-1) edge (m-3-2)
;
\end{tikzpicture}
\]
are well defined and inverse to each other. 
\end{myproof}

\begin{figure}[H]
\centering 
\scriptsize
\begin{tikzpicture}[description/.style={fill=white,inner sep=2pt}]
\matrix (m) [matrix of math nodes, row sep=0.5em, column sep=3em, text height=3.5ex, text depth=1.5ex]
{ {\color{red}{(\F,\F\cap\G)}} & {\color{red}{(\F,\G)}} & {} & {\color{red}{(\F,\G)}} & {\color{blue}{\G}} \\
  {\Large\bm{\tilde{\mathfrak{F}}_{\mathsf{AR}}(\C)}} & {\Large\bm{\tilde{\mathfrak{C}}_{\mathsf{AR}}(\C)}} & {} & {\Large\bm{\tilde{\mathfrak{P}}_{\mathsf{AR}}(\C)}} & {\Large\bm{\tilde{\mathfrak{G}}\op_{\mathsf{AR}(\C)}}} \\
  {\color{red}{(\X,\omega)}} & {\color{red}{(\X,\omega^\wedge)}} & {} & {\color{red}{({}^\perp\G,\G)}} & {\color{blue}{\G}} \\ };
\path[red,|->]
(m-1-2) edge node[above] {$\tilde{\Psi}_{\mathsf{AR}}$} (m-1-1)
(m-3-1) edge node[below] {$\tilde{\Phi}_{\mathsf{AR}}$} (m-3-2)
(m-3-5) edge node[below] {$\tilde{\Omega}_{\mathsf{AR}}$} (m-3-4)
;
\path[blue,|->]
(m-1-4) edge node[above] {$\tilde{\Theta}_\mathsf{AR}$} (m-1-5)
;
\path[->] 
($(m-2-1.east)+(0,-0.1)$) edge [thick] node[below] {$\tilde{\Phi}_{\mathsf{AR}}$} ($(m-2-2.west)-(0,0.1)$)
($(m-2-4.east)+(0,0.1)$) edge [thick] node[above] {$\tilde{\Theta}_\mathsf{AR}$} ($(m-2-5.west)-(0,-0.1)$)
;
\path[<-] 
($(m-2-1.east)+(0,0.1)$) edge [thick] node[above] {$\tilde{\Psi}_{\mathsf{AR}}$} ($(m-2-2.west)-(0,-0.1)$)
($(m-2-4.east)+(0,-0.1)$) edge [thick] node[below] {$\tilde{\Omega}_{\mathsf{AR}}$} ($(m-2-5.west)-(0,0.1)$)
;
\path[-,font=\scriptsize]
(m-2-2) edge [double, thick, double distance=2pt] (m-2-4)
;
\end{tikzpicture}
\caption[Correspondences between left Frobenius pairs, left AB contexts, complete and hereditary cotorsion pairs, and special pre-enveloping sub-categories in Auslander-Reiten categories]{Correspondences between left Frobenius pairs, left AB contexts, complete and hereditary cotorsion pairs, and special pre-enveloping sub-categories in Auslander-Reiten categories.}
\label{fig:correspondence_KS}
\end{figure}

The following result is a categorical version of \emph{Theorem~\ref{theo:Auslander_Reiten}}. It is a consequence of \emph{Proposition~\ref{prop:AR_condition}} and \emph{Corollary~\ref{coro:Holm_Jorgensen}}.

\begin{corollary}\label{coro:ARKS}
Let $\C$ be an Auslander-Reiten and Krull-Schmidt category. Then $\tilde{\mathfrak{P}}_{\mathsf{AR}}(\C) = \tilde{\mathfrak{P}}(\C)$, and there exists a one-to-one correspondence between the classes $\tilde{\mathfrak{G}}\op_{\mathsf{AR}}(\C)$ and $\tilde{\mathfrak{G}}(\C)$. 
\end{corollary}

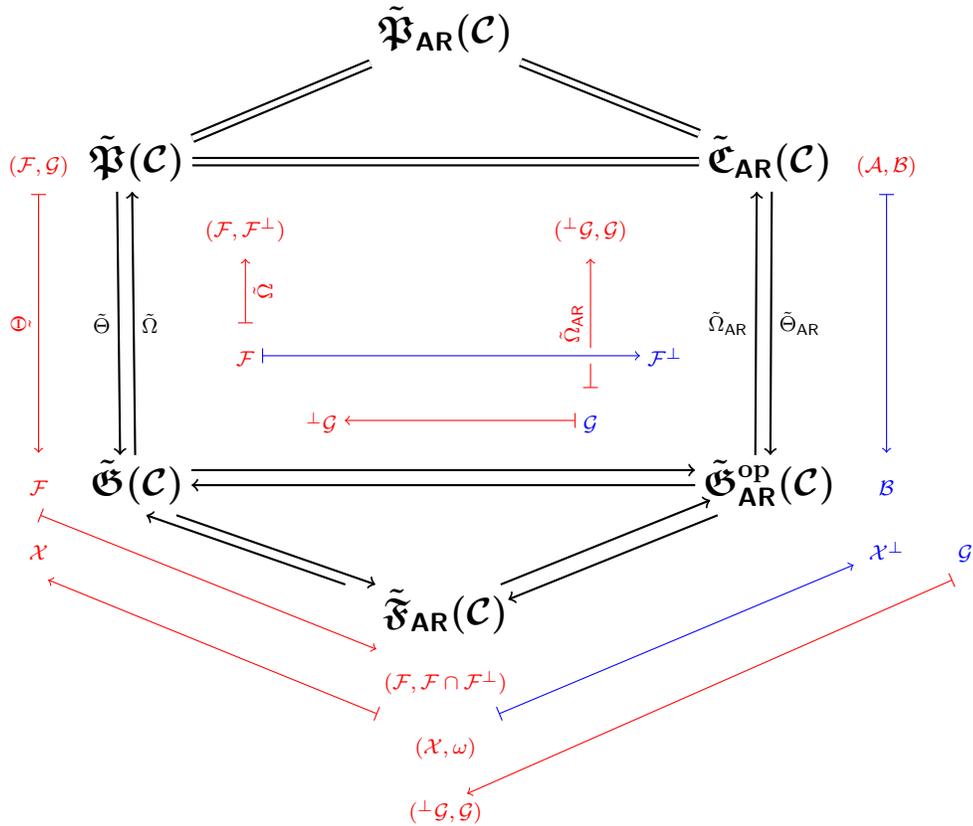
\begin{figure}[H]
\centering 
\scriptsize
\begin{tikzpicture}[description/.style={fill=white,inner sep=2pt}]
\matrix (m) [matrix of math nodes, row sep=0.05em, column sep=0.3em, text height=3.5ex, text depth=1.5ex]
{ {} & {} & {} & {} & {} & {} & {} & {\Large\bm{\tilde{\mathfrak{P}}_{\mathsf{AR}}(\C)}} & {} & {} & {} & {} & {} & {} & {} \\
  {} & {} & {} & {} & {} & {} & {} & {} & {} & {} & {} & {} & {} & {} & {} \\
  {} & {} & {\color{red}{(\F,\G)}} & {\Large\bm{\tilde{\mathfrak{P}}(\C)}} & {} & {} & {} & {} & {} & {} & {} & {\Large\bm{\tilde{\mathfrak{C}}_{\mathsf{AR}}(\C)}} & {\color{red}{(\A,\B)}} & {} & {} \\
  {} & {} & {} & {} & {\color{red}{(\F,\F^\perp)}} & {} & {} & {} & {} & {\color{red}{({}^\perp\G,\G)}} & {} & {} & {} & {} & {} \\
  {} & {} & {} & {} & {} & {} & {} & {} & {} & {} & {} & {} & {} & {} & {} \\
  {} & {} & {} & {} & {\color{red}{\F}} & {} & {} & {} & {} & {} & {\color{blue}{\F^\perp}} & {} & {} & {} & {} \\
  {} & {} & {} & {} & {} & {\color{red}{{}^\perp\G}} & {} & {} & {} & {\color{blue}{\G}} & {} & {} & {} & {} & {} \\
  {} & {} & {\color{red}{\F}} & {\Large\bm{\tilde{\mathfrak{G}}(\C)}} & {} & {} & {} & {} & {} & {} & {} & {\Large\bm{\tilde{\mathfrak{G}}\op_{\mathsf{AR}}(\C)}} & {\color{blue}{\B}} & {} & {} \\
  {} & {} & {\color{red}{\X}} & {} & {} & {} & {} & {} & {} & {} & {} & {} & {\color{blue}{\X^\perp}} & {} & {\color{blue}{\G}} \\
  {} & {} & {} & {} & {} & {} & {} & {\Large\bm{\tilde{\mathfrak{F}}_{\mathsf{AR}}(\C)}} & {} & {} & {} & {} & {} & {} & {} \\
  {} & {} & {} & {} & {} & {} & {} & {\color{red}{(\F,\F\cap\F^\perp)}} & {} & {} & {} & {} & {} & {} & {} \\
  {} & {} & {} & {} & {} & {} & {} & {\color{red}{(\X,\omega)}} & {} & {} & {} & {} & {} & {} & {} \\
  {} & {} & {} & {} & {} & {} & {} & {\color{red}{({}^\perp\G,\G)}} & {} & {} & {} & {} & {} & {} & {} \\ };
\path[red,|->]
(m-3-3) edge node[below,sloped] {$\tilde{\Theta}$} (m-8-3)
(m-6-5) edge node[below,sloped] {$\tilde{\Omega}$} (m-4-5)
(m-7-10) edge (m-7-6)
(m-7-10) edge node[above,sloped] {$\tilde{\Omega}_{\mathsf{AR}}$} (m-4-10)
($(m-8-3.east)+(-0.2,-0.4)$) edge ($(m-11-8.west)-(-0,-0.4)$)
;
\path[red,<-|]
($(m-9-3.east)+(-0.1,-0.4)$) edge ($(m-12-8.west)-(0.4,-0.4)$)
($(m-13-8.east)+(-0.3,0.2)$) edge ($(m-9-14.west)-(-0.3,0.4)$)
;
\path[-latex]
(m-6-5) edge [-,line width=6pt,draw=white] (m-6-11)
;
\path[blue,|->]
(m-6-5) edge (m-6-11)
(m-3-13) edge (m-8-13)
($(m-12-8.east)+(0.2,0.4)$) edge ($(m-9-13.west)-(0.1,0.2)$)
;
\path[->] 
($(m-8-4.east)+(0,0.2)$) edge [thick] ($(m-8-12.west)-(0,-0.2)$)
($(m-3-4.east)+(-1,-0.4)$) edge [thick] node[left] {$\tilde{\Theta}$} ($(m-8-4.west)-(-0.5,-0.4)$)
($(m-3-12.east)+(-0.9,-0.4)$) edge [thick] node[right] {$\tilde{\Theta}_{\mathsf{AR}}$} ($(m-8-12.west)-(-1,-0.4)$)
($(m-8-4.east)+(-0.2,-0.4)$) edge [thick] ($(m-10-8.west)-(-0,-0.4)$)
($(m-10-8.east)+(-0.2,0.4)$) edge [thick] ($(m-8-12.west)-(-0.2,0.2)$)
;
\path[<-] 
($(m-8-4.east)+(0,0)$) edge [thick] ($(m-8-12.west)-(0,0)$)
($(m-3-4.east)+(-0.8,-0.4)$) edge [thick] node[right] {$\tilde{\Omega}$} ($(m-8-4.west)-(-0.7,-0.4)$)
($(m-3-12.east)+(-1.1,-0.4)$) edge [thick] node[left] {$\tilde{\Omega}_{\mathsf{AR}}$} ($(m-8-12.west)-(-0.8,-0.4)$)
($(m-8-4.east)+(-0.6,-0.4)$) edge [thick] ($(m-10-8.west)-(0.4,-0.4)$)
($(m-10-8.east)+(-0.1,0.2)$) edge [thick] ($(m-8-12.west)-(-0.3,0.4)$)
;
\path[-,font=\scriptsize]
(m-1-8) edge [double, thick, double distance=2pt] (m-3-4)
(m-1-8) edge [double, thick, double distance=2pt] (m-3-12)
(m-3-4) edge [double, thick, double distance=2pt] (m-3-12)
;
\end{tikzpicture}
\caption[Correspondences between left Frobenius pairs, left AB contexts, and perfect and hereditary cotorsion pairs, enveloping classes, and covering classes in Auslander-Reiten and Krull-Schmidt category]{Correspondences between left Frobenius pairs, left AB contexts, and perfect and hereditary cotorsion pairs, enveloping classes, and covering classes in Auslander-Reiten and Krull-Schmidt category.}
\label{fig:correspondence_ARKS}
\end{figure}

\begin{proof}[{\it {\bf Proof of Corollary~\ref{coro:ARKS}}}] 

\

First, note that we have one-to-one correspondences $\tilde{\mathfrak{G}}\op_{\mathsf{AR}}(\C) \leftrightarrow \tilde{\mathfrak{P}}_{\mathsf{AR}}(\C)$ and $\tilde{\mathfrak{P}}(\C) \leftrightarrow \tilde{\mathfrak{G}}(\C)$. So it suffices to show $\tilde{\mathfrak{P}}_{\mathsf{AR}}(\C) = \tilde{\mathfrak{P}}(\C)$. Let $(\F,\G) \in \tilde{\mathfrak{P}}(\C)$. Then $\G$ is a left saturated enveloping class in $\C$, and $\F = {}^\perp\G$ satisfies $\F^\wedge = \C$. This last equality implies $\G \subseteq \Inj^{< \infty}(\C)$.  On the other hand, since $\G$ is enveloping and contains the injective modules, it follows that every object has an injective $\G$-pre-envelope. By the \emph{Wakamatsu Lemma} \emph{\cite[Lemma 5.12]{GT}}, these $\G$-pre-envelopes can be constructed in such a way that their cokernels are in ${}^\perp\G = {}^{\perp_1}\G = \F$. Hence, $\G$ is a special pre-enveloping class, and so $(\F,\G) \in \tilde{\mathfrak{P}}_{\mathsf{AR}}(\C)$.

Now let $(\F,\G) \in \tilde{\mathfrak{P}}_{\mathsf{AR}}(\C)$. Then $\G$ is a right saturated special pre-enveloping class in $\C$, contained in $\Inj^{< \infty}(\C)$, and so $\F^\wedge = \C$. By \emph{Lemma~\ref{lem:correspondence_with_G}}, it follows that $\G$ is enveloping and $\F$ is covering in $\C$. Hence, $(\F,\G) \in \tilde{\mathfrak{P}}(\C)$. 
\end{proof}

A cotorsion pair in $\mathsf{Mod}(R)$ is said to be \textbf{cotilting} if it is of the form $({}^\perp C, ({}^\perp C)^\perp)$, for some cotilting module $C$. If $C$ is basic, then $({}^\perp C, ({}^\perp C)^\perp)$ is called a \textbf{cotilting cotorsion pair}.

The following result gives us a characterization of perfect and hereditary cotorsion pairs in $\mathsf{mod}(\Lambda)$, as a consequence of \emph{Corollary~\ref{coro:ARKS}}.

\begin{corollary}\label{coro:coro_del_AR}
Let $\Lambda$ be an Artin $R$-algebra. Then a cotorsion pair $(\F,\G)$ in $\mathsf{mod}(\Lambda)$ is perfect and hereditary if, and only if, it is basic cotilting. Moreover, for every basic cotilting module $C$, one has the equality $\add(C)^\wedge = ({}^\perp C)^\perp$. 
\end{corollary}

%%%%%%%%%%%%%%%%%%%%%%%%%%%%%%%%%%%%%%%%%%%%%%%%%%%%
%%%%%%%%%%%%%%%%%%%%%%%%%%%%%%%%%%%%%%%%%%%%%%%%%%%%

\subsection{Strong Frobenius pairs vs. Hovey triples vs. AB model structures}

We devote the last part of this section to complement the correspondences studied before, now involving the AB model structures in the interplay. We restrict out attention to the following sub-class of $\mathfrak{F}$,
\begin{align*}
\text{s}\mathfrak{F} & := \{ (\X, \omega) \subseteq \C^2 \mbox{ : $(\X,\omega)$ is a strong left Frobenius pair in $\C$ such that $\Proj(\C) \subseteq \X^\wedge$} \}
\end{align*}
and show how this class is in one-to-one correspondence with the AB model structures.

\begin{proposition}\label{prop:Frobenius_are_hereditary}
Let $\C$ be an abelian category with enough projectives. If $(\X,\omega)$ is a strong left Frobenius pair in $\C$ with $\Proj(\C) \subseteq \X^\wedge$, then the $\X^\wedge$-cotorsion pairs $(\X,\omega^\wedge)$ and $(\omega,\X^\wedge)$ in $\C$, from \emph{Theorems~\ref{CP2}} and \emph{\ref{theorem:projective_pair}}, are both left strong hereditary in $\C$.  
\end{proposition}

\begin{myproof}
By \emph{Proposition~\ref{prop:AB_Hovey_triples_hereditaries}}, we already know that $(\X,\omega^\wedge)$ and $(\omega,\X^\wedge)$ are hereditary $\X^\wedge$-cotorsion pairs in $\C$. On the one hand, by \emph{Proposition~\ref{prop:same_projectives}}, we have $\Proj(\C) = \Proj(\X^\wedge)$. On the other hand, $\Proj(\X^\wedge) \subseteq \X, \omega$. Hence, $\Proj(\C) \subseteq \X, \omega$. The fact that $\X$ and $\omega$ are resolving in $\C$ follows by \emph{Remark~\ref{rem:resolving_thick}}.
\end{myproof}

Consider the following class:
\begin{align*}
\mathfrak{T} & := \left\{ (\X,\omega) \subseteq \C^2 \mbox{ : } \begin{array}{ll} \mbox{$\omega \subseteq \X$ is closed under direct summands in $\X$, $\X^\wedge$ is an exact} \\ \mbox{sub-category of $\C$, and $(\X,\X^\wedge,\omega^\wedge)$ is a left strong hereditary} \\ \mbox{Hovey triple in $\X^\wedge$} \end{array} \right\}.
\end{align*}

In the next lines, we prove that the classes $\text{s}\mathfrak{F}$ and $\mathfrak{T}$ coincide. We begin with the following property of Hovey triples.

\begin{proposition}\label{prop:property_Hovey_triple}
Let $\C$ be an abelian category with enough projectives, and $\Sb$ be an thick sub-category of $\C$. If $(\mathcal{F,S,W})$ is a left strong hereditary Hovey triple in $\mathcal{S}$, then $(\F, \F \cap \mathcal{W})$ is a strong left Frobenius pair in $\mathcal{C}$, and $\F \cap \mathcal{W} = \Proj(\C)$. 
\end{proposition}

\begin{myproof}
Set $\omega := \F \cap \mathcal{W}$. By hypothesis, we have that $(\F, \mathcal{W})$ is a left strong hereditary $\Sb$-cotorsion pair in $\C$. By \emph{Theorem~\ref{CP3}}, we have that $(\F, \omega)$ is a left Frobenius pair in $\C$. On the other hand, by condition $\scpcinco$ for the $\Sb$-cotorsion pair $(\omega, \Sb)$, we have that for every $F \in \F$ there exists a short exact sequence 
\[
0 \to F' \to W \to F \to 0
\]
where $W \in \omega$ and $F' \in \Sb$. Using the fact that $\F$ is closed under kernel of epimorphisms in $\F$, we have $F' \in \F$. It follows that $\omega$ is a relative projective generator in $\F$. It suffices to show $\pd_{\F}(\omega) = 0$ in order to conclude that $\F$ is also $\F$-projective, and hence $(\F,\omega)$ will be a strong left Frobenius pair in $\C$. This will follow after proving $\F \cap \mathcal{W} = \Proj(\C)$. 

From the left strong hereditary $\Sb$-cotorsion pair $(\omega, \Sb)$ in $\C$, it is clear that $\omega = \Proj(\Sb)$. On the other hand, since $\C$ has enough projectives, we have by \emph{Proposition~\ref{prop:same_projectives}} that $\Proj(\Sb) = \Proj(\C)$. Therefore, the result follows.  
\end{myproof}

The previous proposition is also valid if we replace $\Sb$ by an exact sub-category $\mathcal{E} \subseteq \C$. However, it is stated and proven in terms of $\Sb$ due to the simplicity of its proof and the purposes of this paper.  

We have the following result.

\begin{theorem}\label{theo:correspondence_Hovey_pairs_vs_Frobenius_pairs}
Let $\C$ be an abelian category with enough projectives. Then $\text{s}\mathfrak{F} = \mathfrak{T}$.

Dually, if $\C$ is an abelian category with enough injectives, then the classes $\text{s}\mathfrak{F}\op$ and $\mathfrak{T}\op$ are equal, where:
\begin{align*}
\text{s}\mathfrak{F}\op & := \{ (\nu, \Y) \subseteq \C^2 \mbox{ \emph{: $(\nu,\Y)$ is a strong right Frobenius pair in $\C$ such that $\Inj(\C) \subseteq \Y^\vee$}} \}, \\
\mathfrak{T}\op & := \left\{ (\nu,\Y) \subseteq \C^2 \mbox{ \emph{:} } \begin{array}{ll} \mbox{\emph{$\nu \subseteq \Y$ is closed under direct summands in $\Y$, $\Y^\vee$ is an exact}} \\ \mbox{\emph{sub-category of $\C$, and $(\Y^\vee,\Y,\nu^\vee)$ is a right strong hereditary}} \\ \mbox{\emph{Hovey triple in $\Y^\vee$}} \end{array} \right\}.
\end{align*}
\end{theorem}

\begin{myproof}
We only prove the equality $\text{s}\mathfrak{F} = \mathfrak{T}$. Let $(\X,\omega) \in \text{s}\mathfrak{F}$, that is, $(\X,\omega)$ is a strong left Frobenius pair in $\C$ such that $\Proj(\C) \subseteq \X^\wedge$. By Proposition~\ref{prop:Frobenius_are_hereditary}, the $\X^\wedge$-cotorsion pairs $(\X,\omega^\wedge)$ and $(\omega, \X^\wedge) = (\X \cap \omega^\wedge, \X^\wedge)$ are left strong hereditary. Then the Hovey triple $(\X, \X^\wedge, \omega^\wedge)$ is left strong hereditary, and hence $(\X,\omega) \in \mathfrak{T}$. 

Now let $(\X,\omega) \in \mathfrak{T}$, that is, $(\X, \X^\wedge, \omega^\wedge)$ is a left strong hereditary Hovey triple in the exact category $\X^\wedge$ such that $\omega \subseteq \X$ is closed under direct summands in $\X$. Then by \emph{Proposition~\ref{prop:property_Hovey_triple}}, we have $(\X, \X \cap \omega^\wedge)$ is a strong left Frobenius pair in $\C$. It is only left to show that $\omega = \X \cap \omega^\wedge$. The inclusion $\omega \subseteq \X \cap \omega^\wedge$ is clear. Now suppose $X \in \X \cap \omega^\wedge$. Since $X \in \omega^\wedge$, there exists a short exact sequence
\[
0 \to W' \to W \to X \to 0
\]
with $W \in \omega$ and $W' \in \omega^\wedge$. On the other hand, $X \in \X$ and $(\X, \omega^\wedge)$ is a cotorsion pair in $\X^\wedge$, and so the previous sequence splits (as a short exact sequence in $\X^\wedge$), which implies that $X$ is a direct summand of $W \in \omega$, and so $X \in \omega$. Therefore, $\X \cap \omega^\wedge \subseteq \omega$. 
\end{myproof}

The following result is a consequence of \emph{Proposition~\ref{prop:property_Hovey_triple}} and \emph{Theorem~\ref{theo:correspondence_Hovey_pairs_vs_Frobenius_pairs}}.

\begin{corollary}\label{coro:heart_equals_proj}
Let $(\X,\omega)$ be a strong left Frobenius pair in $\C$ with enough projectives. Then the inclusion $\Proj(\C) \subseteq \X^\wedge$ implies $\omega = \Proj(\C)$. 
\end{corollary}

To conclude this section, we show that there exists a one-to-one correspondence between $\text{s}\mathfrak{F} = \mathfrak{T}$ and the following collection of exact model structures:
\begin{align*}
\mathfrak{M} & := \left\{ (\mathcal{M}, \Sb) \mbox{ : } \begin{array}{ll} \mbox{$\Sb$ is a thick sub-category of $\C$ and $\mathcal{M} = (\mathcal{Q,R,T})$ is a projective exact} \\ \mbox{model structure on $\Sb$ such that $\mathcal{Q}$ is resolving in $\C$, and $\mathcal{T} \subseteq \mathcal{Q}^\wedge$} \end{array} \right\}.
\end{align*}

\begin{theorem}\label{theo:Frobenius_pairs_vs_model_structures} 
Let $\C$ be an abelian category with enough projectives. Then the mapping
\begin{align*}
\Xi \colon \text{s}\mathfrak{F} & \to \mathfrak{M} \\
(\X, \omega) & \mapsto (\mathcal{M}^{\rm proj}_{\rm AB}(\X,\omega), \X^\wedge),
\end{align*}
where $\mathcal{M}^{\rm proj}_{\rm AB}(\X,\omega)$ is the projective AB model structure on $\X^\wedge$ from \emph{Theorem~\ref{theo:Model_category_Frobenius}}, defines a one-to-one correspondence between the classes $\text{s}\mathfrak{F}$ and $\mathfrak{M}$. 

Dually, if $\C$ is an abelian category with enough injectives, then the mapping 
\begin{align*}
\Xi\op \colon \text{s}\mathfrak{F}\op & \to \mathfrak{M}\op \\
(\nu,\Y) & \mapsto (\mathcal{M}^{\rm inj}_{\rm AB}(\nu,\Y), \Y^\vee),
\end{align*}
where $\mathcal{M}^{\rm inj}_{\rm AB}(\nu,\Y)$ is the injective AB model structure on $\Y^\vee$ from \emph{Theorem~\ref{theo:Model_category_Frobenius}}, defines a one-to-one correspondence between the classes $\text{s}\mathfrak{F}\op$ and 
\begin{align*}
\mathfrak{M}\op & := \left\{ (\mathcal{M}, \Sb) \mbox{ \emph{:} } \begin{array}{ll} \mbox{\emph{$\Sb$ is a thick sub-category of $\C$ and $\mathcal{M} = (\mathcal{Q,R,T})$ is an injective exact}} \\ \mbox{\emph{model structure on $\Sb$ such that $\mathcal{R}$ is coresolving in $\C$, and $\mathcal{T} \subseteq \mathcal{R}^\vee$}} \end{array} \right\}.
\end{align*}
\end{theorem}

\newpage

\begin{myproof}
We only prove the statement concerning $\Xi$. First, note that the map $\Xi$ is well defined since the exact model structure $\mathcal{M}^{\rm proj}_{\rm AB}(\X,\omega)$ on $\X^\wedge$ is unique by \emph{Hovey-Gillespie Correspondence}, and $\X$ is resolving in $\C$ by \emph{Theorem~\ref{theo:correspondence_Hovey_pairs_vs_Frobenius_pairs}}. 

Now we construct an inverse for $\Xi$. Let $\Gamma$ be the map:
\begin{align*}
\Gamma \colon \mathfrak{M} & \to \text{s}\mathfrak{F} \\
(\mathcal{M}, \Sb) & \mapsto (\mathcal{Q}, \mathcal{Q} \cap \mathcal{T}),
\end{align*} 
where $\mathcal{Q}$, $\mathcal{R}$ and $\mathcal{T}$ denote the classes of cofibrant, fibrant and trivial objects of $\mathcal{M}$. We check $\Gamma$ is well defined. If $\mathcal{M} = (\mathcal{Q,R,T})$ is a projective exact model structure on $\Sb$, then $\mathcal{R} = \Sb$, and by \emph{Hovey-Gillespie Correspondence} we have that $(\mathcal{Q}, \Sb, \mathcal{T})$ is a Hovey triple. On the other hand, the cotorsion pair $(\mathcal{Q} \cap \mathcal{T}, \Sb)$ in $\Sb$ is clearly left hereditary in $\Sb$, and $(\mathcal{Q,T})$ is also a left hereditary cotorsion pair in $\Sb$ since $\mathcal{Q}$ is resolving in $\Sb$. Since $\Sb$ is thick, we have $\mathcal{Q}$ and $\mathcal{Q} \cap \mathcal{T}$ are both pre-resolving in $\C$. In order to show that the Hovey triple $(\mathcal{Q}, \Sb, \mathcal{T})$ is left strong hereditary in $\Sb$ and apply \emph{Proposition~\ref{prop:property_Hovey_triple}} to conclue that $(\mathcal{Q}, \mathcal{Q} \cap \mathcal{T})$ is a strong left Frobenius pair in $\C$, it is only left to show that $\Proj(\C) \subseteq \mathcal{Q}, \mathcal{Q} \cap \mathcal{T}$. 

By definition of $\mathfrak{M}$, we have $\Proj(\C) \subseteq \mathcal{Q}$. On the other hand, note that $\mathcal{Q} \subseteq \Sb$. Then by \emph{Proposition~\ref{prop:same_projectives}}, we have that $\Proj(\C) = \Proj(\Sb)$, where $\Proj(\Sb) = \mathcal{Q} \cap \mathcal{T}$ since $(\mathcal{Q} \cap \mathcal{T}, \Sb)$ is a cotorsion pair in $\Sb$. It follows $(\mathcal{Q}, \mathcal{Q} \cap \mathcal{T}) \in \text{s}\mathfrak{F}$.

Finally, to check that $\Xi$ and $\Gamma$ are inverse to each other, we need to check the equalities $\Sb = \mathcal{Q}^\wedge$ and $\mathcal{T} = (\mathcal{Q} \cap \mathcal{T})^\wedge$, for every $(\mathcal{M} = (\mathcal{Q,R,T}), \Sb) \in \mathfrak{M}$.
\begin{itemize}[itemsep=2pt,topsep=0pt]
\item[$\bullet$] \underline{Proof of the equality $\Sb = \mathcal{Q}^\wedge$}: Note that $\mathcal{Q}^\wedge \subseteq \Sb$ since $\Sb$ is thick. Now let $S \in \Sb$. Since the pair $(\mathcal{Q,T})$ is complete in $\Sb$, there exists a short exact sequence
\[
0 \to T \to Q \to S \to 0
\]
with $Q \in \mathcal{Q}$ and $T \in \mathcal{T} \subseteq \mathcal{Q}^\wedge$. It follows $S \in \mathcal{Q}^\wedge$, and so $\Sb = \mathcal{Q}^\wedge$.  

\item[$\bullet$] \underline{Proof of the equality $\mathcal{T} = (\mathcal{Q} \cap \mathcal{T})^\wedge$}: The equality $\Sb = \mathcal{Q}^\wedge$ proven above and \emph{Theorem~\ref{CP2}} imply that 
\[
(\mathcal{Q} \cap \mathcal{T})^\wedge = \mathcal{Q}^\perp \cap \mathcal{Q}^\wedge = \mathcal{Q}^\perp \cap \Sb.
\] 
On the other hand, since $(\mathcal{Q,T})$ is a hereditary cotorsion pair in $\Sb$, we have 
\[
\mathcal{T} = \mathcal{Q}^{\perp_{1,\Sb}} = \mathcal{Q}^{\perp_1} \cap \Sb = \mathcal{Q}^\perp \cap \Sb.
\]
Hence, $(\mathcal{Q} \cap \mathcal{T})^\wedge = \mathcal{T}$.   
\end{itemize}

We have 
\begin{align*}
\Xi \circ \Gamma(\mathcal{M,S}) & = \Xi(\mathcal{Q}, \mathcal{Q} \cap \mathcal{T}) = (\mathcal{M}_{\rm AB}(\mathcal{Q}, \mathcal{Q} \cap \mathcal{T}), \mathcal{Q}^\wedge),
\end{align*}
where $(\mathcal{Q},\mathcal{Q}^\wedge, (\mathcal{Q} \cap \mathcal{T})^\wedge)$ is the Hovey triple correspinding to $\mathcal{M}_{\rm AB}(\mathcal{Q}, \mathcal{Q} \cap \mathcal{T})$. From the equalities $\mathcal{Q}^\wedge = \Sb = \mathcal{R}$ and $(\mathcal{Q} \cap \mathcal{T})^\wedge = \mathcal{T}$, we have:
\begin{align*}
\Xi \circ \Gamma(\mathcal{M,S}) & = (\mathcal{M}_{\rm AB}(\mathcal{Q}, \mathcal{Q} \cap \mathcal{T}), \mathcal{Q}^\wedge) = ((\mathcal{Q, R, T}), \Sb) = (\mathcal{M,S}).
\end{align*}
Hence $\Xi \circ \Gamma = {\rm id}_{\mathfrak{M}}$. 

Now let $(\X, \omega) \in \text{s}\mathfrak{F}$. We have:
\begin{align*}
\Gamma \circ \Xi(\X,\omega) & = \Gamma(\mathcal{M}_{\rm AB}(\X,\omega), \X^\wedge) = (\X, \X \cap \omega^\wedge).
\end{align*}
Since $(\X,\omega)$ is a left Frobenius pair, we have $\X \cap \omega^\wedge = \omega$ by \emph{Theorem~\ref{CP2}}. Then 
\begin{align*}
\Gamma \circ \Xi(\X,\omega) & = (\X, \X \cap \omega^\wedge) = (\X, \omega).
\end{align*}
Hence, $\Gamma \circ \Xi = {\rm id}_{\text{s}\mathfrak{F}}$. 
\end{myproof}

We close this section by complementing the correspondence given in \emph{Theorem~\ref{theo:correspondence_with_AB-context}} for abelian categories with enough projectives, when we restrict to the sub-class $\text{s}\mathfrak{F}$.

\begin{corollary}
Let $\C$ be an abelian category with enough projectives. Then there exists a one-to-one correspondence between the class $\text{s}\mathfrak{F}$, $\mathfrak{M}$, and:
\begin{align*}
\text{s}\mathfrak{P} & := \left\{ (\F,\G) \subseteq \C^2 \mbox{ \emph{:} } \begin{array}{ll} \mbox{\emph{$(\F,\G)$ is a $\Thick(\F)$-cotorsion pair in $\C$ with $\id_{\F}(\G) = 0$}} \\ \mbox{\emph{and $\F \cap \G = \Proj(\C)$}} \end{array} \right\}.
\end{align*}
Dually, if $\C$ is an abelian category with enough injectives, then there exists a one-to-one correspondence between the class $\text{s}\mathfrak{F}$, $\mathfrak{M}\op$, and:
\begin{align*}
\text{s}\mathfrak{P}\op & := \left\{ (\F,\G) \subseteq \C^2 \mbox{ \emph{:} } \begin{array}{ll} \mbox{\emph{$(\F,\G)$ is a $\Thick(\G)$-cotorsion pair in $\C$ with $\pd_{\G}(\F) = 0$}} \\ \mbox{\emph{and $\F \cap \G = \Inj(\C)$}} \end{array} \right\}.
\end{align*}
\end{corollary}

\begin{myproof}
We only prove that $\text{s}\mathfrak{F}$ and $\text{s}\mathfrak{P}$ are in one-to-one correspondence. Consider the map $\Phi \colon \mathfrak{F} \to \mathfrak{P}$ from \emph{Theorem~\ref{theo:correspondence_with_AB-context}}. Let $(\X,\omega) \in \text{s}\mathfrak{F} \subseteq \mathfrak{F}$. We show that $\Phi(\X,\omega) = (\X,\omega^\wedge) \in \text{s}\mathfrak{P}$. First, we already know that $(\X,\omega^\wedge)$ is a $\X^\wedge$-cotorsion pair in $\C$ with $\id_{\X}(\omega^\wedge) = 0$. On the other hand, since $\C$ has enough projectives, we can apply \emph{Corollary~\ref{coro:heart_equals_proj}} and obtain $\X \cap \omega^\wedge = \omega = \Proj(\C)$, thus proving $(\X,\omega^\wedge) \in \text{s}\mathfrak{P}$. This implies that the restriction of $\Phi$ on $\text{s}\mathfrak{F}$ gives rise to a mapping $\text{s}\Phi := \Phi|_{\text{s}\mathfrak{F}} \colon \text{s}\mathfrak{F} \to \text{s}\mathfrak{P}$. 

In order to show that the mapping $\text{s}\Phi$ defines a one-to-one correspondence between $\text{s}\mathfrak{F}$ and $\text{s}\mathfrak{P}$, it suffices to show that the restriction of $\Psi \colon \mathfrak{P} \to \mathfrak{F}$ (the inverse of $\Phi$ in \emph{Theorem~\ref{theo:correspondence_with_AB-context}}) on $\text{s}\mathfrak{P}$, has its image in $\text{s}\mathfrak{F}$. Let $(\F,\G) \in \text{s}\mathfrak{P}$, that is, $(\F,\G)$ is a $\Thick(\F)$-cotorsion pair in $\C$ with $\id_{\F}(\G) = 0$ and $\F \cap \G = \Proj(\C)$. On the one hand, we already know that $\Psi(\F,\G) = (\F, \F \cap \G)$ is a left Frobenius pair in $\C$. So it is only left to show that $\omega := \F \cap \G$ is an $\F$-projective relative generator in $\F$ with $\Proj(\C) \subseteq \F^\wedge$. This follows by the facts that $\F \cap \G = \Proj(\C)$, that $\C$ has enough projectives, and that $\F$ is left thick. 
\end{myproof}

\begin{figure}[H]
\centering 
\scriptsize
\begin{tikzpicture}[description/.style={fill=white,inner sep=2pt}]
\matrix (m) [matrix of math nodes, row sep=0.5em, column sep=3em, text height=3.5ex, text depth=1.5ex]
{ {\color{red}{(\X,\omega^\wedge)}} & {\color{red}{(\X,\omega)}} & {} & {\color{red}{(\X,\omega)}} & {\color{red}{(\mathcal{M}^{\rm proj}_{\rm AB}(\X,\omega),\X^\wedge)}} \\
  {\Large\bm{\text{s}\mathfrak{P}}} & {\Large\bm{\text{s}\mathfrak{F}}} & {} & {\Large\bm{\mathfrak{T}}} & {\Large\bm{\mathfrak{M}}} \\
  {\color{red}{(\F,\G)}} & {\color{red}{(\F,\F\cap\G)}} & {} & {\color{red}{(\mathcal{Q},\mathcal{Q}\cap\mathcal{T})}} & {\color{red}{(\mathcal{M},\Sb)}} \\ };
\path[red,|->]
(m-1-2) edge node[above] {$\text{s}\Phi$} (m-1-1)
(m-1-4) edge node[above] {$\Xi$} (m-1-5)
(m-3-1) edge node[below] {$\text{s}\Psi$} (m-3-2)
(m-3-5) edge node[below] {$\Gamma$} (m-3-4)
;
\path[->] 
($(m-2-1.east)+(0,-0.1)$) edge [thick] node[below] {$\text{s}\Psi$} ($(m-2-2.west)-(0,0.1)$)
($(m-2-4.east)+(0,0.1)$) edge [thick] node[above] {$\Xi$} ($(m-2-5.west)-(0,-0.1)$)
;
\path[<-] 
($(m-2-1.east)+(0,0.1)$) edge [thick] node[above] {$\text{s}\Phi$} ($(m-2-2.west)-(0,-0.1)$)
($(m-2-4.east)+(0,-0.1)$) edge [thick] node[below] {$\Gamma$} ($(m-2-5.west)-(0,0.1)$)
;
\path[-,font=\scriptsize]
(m-2-2) edge [double, thick, double distance=2pt] (m-2-4)
;
\end{tikzpicture}
\caption[Correspondences between strong left Frobenius pairs and projective AB model structures on abelian categories with enough projectives]{Correspondences between strong left Frobenius pairs and projective AB model structures on abelian categories with enough projectives.}
\label{fig:correspondence_AB_model_projective}
\end{figure}

\begin{figure}[H]
\centering 
\scriptsize
\begin{tikzpicture}[description/.style={fill=white,inner sep=2pt}]
\matrix (m) [matrix of math nodes, row sep=0.5em, column sep=3em, text height=3.5ex, text depth=1.5ex]
{ {\color{blue}{(\nu^\vee,\Y)}} & {\color{blue}{(\nu,\Y)}} & {} & {\color{blue}{(\nu,\Y)}} & {\color{blue}{(\mathcal{M}^{\rm inj}_{\rm AB}(\nu,\Y),\Y^\vee)}} \\
  {\Large\bm{\text{s}\mathfrak{P}\op}} & {\Large\bm{\text{s}\mathfrak{F}\op}} & {} & {\Large\bm{\mathfrak{T}\op}} & {\Large\bm{\mathfrak{M}\op}} \\
  {\color{blue}{(\F,\G)}} & {\color{blue}{(\F\cap\G,\G)}} & {} & {\color{blue}{(\mathcal{R}\cap\mathcal{T},\mathcal{R})}} & {\color{blue}{(\mathcal{M},\Sb)}} \\ };
\path[blue,|->]
(m-1-2) edge node[above] {$\text{s}\Phi\op$} (m-1-1)
(m-1-4) edge node[above] {$\Xi\op$} (m-1-5)
(m-3-1) edge node[below] {$\text{s}\Psi\op$} (m-3-2)
(m-3-5) edge node[below] {$\Gamma\op$} (m-3-4)
;
\path[->] 
($(m-2-1.east)+(0,-0.1)$) edge [thick] node[below] {$\text{s}\Psi\op$} ($(m-2-2.west)-(0,0.1)$)
($(m-2-4.east)+(0,0.1)$) edge [thick] node[above] {$\Xi\op$} ($(m-2-5.west)-(0,-0.1)$)
;
\path[<-] 
($(m-2-1.east)+(0,0.1)$) edge [thick] node[above] {$\text{s}\Phi\op$} ($(m-2-2.west)-(0,-0.1)$)
($(m-2-4.east)+(0,-0.1)$) edge [thick] node[below] {$\Gamma\op$} ($(m-2-5.west)-(0,0.1)$)
;
\path[-,font=\scriptsize]
(m-2-2) edge [double, thick, double distance=2pt] (m-2-4)
;
\end{tikzpicture}
\caption[Correspondences between strong right Frobenius pairs and injective AB model structures on abelian categories with enough injectives]{Correspondences between strong right Frobenius pairs and injective AB model structures on abelian categories with enough injectives.}
\label{fig:correspondence_AB_model_injective}
\end{figure}

%%%%%%%%%%%%%%%%%%%%%%%%%%%%%%%%%%%%%%%%%%%%%%%%%%%%%%%%%%%%%
%%%%%%%%%%%%%%%%%%%%%%%%%%%%%%%%%%%%%%%%%%%%%%%%%%%%%%%%%%%%%

\bibliographystyle{alpha}
\bibliography{bibliorcp}

\end{document}